\documentclass[11pt]{amsart}
\usepackage[T1]{fontenc}
\usepackage[utf8]{inputenc}
\usepackage{lmodern}
\usepackage{amssymb,mathtools}
\usepackage[margin=2.6cm]{geometry}
\usepackage{enumitem}
\usepackage{array}
\usepackage{longtable}
\usepackage{tikz}
\usetikzlibrary{positioning,arrows.meta}
\usepackage{microtype}
\usepackage[dvipsnames]{xcolor}
\usepackage[colorlinks=true,linkcolor=NavyBlue,citecolor=NavyBlue,urlcolor=NavyBlue]{hyperref}

\newtheorem{theorem}{Theorem}[section]
\newtheorem{proposition}[theorem]{Proposition}
\newtheorem{lemma}[theorem]{Lemma}
\newtheorem{corollary}[theorem]{Corollary}
\theoremstyle{definition}
\newtheorem{definition}[theorem]{Definition}
\theoremstyle{remark}
\newtheorem{remark}[theorem]{Remark}
\numberwithin{equation}{section}

\DeclareUrlCommand\lean{\urlstyle{tt}}
\makeatletter\g@addto@macro\UrlBreaks{\do\_}\makeatother
\newcommand{\F}{\mathcal F}
\newcommand{\U}{\mathcal U}
\newcommand{\Ev}{\mathcal E}
\newcommand{\Mlcm}{\mathcal M}
\newcommand{\Nb}{\mathcal N}
\newcommand{\lc}{\operatorname{lc}}
\newcommand{\eps}{\varepsilon}
\DeclareMathOperator{\lcm}{lcm}
\newcommand{\AImodels}{Claude Opus 5.0 and Claude Opus 5.5}
\newcommand{\leanDOI}{\href{https://doi.org/10.5281/zenodo.23188522}{10.5281/zenodo.23188522}}
\newcommand{\leanSHA}{5c49ec66c021f0923c4a20611405b2d5ef6e9e6e6f2d1ecf018abcf603f1b7bf}
\newcommand{\leanDOItwo}{\href{https://doi.org/10.5281/zenodo.24017156}{10.5281/zenodo.24017156}}
\newcommand{\leanDOIall}{\href{https://doi.org/10.5281/zenodo.23188521}{10.5281/zenodo.23188521}}

\title[Exceptional sets of $K/n=1/x+1/y+1/z$]{The exponent $3/4$ for the exceptional sets of the equations $K/n=1/x+1/y+1/z$}
\author{Dmytro Trybushnyi}
\address{Karlsruhe Institute of Technology, Karlsruhe, Germany}
\email{Dmytro.Trybushnyi@kit.edu}
\date{}

\begin{document}
\begin{abstract}
For $K\ge4$ let $E_K(N)$ be the number of integers $2\le n\le N$ for which $K/n=1/x+1/y+1/z$ has no solution in
positive integers; $K=4$ is the Erd\H{o}s--Straus equation. We prove that
\[
E_K(N)\le N\exp\Bigl(-c\Bigl(\frac{(\log N)^3}K\Bigr)^{1/4}\Bigr)
\]
with absolute constants $c,B>0$, uniformly for $4\le K\le(\log N)^3/B$. This improves the exponent $2/3$ of Vaughan's
bound of 1970. The bound is non-trivial for $K$ up to a constant times $(\log N)^3$, which is, up to a factor
$(\log\log N)^{2+o(1)}$, the scale at which, by work of Pomerance and Weingartner, exceptional primes become common. The sieve is combinatorial: a
truncated inclusion--exclusion over the classes of parametric solutions, whose second moment is the partition function
of a two-copy polymer gas, bounded by a real form of the Koteck\'y--Preiss criterion. Averaging over the residue
classes of the small primes gives the family its full mass, and summing each cluster in a random order bounds the
cluster sums. Since a divisor of an exception is an exception, nothing about the small primes themselves is needed, and
the sieve can use the classes of the whole family of parametric solutions, one triple per class. All constants are
explicit, and the whole proof has been checked in the proof assistant Lean~4.
\end{abstract}
\maketitle

\section{Introduction}\label{sec:intro}

Fix an integer $K\ge1$ and consider
\begin{equation}\label{eq:ES}
\frac Kn=\frac1x+\frac1y+\frac1z,\qquad x,y,z\in\mathbb Z_{>0}.
\end{equation}
For $K\le3$ it is solvable for every $n$, so the question starts at $K=4$. The Erd\H{o}s--Straus conjecture
\cite{Erdos1950} asserts solvability for $K=4$ and every $n\ge2$. Sierpi\'nski \cite{Sierpinski} conjectured the same for
$K=5$, and Schinzel for every fixed $K$ and all large $n$ (see \cite{Vaughan}). For $K=4$ the conjecture has been
verified up to $10^{18}$ \cite{Salez,MihneaDumitru}, and polynomial identities settle every class modulo $840$ except six
\cite{Mordell}. The problem is one of density. Call $n$ \emph{representable} if \eqref{eq:ES} is solvable, and put
\[
E_K(N)=\#\{2\le n\le N:\ n\text{ is not representable}\}.
\]

Density-zero results appeared around 1970 \cite{Webb,Vaughan}, and for more than fifty years the strongest bound
was Vaughan's \cite{Vaughan},
\begin{equation}\label{eq:Vaughan}
E_K(N)\ll_KN\exp\bigl(-c_K(\log N)^{2/3}\bigr),
\end{equation}
made uniform in $K\le(\log N)^2$ by Pomerance and Weingartner \cite[Theorem~1.3]{PW}, in the form
$E_K(N)\le N\exp(-C((\log N)^2/\varphi(K))^{1/3})$. Dahan \cite[Proposition~4.32]{Dahan} shows that the large-sieve
scheme cannot pass the exponent $2/3$ for witness systems of Vaughan's kind. In concurrent work, Dahan \cite{Dahan57}
announces $E_4(N)\ll N\exp(-c(\log N)^{5/7}(\log\log N)^{2/7})$, by a sieve with composite witnesses. Elsholtz and
Tao \cite{ElsholtzTao} count solutions rather than exceptions.

\subsection{The result}

\begin{theorem}\label{thm:A}
Let $K\ge4$ and $N\ge2$, and suppose that $h=\lfloor\log_2N\rfloor\ge4096$.
\begin{enumerate}[label=\textup{(\alph*)},leftmargin=*]
\item If $Kh\ge2^{66}$ and $2^{85}K\le h^3$, then
\[
E_K(N)\le N\exp\Bigl(-\frac{(\log N)^{3/4}}{24\,(2^{32}K)^{1/4}}\Bigr).
\]
\item If $Kh\ge2^{61}$ and $2^{90}K\le h^3$, the same holds with $2^{37}K$ in place of $2^{32}K$.
\end{enumerate}
\end{theorem}

Thus $E_K(N)\le N\exp(-c((\log N)^3/K)^{1/4})$ with $c=1/(24\cdot2^8)=1/6144$, for every $K$ with $4\le K\le h^3/2^{85}$
as soon as $h\ge\max(2^{12},2^{66}/K)$. Neither part contains the other: (a) has the better constant, and (b) the smaller
threshold for small $K$. For the Erd\H{o}s--Straus equation the conditions reduce to $h\ge2^{64}$ in (a) and $h\ge2^{59}$
in (b).

\begin{corollary}\label{cor:ES}
If $N\ge2^{2^{64}}$, then $E_4(N)\le N\exp\bigl(-(\log N)^{3/4}/8689\bigr)$. If $N\ge2^{2^{59}}$, then
$E_4(N)\le N\exp\bigl(-(\log N)^{3/4}/20666\bigr)$.
\end{corollary}

Part (b) uses only the triples with $a,c\le d$ in the parametric solutions below. Part (a) sieves with the classes of
the whole family, one triple per class, and counts the mass on a larger core (\S\ref{sec:full}). This gains a factor
$\frac{13}4$ in the mass; together with sharper numerical constants, $M$ drops by $2^5$. Since
$\varphi(K)\gg K/\log\log K$, the exponent $((\log N)^3/K)^{1/4}$ of Theorem~\ref{thm:A} exceeds the exponent
$((\log N)^2/\varphi(K))^{1/3}$ of \cite[Theorem~1.3]{PW} by a factor $\gg(K\log N)^{1/12}(\log\log K)^{-1/3}$, so the
new bound is stronger in $K$ as well as in $N$, apart from the constants.

The bound is non-trivial up to $K\asymp(\log N)^3$, and this is the natural limit. By
\cite[Theorem~3.1 and its proof]{PW}, a positive proportion of the primes in $(N/2,N]$ are exceptional once
$K\ge(\log N)^3(\log\log N)^{2+o(1)}$, so that $E_K(N)\gg N/\log N$; Theorem~\ref{thm:A} gives
$E_K(N)\le N/(\log N)^2=o(N/\log N)$ for $K\le(c/2)^4(\log N)^3/(\log\log N)^4$. In the heuristic of \cite[\S1]{PW} a prime $p$ is exceptional with probability
$\exp(-(\log p)^3/K)$; Theorem~\ref{thm:A} replaces the exponent $(\log p)^3/K$ of this heuristic, at $p\approx N$, by a
constant times its fourth root. In an independent model of the classes used here, the exponent $\frac34$ and the power
$K^{-1/4}$ are also the limit of square-sieve weights (Proposition~\ref{prop:barrier}).

\subsection{The strategy}\label{sec:plan}

Everything rests on one family of identities (Lemma~\ref{lem:typeII}): if $a,c,d\ge1$, $u=Kacd-1$ and $u\mid nc+a$,
then, with $b=(nc+a)/u$,
\[
\frac Kn=\frac1{abd}+\frac1{acdn}+\frac1{bcdn}.
\]
So for each triple $t=(a,c,d)$ the class $E_t=\{n:u_t\mid nc+a\}$, of density $1/u_t$, consists of representable
numbers, and a non-representable $n$ avoids all of them. We sieve by these classes.

\emph{The balance.} Let $\F$ be a family of triples with $u_t\le2^L$ and \emph{mass} $\mu=\sum_t1/u_t$. If the classes
were independent, a proportion $\approx e^{-\mu}$ of the integers would avoid them all. To make this rigorous on
$[1,N]$ we use a truncated inclusion--exclusion of depth $k\asymp\mu$, which costs an error $(\#\F)^{2k}=2^{O(\mu L)}$.
So we need
\[
\log N\gg\mu L,\qquad\text{and we gain}\qquad e^{-c\mu}.
\]
The triples with $u_t\le2^L$ have mass of order $L^3/K$, one logarithm for each of the free scales of $a$, $c$ and
$d$. If a mass $\mu\asymp L^3/K$ were usable, then $L^4\asymp K\log N$, and the gain is
$\exp(-c((\log N)^3/K)^{1/4})$. This is the whole source of the exponent, and also its limit (\S\ref{sec:disc}).

\emph{The obstruction.} Inclusion--exclusion over correlated classes is controlled by a cluster expansion, and its
convergence criterion needs the moduli to be $Q$-rough, that is, free of prime factors $\le Q$. Roughness keeps only the
fraction $\asymp1/\log Q$ of the mass, and the cluster sums must be bounded with $\log Q$ small. A step-by-step bound
does not allow this: it loses a power of $\log\log N$. Three ideas resolve the difficulty.

\begin{enumerate}[leftmargin=*,label=(\roman*)]
\item \emph{A real Koteck\'y--Preiss bound} (\S\ref{sec:KP}). For real activities, Dobrushin's induction gives the upper
bound we need in a few lines, without convergence of any expansion.
\item \emph{Averaging over the small primes} (\S\S\ref{sec:avg}--\ref{sec:mass}). We split $[1,N]$ into the unit
classes modulo the product $P_0$ of the small primes not dividing $K$, and condition each class $E_t$ on the small
part of its modulus. The roughness loss $1/\log Q$ is recovered, up to a constant factor, by a divisor sum of order $\log y_0$, and the
mass becomes $\asymp L^3/K$.
\item \emph{Clusters in a random order} (\S\ref{sec:free}). The weight of a cluster factors over its graph of shared
primes. When its vertices are summed in an order chosen for its shape, a vertex is cheap if it has many earlier
neighbours or many later ones. In a random order the number of earlier neighbours is uniformly distributed, so on
average every vertex is cheap.
\end{enumerate}

\begin{figure}[ht]
\centering
\small
\begin{tikzpicture}[node distance=6mm and 6mm,
  box/.style={draw,rounded corners=2pt,align=flush center,inner sep=3pt,text width=40mm},
  arr/.style={-{Latex[length=2mm]}}]
\node[box] (sieve) {Sieve weight $F_k$, transfer\\ cost $(\#\Ev+1)^{2k}$ \,(\S\ref{sec:sieve})};
\node[box,right=of sieve] (gas) {Two-copy polymer gas,\\ real Koteck\'y--Preiss (\S\ref{sec:KP})};
\node[box,right=of gas] (local) {One level: $2^{-A}$ per\\ class, $A\asymp L^3/M$ (\S\ref{sec:level})};
\node[box,below=of sieve] (key) {Key lemma, cofactor\\ exponent $j$ (\S\ref{sec:arith})};
\node[box,below=of gas] (clus) {Clusters in a random\\ order (\S\ref{sec:free})};
\node[box,below=of local] (level) {Level $L^4\asymp M\log N$,\\ $M=2^{37}K$ or $2^{32}K$ (\S\ref{sec:level})};
\node[box,below=of key] (avg) {Averaging over classes\\ mod $P_0$ (\S\ref{sec:avg})};
\node[box,below=of clus] (mass) {Mass $\ge L^3/M$\\ (\S\ref{sec:mass})};
\node[box,below=of level] (boot) {Divisors $n=sm$:\\ $\sum_s1/s=P_0/\varphi(P_0)$ (\S\ref{sec:boot})};
\draw[arr] (sieve) -- (gas); \draw[arr] (gas) -- (local);
\draw[arr] (key) -- (clus); \draw[arr] (clus) -- (gas);
\draw[arr] (avg) -- (key); \draw[arr] (avg) -- (mass);
\draw[arr] (mass.east) -- ([xshift=3mm]mass.east) |- ([yshift=-2.5mm]local.west);
\draw[arr] (level) -- (local); \draw[arr] (local.east) -- ([xshift=3mm]local.east) |- (boot.east);
\end{tikzpicture}
\caption{Structure of the proof.}\label{fig:plan}
\end{figure}
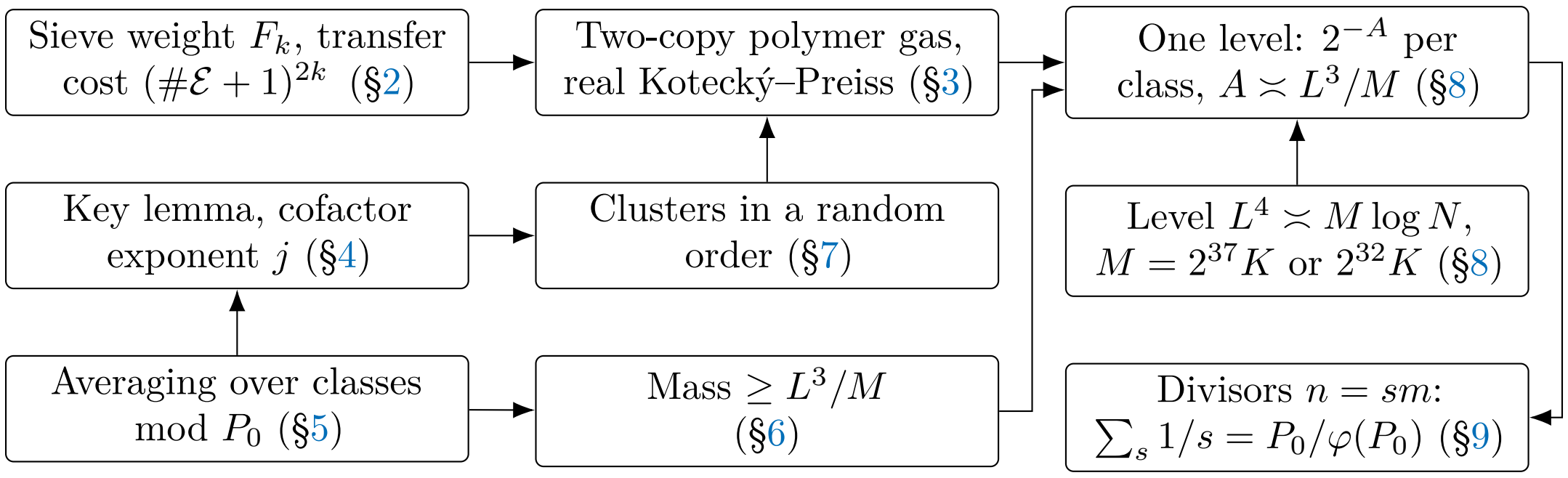

The averaging counts only the $n$ prime to $P_0$; the others are reduced to these by divisors. A representation of
$K/m$ scales to one of $K/(sm)$, so every divisor of an exceptional number is exceptional, and $E_K(N)$ is bounded by a
sum, over the $P_0$-part $s$ of $n$, of counts in unit classes (Lemma~\ref{lem:reduce}). The densities match exactly:
$\sum_s1/s=P_0/\varphi(P_0)$ against the proportion $\varphi(P_0)/P_0$ of unit classes. So nothing about which small
primes are representable is needed, and $K$ enters only through the mass. Figure~\ref{fig:plan} shows the structure of
the proof.

\subsection{Organization and conventions}\label{sec:conv}
Section~\ref{sec:sieve} sets up the sieve and the polymer gas, \S\ref{sec:KP} proves the real Koteck\'y--Preiss bound,
\S\ref{sec:arith} the key lemma on the classes $E_t$, and \S\S\ref{sec:avg}--\ref{sec:mass} construct the conditional
family and bound its mass from below. Section~\ref{sec:free} bounds the cluster sums, \S\ref{sec:level} combines
everything at one level, and \S\ref{sec:boot} proves Theorem~\ref{thm:A}(b). Section~\ref{sec:full} extends the
argument to the whole family of parametric solutions and proves Theorem~\ref{thm:A}(a). Section~\ref{sec:disc} explains the
exponent and proves that it is the limit of square-sieve weights in an independent model. Appendices~\ref{app:primes} and~\ref{app:num} contain the elementary prime
estimates and the numerics of the cluster bound, and Appendix~\ref{app:lean} the correspondence with the formal proof.

Throughout, $K\ge4$ is fixed; it is unrelated to the cutoff $K_0$ of \S\ref{sec:avg}. We write $\omega(m)$ for the
number of distinct prime factors of $m$, $\tau(m)$ for the number of divisors, $\varphi$ for Euler's function and $\#$
for cardinality; $p$ and $q$ denote primes unless stated otherwise. All logarithms are natural unless written
$\log_2$. We put $\ell=\ell(L)=\lfloor\log_2L\rfloor+1$, so that $2^{\ell-1}\le L<2^\ell$. For a set $C$ of events,
$\lc(C)$ is the least common multiple of their moduli. In \S\ref{sec:free} the letters $s$, $j$, $d$, $p$, $m$ are
reused locally for the parameters of a cluster shape, with the meanings defined there; $r$ always denotes a residue
class modulo $P_0$.

The fixed parameters are collected in Table~\ref{tab:par}; all are functions of $\ell$, except $M$, which depends on
$K$, and $T$, $A$, $k$, which depend on $K$ and $L$.

\begin{table}[ht]
\centering\small
\renewcommand{\arraystretch}{1.15}
\begin{tabular}{@{}lll@{}}
Symbol & Value & Role (where defined)\\\hline
$Q$ & $2^{128\ell}$ & roughness of the moduli $v_t$ (\S\ref{sec:avg})\\
$E_0$ & $2^{5\ell}$ & bound for the small part $s_t$ of $u_t$ (\S\ref{sec:avg})\\
$Y$ & $2^{\ell^3}$ & lower bound for $a,c$ (\S\ref{sec:avg})\\
$K_0$ & $3\ell$ & bound for $\omega(v_t)$ (\S\ref{sec:avg})\\
$y_0$, $P_0$ & $2^\ell$, $\prod_{p\le y_0,\,p\nmid K}p$ & small primes and their product (\S\ref{sec:avg})\\
$r_B$, $J_0^\flat$ & $2\lfloor\ell/2\rfloor+43$, $30+\ell+(128\ell+1)r_B$ & Brun's sieve: depth and length (\S\ref{sec:brun})\\
$G_0$, $H_0$, $J_1$ & $16\ell$, $3\ell+4$, $J_0^\flat+5\ell+1$ & margins of the scales (Proposition~\ref{prop:massF0})\\
$\Gamma$ & $\ell^3-2\ell-7$ & cap of the step costs (\S\ref{sec:num}, Appendix~\ref{app:num})\\
$T_0$, $H$ & $J_0^\flat+10\ell+64$, $3\ell+4$ & margins of the wide scales (\S\ref{sec:full})\\
$M$ & $2^{37}K$ or $2^{32}K$ & mass $\ge L^3/M$ (\S\S\ref{sec:mass}, \ref{sec:full})\\
$T$, $A$, $k$ & $\lfloor L^3/M\rfloor$, $\lfloor\frac23\max(T-1,0)\rfloor$, $15T$ & mass used, gain, depth (Proposition~\ref{prop:oneclass})
\end{tabular}
\caption{The parameters.}\label{tab:par}
\end{table}

\subsection{Formal verification}\label{sec:formalintro}
The proof in this paper is complete and does not depend on a computer. Independently, we formalized every lemma,
proposition, corollary and theorem of this paper in Lean~4 \cite{Lean4} with Mathlib \cite{mathlib}, with all explicit
constants; in a few places the formal statement differs slightly from the text, in ways listed in
Appendix~\ref{app:lean} that change no conclusion. The kernel checks every proof, with no axioms beyond \texttt{propext}, \texttt{Classical.choice} and
\texttt{Quot.sound}, and the numerical inequalities are proved in the kernel from rational bounds on $\log2$ and $e$,
without floating-point evaluation. The headline statements, Theorem~\ref{thm:A}, Corollary~\ref{cor:ES} and
Proposition~\ref{prop:barrier}, are restated in a separate audit file from Mathlib definitions only, so that a reader can
check that the formal statements say what the paper claims. Appendix~\ref{app:lean} describes the package, states
where it is deposited, and maps each statement to its declarations.

\subsection{AI Declaration}\label{sec:AI}
The author used Claude (Anthropic; \AImodels) for searching the literature, code generation (the Lean~4 formalization and
the programs for the numerical experiments), testing hypotheses, ruling out other approaches, devising proof and
computational strategies, drafting proofs of individual lemmas, checking the results, and drafting, revising and
proofreading the text. The author posed the problem, directed the work, reviewed the arguments, and takes full
responsibility for the content. The final results were verified without the use of AI: every step of the proofs of
Theorem~\ref{thm:A}, with all constants, and of Proposition~\ref{prop:barrier} is checked by Lean's kernel, and the
author checked that the statements of the audit file, which use only Mathlib definitions, are those of
Theorem~\ref{thm:A}, Corollary~\ref{cor:ES} and Proposition~\ref{prop:barrier} (\S\ref{sec:formalintro}).

\section{Events, the sieve weight and the two-copy gas}\label{sec:sieve}

\begin{lemma}[Type II solutions]\label{lem:typeII}
Let $t=(a,c,d)$ with $a,c,d\ge1$ and $u=Kacd-1$. Suppose $n\ge1$ and $u\mid nc+a$. Then $n$ is representable.
\end{lemma}

\begin{proof}
Put $b=(nc+a)/u$, a positive integer. Then $b(Kacd-1)=nc+a$, that is $Kabcd=a+b+cn$, and dividing by
$abcdn$ gives
\[
\frac Kn=\frac1{abd}+\frac1{acdn}+\frac1{bcdn}.\qedhere
\]
\end{proof}

\emph{Events.} An \emph{event} is a triple $\eps=\langle m,c,a\rangle$ of integers with $m\ge1$ and
$\gcd(m,c)=1$. It \emph{holds} at $j\in\mathbb Z$ if $m\mid jc+a$, which is a single residue class of $j$ modulo
$m$; $m=m_\eps$ is its \emph{modulus}. For a triple $t=(a,c,d)$, $E_t=\langle u_t,c,a\rangle$ is an event,
since $u_t=Kacd-1$ is prime to $c$, and by Lemma~\ref{lem:typeII} a non-representable $n\ge1$ lies in no $E_t$.
If $R$ is a finite set of events, then by the Chinese remainder theorem the set of $j$ at which every event of
$R$ holds is either empty or a single class modulo $\lc(R)=\lcm_{\eps\in R}m_\eps$.

\emph{The sieve weight.} Let $\Ev$ be a finite family of events, and put
$S(j)=\{\eps\in\Ev:\ \eps\text{ holds at }j\}$. A set $T\subseteq\Ev$ is \emph{independent} if its moduli are
pairwise coprime. Let $I_i(j)$ be the number of independent $i$-subsets of $S(j)$, and
\[
F_k(j)=\sum_{i=0}^k(-1)^iI_i(j).
\]
If no event of $\Ev$ holds at $j$, then $S(j)=\emptyset$ and $F_k(j)=I_0(j)=1$. The square $F_k^2$ is a Selberg-type \cite{HalberstamRichert}
weight: it is nonnegative, equal to $1$ where no event holds, and its average is an alternating sum over pairs of
independent sets.

\begin{lemma}[Transfer]\label{lem:transfer}
Let $\Mlcm=\lc(\Ev)$. Then for every $N\ge0$,
\[
\Mlcm\sum_{0\le j\le N}F_k(j)^2\le N\sum_{0\le j<\Mlcm}F_k(j)^2+\Mlcm(\#\Ev+1)^{2k}.
\]
\end{lemma}

\begin{proof}
Expanding the square,
\[
F_k(j)^2=\sum_{T,T'}(-1)^{|T|+|T'|}\,\mathbf 1[\text{every }\eps\in T\cup T'\text{ holds at }j],
\]
where $T,T'$ run over independent subsets of $\Ev$ of size $\le k$. The set of $j$ at which all events of
$R=T\cup T'$ hold is empty or a single class modulo $\lc(R)$, and $\lc(R)\mid\Mlcm$. For a single class, the count in
$[0,N]$ differs from $N/\lc(R)$ by at most $1$, and the count in $[0,\Mlcm)$ is exactly $\Mlcm/\lc(R)$. Multiply by $\Mlcm$ and
sum. The number of pairs $(T,T')$ is at most $(\sum_{i\le k}\binom{\#\Ev}{i})^2\le(\#\Ev+1)^{2k}$.
\end{proof}

Expanding the square once more gives the \emph{square expansion}
\[
\sum_{j<\Mlcm}F_k(j)^2=\sum_{i,i'\le k}(-1)^{i+i'}C_{ii'},\qquad
C_{ii'}=\sum_{j<\Mlcm}I_i(j)I_{i'}(j).
\]

\emph{The gas.} We now encode the square expansion as a polymer gas, so that the tail beyond degree $k$ is controlled
by the same gas at positive activities. In the rest of this section every modulus $m_\eps$, $\eps\in\Ev$, is assumed
to be at least $2$; for the families used below this holds by (E2) of \S\ref{sec:avg}. Let $\U=\{1,2\}\times\Ev$. Two
vertices $v,w\in\U$ are \emph{adjacent}, $v\sim w$, if $\gcd(m_v,m_w)\neq1$; this includes $(1,\eps)\sim(2,\eps)$. For
$S\subseteq\U$ write $S=(\{1\}\times T)\cup(\{2\}\times T')$ and
\begin{multline*}
W_{y,z}(S)=\mathbf 1[T,T'\text{ independent}]\ \mathbf 1[\text{every component of }S\text{ has}\le2k\text{ vertices}]\\
\times\ y^{|T|}z^{|T'|}\ \frac{\#\{0\le j<\Mlcm:\ \text{every }\eps\in T\cup T'\text{ holds at }j\}}{\Mlcm}.
\end{multline*}
Put $Z(y,z)=\sum_{S\subseteq\U}W_{y,z}(S)$. Two sets are \emph{separated} if no vertex of one is
adjacent to a vertex of the other.

\begin{lemma}\label{lem:gas}
\leavevmode
\begin{enumerate}[label=\textup{(\alph*)}]
\item If $S_1,S_2\subseteq\U$ are separated, then $W(S_1\cup S_2)=W(S_1)W(S_2)$.
\item Consequently $Z(y,z)=\sum_{\mathcal G}\prod_{C\in\mathcal G}W_{y,z}(C)$. The sum runs over families
$\mathcal G$ of pairwise separated connected subsets of $\U$ (\emph{polymers}). Thus $Z$ is the partition
function of a hard-core polymer gas in which two polymers are incompatible unless separated. Incompatibility is
symmetric, and it is reflexive because every vertex is adjacent to itself, as $\gcd(m_v,m_v)=m_v\ge2$.
\item Write $\Mlcm Z(y,z)=\sum_{i,i'}c_{ii'}y^iz^{i'}$. Then every $c_{ii'}\ge0$, and $c_{ii'}=C_{ii'}$ for
$i,i'\le k$.
\end{enumerate}
\end{lemma}

\begin{proof}
(a) The moduli of $S_1$ and of $S_2$ are coprime across, so by the Chinese remainder theorem the density of
common solutions is multiplicative. Independence within a copy holds for $S_1\cup S_2$ if and only if it holds for
$S_1$ and $S_2$. The components of $S_1\cup S_2$ are those of $S_1$ and of $S_2$.

(b) Decompose $S$ into its connected components, which are pairwise separated, and apply (a).

(c) $c_{ii'}$ sums nonnegative counts. If $|T|=i\le k$ and $|T'|=i'\le k$, every component of $S$ has at
most $2k$ vertices, so the size condition is automatic and the sum is $C_{ii'}$.
\end{proof}

\begin{lemma}[Truncation by positivity]\label{lem:tail}
With $S_k=\sum_{i,i'\le k}(-1)^{i+i'}c_{ii'}$,
\[
2^{k+1}\bigl(S_k-\Mlcm Z(-1,-1)\bigr)\le\Mlcm Z(2,2).
\]
\end{lemma}

\begin{proof}
$S_k-\Mlcm Z(-1,-1)$ is a sum of terms $\pm c_{ii'}$ with $\max(i,i')>k$. Each satisfies
$\pm c_{ii'}\le c_{ii'}\le2^{i+i'-k-1}c_{ii'}$, since $i+i'\ge k+1$ and $c_{ii'}\ge0$.
\end{proof}

\begin{proposition}\label{prop:cond}
Let $A,B\ge0$ be integers with $2^AZ(-1,-1)\le1$, $Z(2,2)\le2^B$ and $A+B\le k+1$. Then for every $N'\ge0$,
\[
2^A\,\#\{0\le j\le N':\ \text{no event of }\Ev\text{ holds at }j\}\le2N'+2^A(\#\Ev+1)^{2k}.
\]
\end{proposition}

\begin{proof}
By Lemma~\ref{lem:gas}(c) and Lemma~\ref{lem:tail},
\[
2^A\sum_{j<\Mlcm}F_k(j)^2=2^AS_k\le 2^A\Mlcm Z(-1,-1)+2^{A-k-1}\Mlcm Z(2,2)\le\Mlcm+\Mlcm=2\Mlcm.
\]
Insert this into Lemma~\ref{lem:transfer} with $N=N'$, and use $F_k(j)=1$ at every $j$ where no event holds.
\end{proof}

\section{A real Koteck\'y--Preiss bound}\label{sec:KP}

Let $\mathcal P$ be a finite set of polymers with a reflexive, symmetric incompatibility relation
$\not\sim$, and let $w:\mathcal P\to\mathbb R$ be activities of arbitrary sign. For
$\Lambda\subseteq\mathcal P$ put
\[
Z(\Lambda)=\sum_{S\subseteq\Lambda\ \text{compatible}}\ \prod_{\gamma\in S}w(\gamma),
\]
and let $\Nb(\gamma)=\{\gamma'\ne\gamma:\gamma'\not\sim\gamma\}$. Splitting according to whether
$\gamma\in S$ gives, for $\gamma\notin\Lambda$,
\begin{equation}\label{eq:rec}
Z(\Lambda\cup\{\gamma\})=Z(\Lambda)+w(\gamma)\,Z(\Lambda\setminus\Nb(\gamma)).
\end{equation}

\begin{proposition}[Dobrushin's criterion]\label{prop:dob}
Let $c:\mathcal P\to[0,\infty)$ and $S(\gamma)=\sum_{\gamma'\in\Nb(\gamma)}c(\gamma')$. Suppose
\[
|w(\gamma)|\,e^{S(\gamma)}\le1-e^{-c(\gamma)}\qquad\text{for all }\gamma.
\]
Then $Z(\Lambda)>0$ for all $\Lambda$, and
\[
Z(\mathcal P)\le\exp\sum_\gamma\bigl(w(\gamma)+|w(\gamma)|(e^{S(\gamma)}-1)\bigr).
\]
\end{proposition}

\begin{proof}
We prove by strong induction on $\Lambda$ that
\begin{equation}\label{eq:Q}
Z(\Lambda)>0\quad\text{and}\quad|Z(\Lambda)-Z(\Lambda\setminus\gamma)|\le(1-e^{-c(\gamma)})\,Z(\Lambda\setminus\gamma)\quad(\gamma\in\Lambda).
\end{equation}

\emph{Consequences of \eqref{eq:Q}.} It gives
$e^{-c(\gamma)}Z(\Lambda\setminus\gamma)\le Z(\Lambda)\le(2-e^{-c(\gamma)})Z(\Lambda\setminus\gamma)$.
Since $2-e^{-x}\le e^{x}$, removing one polymer changes $Z$ by a factor in $[e^{-c},e^{c}]$. If
\eqref{eq:Q} holds for all subsets of $\Lambda'$, then iterating shows that removing a set $D$ from $\Lambda'$
changes $Z$ by a factor in $[e^{-\sum_Dc},e^{\sum_Dc}]$.

\emph{Inductive step.} Let $\gamma\in\Lambda$ and $\Lambda'=\Lambda\setminus\gamma$. By \eqref{eq:rec},
\[
|Z(\Lambda)-Z(\Lambda')|=|w(\gamma)|\,Z(\Lambda'\setminus\Nb(\gamma))\le|w(\gamma)|e^{S(\gamma)}Z(\Lambda')\le(1-e^{-c(\gamma)})Z(\Lambda').
\]
Positivity then follows from $Z(\Lambda)\ge e^{-c(\gamma)}Z(\Lambda')>0$, with $Z(\emptyset)=1$.

\emph{The bound.} Insert polymers one at a time. By \eqref{eq:rec},
$Z(\Lambda\cup\gamma)=Z(\Lambda)(1+w(\gamma)\rho)$ with $\rho\in[e^{-S(\gamma)},e^{S(\gamma)}]$. For either
sign of $w$ we have $w\rho\le w+|w|(e^S-1)$; for $w<0$ this uses $e^{-S}\ge2-e^S$. Finally
$1+x\le e^x$.
\end{proof}

\begin{corollary}[Koteck\'y--Preiss, real form]\label{cor:KP}
Let $a:\mathcal P\to[0,\infty)$ satisfy
\[
\sum_{\gamma'\not\sim\gamma}|w(\gamma')|e^{a(\gamma')}\le a(\gamma)\qquad\text{for all }\gamma,
\]
where the sum includes $\gamma'=\gamma$. Then
$Z(\mathcal P)\le\exp\bigl(\sum_\gamma w(\gamma)+\sum_\gamma|w(\gamma)|(e^{a(\gamma)}-1)\bigr)$.
\end{corollary}

\begin{proof}
Take $c=|w|e^{a}$. The hypothesis gives $S(\gamma)\le a(\gamma)-c(\gamma)$. Hence
\[
|w|e^{S}\le|w|e^{a-c}=ce^{-c}\le1-e^{-c},
\]
using $c\le e^c-1$. Also $e^S-1\le e^a-1$. Apply Proposition~\ref{prop:dob}.
\end{proof}

Koteck\'y and Preiss \cite{KP} prove convergence of the full cluster expansion, for complex activities, under the
hypothesis of Corollary~\ref{cor:KP}. For real activities the bound of the corollary is all we need, and its proof
follows Dobrushin \cite{Dobrushin}; see also \cite{FernandezProcacci}. The induction is the one behind the Lov\'asz
local lemma; for the hard-core gas with negative activities, the local lemma in its Shearer and cluster-expansion forms
is treated in \cite{Shearer,ScottSokal,Bissacot}. Applied to the gas of \S\ref{sec:sieve}, the
first-order term $\sum W_{-1,-1}$ is $-2\mu+\mu=-\mu$, from the singletons and the double sites
(\S\ref{sec:free}), so the whole problem is to show that the larger polymers are small.

\section{The key lemma}\label{sec:arith}

In this section $L\ge1$ and $Y\ge1$ are arbitrary integers. A triple $t=(a,c,d)$ is \emph{admissible} if
$Y\le\min(a,c)$, $\max(a,c)\le d$, $\gcd(a,c)=1$ and $u_t\le2^L$. Put
\[
c_Y=\frac{L+1}{4Y^2}+\frac{224+8(L+1)^2}{Y}.
\]
A triple is determined by its modulus and its class. Indeed, let $(a,c,d)$ and $(a',c',d')$ have $\gcd(a,c)=\gcd(a',c')=1$,
$a,c\le d$, $a',c'\le d'$ and the same modulus $u=Kacd-1=Ka'c'd'-1$, and suppose their classes meet, so that
$u\mid ac'-a'c$. From $a^2\le acd$ and the like we get $|ac'-a'c|<2acd<u$, so $ac'=a'c$, and coprimality gives
$(a,c,d)=(a',c',d')$. The proof uses the quantitative versions of this rigidity: the key lemma below, and
Lemma~\ref{lem:nodup} for the conditional family. The key lemma drives all the cluster bounds of \S\ref{sec:free}.

\begin{lemma}\label{lem:lattice}
Let $D\ge1$, $\rho\in\mathbb Z$, and $\Lambda=\{(x,y)\in\mathbb Z^2: D\mid x+\rho y\}$. Let $w$ be a
nonzero vector of $\Lambda$ with minimal sup-norm $\lambda_1=|w|_\infty$. Then:
\begin{enumerate}[label=\textup{(\alph*)}]
\item $\lambda_1^2\le D$;
\item every $p\in\Lambda$ with $\det(p,w)=0$ is an integer multiple of $w$;
\item every $p\in\Lambda$ with $\det(p,w)\ne0$ satisfies $D\le2|p|_\infty\lambda_1$;
\item for $A,B\ge1$, $\#\bigl(\Lambda\cap[A,2A)\times[B,2B)\bigr)\le2AB/D+3\max(A,B)/\lambda_1+1$.
\end{enumerate}
\end{lemma}

\begin{proof}
(a) Among the $(\lfloor\sqrt D\rfloor+1)^2>D$ pairs $(x,y)\in[0,\sqrt D]^2$, two have equal
$x+\rho y$ modulo $D$. Their difference is a nonzero point of $\Lambda$ with sup-norm $\le\sqrt D$.

(b) Write $w=g w_0$ with $w_0$ primitive. A point $p$ with $\det(p,w)=0$ is a multiple $s w_0$. If $g\nmid s$, choose an
integer $m$ with $|s-mg|<g$. Then $p-mw=(s-mg)w_0$ is a nonzero point of $\Lambda$
shorter than $w$, which contradicts minimality. So $g\mid s$.

(c) $\Lambda$ has index $D$ in $\mathbb Z^2$, so $\det(p,w)\in D\mathbb Z$. Also
$|\det(p,w)|\le2|p|_\infty|w|_\infty$.

(d) Group the points by the value $\det(\cdot,w)$. Each fibre is a translate of a subset of
$\mathbb Zw$, so it has at most $\min(A/|w_x|,B/|w_y|)+1$ points. The values of $\det(\cdot,w)$ on the
box are multiples of $D$ in an interval of length at most $A|w_y|+B|w_x|$. Multiplying the number of
values by the fibre bound gives $2AB/D$, plus $\min(A/|w_x|,B/|w_y|)\le\max(A,B)/\lambda_1$, plus
$(A|w_y|+B|w_x|)/D\le2\max(A,B)\lambda_1/D\le2\max(A,B)/\lambda_1$ by (a), plus $1$.
\end{proof}

\begin{lemma}[Key lemma]\label{lem:key}
Let $Q\ge1$, $j\ge1$, $D\ge2$ and $\rho\in\mathbb Z$. Then
\[
\sum_{t}\frac D{u_t}\le\frac2{(Q+1)^j}+c_Y+\frac{16(L+1)}{\sqrt D}+\frac{(L+1)^3}{2D},
\]
the sum running over admissible $t$ with $D\mid u_t$, $D\mid a+\rho c$ and $u_t\ge(Q+1)^jD$.
\end{lemma}

The idea is as follows. The pairs $(a,c)$ lie in the lattice $\{D\mid x+\rho y\}$ of index $D$. For fixed $(a,c)$ the
moduli form one progression modulo $KacD$, so only the smallest one matters. On the ray of the shortest lattice vector
there is at most one coprime pair, and its smallest modulus is at least $(Q+1)^jD$; this gives the \emph{ray term}
$2(Q+1)^{-j}$, the reciprocal of the guaranteed cofactor $u_t/D$. Off the ray the smallest modulus is at least
$2ac\max(a,c)$, and a dyadic count of lattice points gives the rest.

\begin{proof}
Put $Q'=(Q+1)^j-1\ge1$, so $u_t\ge(Q'+1)D$ and in particular $u_t\ne D$. The pairs $(a,c)$ of the sum lie in the
lattice $\Lambda$ of Lemma~\ref{lem:lattice}. Let $w$ and $\lambda_1$ be as there.

\emph{Fibres.} Fix $(a,c)$. Then $D\mid Kacd-1$ fixes $d$ modulo $D$, so the moduli $u$ lie in one class
modulo $KacD$. They lie in $[u_0,2^L]$, where $u_0$ is the smallest. Hence, as $K\ge4$,
\[
\sum D/u\le D/u_0+(L+1)/(Kac)\le D/u_0+(L+1)/(4ac),
\]
by $\sum_{m\le2^L}1/(u_0+gm)\le1/u_0+(L+1)/g$.

\emph{The ray.} At most one coprime positive pair lies on $\mathbb Rw$: by Lemma~\ref{lem:lattice}(b) such a pair is
$sw$ with $s\in\mathbb Z$, and primitivity forces $s=\pm1$. For it
$D/u_0\le1/Q'\le2/(Q+1)^j$, so it contributes at most $2/(Q+1)^j+(L+1)/(4Y^2)$.

\emph{Off the ray.} For an off-ray pair the smallest modulus satisfies $u_0\ge2ac\max(a,c)$ (as $\max(a,c)\le d$), so
by the fibre bound the pair contributes at most $D/(2ac\max(a,c))+(L+1)/(4ac)$. Split the pairs into dyadic boxes
$[A,2A)\times[B,2B)$ with $A,B\in\{2^i:0\le i\le L\}$, and write $M_1=\max(A,B)$, $m_1=\min(A,B)$. A pair in the box
contributes at most $D/(2ABM_1)+(L+1)/(4AB)$. A box containing an off-ray point satisfies $M_1>X:=D/(4\lambda_1)$, since
$D\le2|p|_\infty\lambda_1<4M_1\lambda_1$ by Lemma~\ref{lem:lattice}(c), and $m_1>Y/2$, since $a,c\ge Y$. By
Lemma~\ref{lem:lattice}(d) the box has at most $2AB/D+3M_1/\lambda_1+1$ points, so it contributes at most
\[
\Bigl(\frac{2AB}D+\frac{3M_1}{\lambda_1}+1\Bigr)\Bigl(\frac{D}{2ABM_1}+\frac{L+1}{4AB}\Bigr)
=\frac1{M_1}+\frac{L+1}{2D}+\frac{3D}{2\lambda_1AB}+\frac{3(L+1)}{4\lambda_1m_1}+\frac D{2ABM_1}+\frac{L+1}{4AB}.
\]
Summing the six terms over the boxes with $M_1>X$ and $m_1>Y/2$ (at most $2(L+1)$ boxes have a given value of $M_1$ or of
$m_1$, and there are at most $(L+1)^2$ boxes), with $\sum_{2^i\ge X}2^{-i}\le2/X$, $\sum_{2^i>Y/2}2^{-i}\le4/Y$ and
$\sum_{2^i\ge X}4^{-i}\le\frac43X^{-2}$, the off-ray points contribute at most the following; the fourth and sixth terms
of the box bound give $6(L+1)^2/Y$ and $2(L+1)^2/Y$, which are combined in the last term:
\[
\frac{4(L+1)}X+\frac{(L+1)^3}{2D}+\frac{24D}{\lambda_1XY}+\frac{8D}{X^2Y}+\frac{8(L+1)^2}Y.
\]
With $X=D/(4\lambda_1)$ and $1\le\lambda_1\le\sqrt D$ (Lemma~\ref{lem:lattice}(a)): $4(L+1)/X\le16(L+1)/\sqrt D$,
$24D/(\lambda_1XY)=96/Y$ and $8D/(X^2Y)=128\lambda_1^2/(DY)\le128/Y$. Together with the ray this gives the bound, with
$c_Y=\frac{L+1}{4Y^2}+\frac{96+128+8(L+1)^2}Y$.
\end{proof}

If every $u_t$ in the sum has no divisor in $[2,Q]$, the condition $u_t\ge(Q+1)D$ follows from $u_t\ne D$, since the
cofactor $u_t/D\ne1$ divides $u_t$; so the case $j=1$ applies. For $j=1$ the proof gives the ray term $1/Q'=1/Q$, slightly
better than $2/(Q+1)$; the formal proof uses this form. The exponent $j$ is what \S\ref{sec:free} feeds in: a
vertex of a cluster with $j$ later neighbours or private primes has cofactor at least $(Q+1)^j$.

\section{Averaging over the small primes}\label{sec:avg}

From now on $L\ge2^{23}$, so that $\ell=\ell(L)\ge24$, and $K\le2^{3\ell}$. The parameters are
\[
Q=2^{128\ell},\qquad E_0=2^{5\ell},\qquad Y=2^{\ell^3},\qquad K_0=3\ell,\qquad y_0=2^\ell,\qquad
P_0=\prod_{p\le y_0,\ p\nmid K}p .
\]
Every prime factor of $P_0$ is at most $y_0\le Q$ and does not divide $K$. Nothing is assumed about its
representability. For a triple $t$ let
\[
s_t=\gcd(u_t,P_0),\qquad v_t=u_t/s_t .
\]
Call a positive integer \emph{good} if it is squarefree, has no prime factor $\le Q$, and has at most $K_0$ prime
factors.

\begin{definition}\label{def:condfam}
For $r\in\mathbb Z$, the \emph{conditional family} $\F_r$ is the set of triples $t=(a,c,d)$ with
\begin{itemize}[leftmargin=*]
\item $Y\le\min(a,c)$, $\max(a,c)\le d$, $\gcd(a,c)=1$ and $u_t\le2^L$;
\item $s_t\le E_0$ and $s_t\mid a+rc$;
\item $v_t$ has no divisor in $[2,Q]$.
\end{itemize}
For $\F\subseteq\F_r$ its \emph{mass} is $\mu(\F)=\sum_{t\in\F}1/v_t$.
\end{definition}

So the $Q$-smooth part of $u_t$ is $s_t$: a squarefree divisor of $P_0$ of size at most $E_0$, compatible with
the class $r$. Since $v_t$ is $Q$-rough, $\gcd(v_t,P_0)=1$. By Lemmas~\ref{lem:condev} and~\ref{lem:nodup} below, the mass is
$\sum1/(\text{modulus})$ for the events of $\F$ on the progression $r+P_0j$, not $\sum1/u_t$: the small part of the
modulus is paid once, as a congruence on $(a,c)$, and then disappears.

\begin{lemma}[Divisor reduction]\label{lem:reduce}
For $X\ge0$ put
\[
R(X)=\sum_{\substack{0\le r<P_0\\ \gcd(r,P_0)=1}}\#\bigl\{0\le j\le\lfloor X/P_0\rfloor:\ r+P_0j\ \text{is not representable}\bigr\}.
\]
Then
\[
E_K(N)\le\sum_{\substack{1\le s\le N\\ s\mid P_0^\infty}}R\bigl(\lfloor N/s\rfloor\bigr),
\]
where $s\mid P_0^\infty$ means that every prime factor of $s$ divides $P_0$.
\end{lemma}

\begin{proof}
If $n=sm$ and $K/m=1/x+1/y+1/z$, then $K/n=1/(sx)+1/(sy)+1/(sz)$. So every divisor $m$ of a non-representable $n$ is
non-representable. This includes $m=1$, since $1/x+1/y+1/z\le3<K$. Let $2\le n\le N$ be non-representable, $s$ the
largest divisor of $n$ with $s\mid P_0^\infty$, and $m=n/s$. Then $\gcd(m,P_0)=1$, $m\le\lfloor N/s\rfloor$ and $m$ is
not representable. Write $m=r+P_0j$ with $r=m\bmod P_0$, so $\gcd(r,P_0)=1$ and $j\le\lfloor\lfloor N/s\rfloor/P_0\rfloor$.
The map $n\mapsto(s,r,j)$ is injective, as $n=s(r+P_0j)$.
\end{proof}

If every prime of $P_0$ were representable only $s=1$ would occur. In general the $s$ are paid in \S\ref{sec:boot},
where their reciprocals add up to $P_0/\varphi(P_0)$ and cancel the density $\varphi(P_0)/P_0$ of the unit classes.

\begin{lemma}\label{lem:condev}
Let $t\in\F_r$, $r\ge0$, and $\hat e_t=\langle v_t,\,P_0c_t,\,rc_t+a_t\rangle$. Then $\hat e_t$ is an event, and
for every $j\ge0$, $r+P_0j\in E_t$ if and only if $\hat e_t$ holds at $j$.
\end{lemma}

\begin{proof}
$(r+P_0j)c+a=jP_0c+(rc+a)$. Here $u_t=s_tv_t$ with $\gcd(s_t,v_t)=1$, and $s_t$ divides both $P_0$ and $rc+a$, so it
divides the right side for every $j$; hence $u_t$ divides it if and only if $v_t$ does. Finally
$\gcd(v_t,c)=1$ because $v_t\mid Kacd-1$, and $\gcd(v_t,P_0)=1$.
\end{proof}

\emph{Why conditioning loses nothing.} For a prime $q\mid P_0$ and a unit class $r$, the weight
$q\,\mathbf 1[q\,\|\,u]\,\mathbf 1[q\mid a+rc]+\mathbf 1[q\nmid u]$ has average exactly $1$ over the triples $(a,c,d)$
modulo $q^2$ with $(a,c)\not\equiv(0,0)$, by a direct count, while the roughness condition $q\nmid u$ alone keeps only the
fraction $1-1/q$. At a prime $q\mid K$ the conditioning is impossible, but $q$ never divides a modulus, and roughness
costs nothing there either. The proof uses a simpler lower bound: for each $e\mid P_0$ with $e\le E_0$ it keeps the
triples with $s_t=e$ and $v_t$ rough. The resulting weights are $\varphi(e)/e^2$, because for $p\mid e$ and
$e\mid a+rc$, $p\mid c$ forces $p\mid a$ (\S\ref{sec:mass}).

\subsection{The conditional key lemma}
For $D\ge1$ and $j\ge1$ put
\[
\eps'_D(j)=E_0\Bigl(\frac2{(Q+1)^j}+c_Y\Bigr)+32(L+1)\sqrt{\frac{E_0}D}+\frac{(L+1)^3(1+\log E_0)}{2D},\qquad\eps'_D=\eps'_D(1).
\]

\begin{lemma}[Conditional key lemma]\label{lem:condkey}
Let $r,\rho\in\mathbb Z$, $j\ge1$ and $D\ge2$ with $\gcd(D,P_0)=1$. Summing over $t\in\F_r$ with $D\mid v_t$,
$D\mid a+\rho c$ and $v_t\ge(Q+1)^jD$,
\[
\sum_t\frac D{v_t}\le\eps'_D(j).
\]
For $j=1$ the condition $v_t\ge(Q+1)D$ may be replaced by $v_t\ne D$.
\end{lemma}

\begin{proof}
Group the triples by $e=s_t$, a divisor of $P_0$ with $e\le E_0$. Since $\gcd(D,e)=1$, there is a class $\rho_e$ with
$De\mid a+\rho_ec$ whenever $D\mid a+\rho c$ and $e\mid a+rc$. For these triples $D/v_t=De/u_t$, $De\mid u_t$ and
$u_t/(De)=v_t/D\ge(Q+1)^j$, so Lemma~\ref{lem:key} at the modulus $De$ bounds the group by
$2(Q+1)^{-j}+c_Y+16(L+1)(De)^{-1/2}+(L+1)^3/(2De)$. Sum over $e$ with $\#\{e\le E_0\}\le E_0$,
$\sum_{e\le E_0}e^{-1/2}\le2\sqrt{E_0}$ and $\sum_{e\le E_0}1/e\le1+\log E_0$. For the last claim: $v_t/D\ne1$ divides
$v_t$, which has no divisor in $[2,Q]$.
\end{proof}

\subsection{Duplicates and pruning}

\begin{lemma}[One ray]\label{lem:oneray}
Let $\vartheta,D\ge1$ with $2\vartheta<KY$, and let $(a,c,d)$, $(a',c',d')$ satisfy $Y\le\min(a,c)$, $\max(a,c)\le d$, $Y\le\min(a',c')$, $\max(a',c')\le d'$,
$Kacd-1\le\vartheta D$ and $Ka'c'd'-1\le\vartheta D$. If $D\mid ac'-a'c$, then $ac'=a'c$.
\end{lemma}

\begin{proof}
$KYa^2\le Kacd$ (as $Y\le c$ and $a\le d$), and likewise for $c$, $a'$, $c'$. So all four are at most
$\sqrt{(\vartheta D+1)/(KY)}$, and $0<ac',a'c\le(\vartheta D+1)/(KY)<D$, using $2\vartheta<KY$. Hence $|ac'-a'c|<D$.
\end{proof}

At our parameters $2E_0=2^{5\ell+1}<4\cdot2^{\ell^3}\le KY$ and $E_0<KY^2$.

\begin{lemma}[No duplicates]\label{lem:nodup}
If $t,t'\in\F_r$ have $v_t=v_{t'}=v$ and $v\mid ac'-a'c$, then $t=t'$.
Equivalently, by Lemma~\ref{lem:condev}, two members of $\F_r$ with the same $v$ whose events $\hat e_t$, $\hat e_{t'}$
have a common point are equal.
\end{lemma}

\begin{proof}
Since $u_t=s_tv\le E_0v$, Lemma~\ref{lem:oneray} with $D=v$, $\vartheta=E_0$ gives $ac'=a'c$, and coprimality gives
$(a,c)=(a',c')$. If $d\ne d'$, then $v\mid Kac(d-d')$ and $\gcd(v,Kac)=1$, so $v\le|d-d'|$; say $d\ge v+1$. Then
$s_tv=Kacd-1\ge Kac\,v+Kac-1$, so $s_t\ge Kac\ge KY^2>E_0$, a contradiction. For the equivalence: the class of
$\hat e_t$ is $j\equiv-(rc+a)(P_0c)^{-1}\pmod v$, and two such classes agree if and only if $v\mid ac'-a'c$.
\end{proof}

\emph{Pruning.} A triple $t'\in\F_r$ is a \emph{nested partner} of $t\in\F_r$ if $v_{t'}\mid v_t$, $v_{t'}\ne v_t$ and
$v_{t'}\mid a_tc_{t'}-a_{t'}c_t$; by Lemma~\ref{lem:condev} this says that the events $\hat e_t,\hat e_{t'}$ meet and
the modulus of $\hat e_{t'}$ properly divides that of $\hat e_t$. For finite $\F\subseteq\F_r$ let $\F^\ast$ be the
set of members of $\F$ without a nested partner in $\F$.

\begin{lemma}[Pruning]\label{lem:condprune}
If all members of $\F\subseteq\F_r$ have $v_t\ge2$, then $\mu(\F^\ast)\ge(1-\eps'_Q)\mu(\F)$, and
$\eps'_Q\le2^{-(\ell+8)}$.
\end{lemma}

\begin{proof}
Fix $t'\in\F$. The triples of which $t'$ is a nested partner satisfy $v_{t'}\mid a+\rho c$ with
$\rho\equiv-a_{t'}c_{t'}^{-1}\pmod{v_{t'}}$, and $v_{t'}>Q$ is prime to $P_0$. By Lemma~\ref{lem:condkey} with $j=1$ and
$D=v_{t'}$ they weigh at most $\eps'_{v_{t'}}/v_{t'}\le\eps'_Q/v_{t'}$, since $\eps'_D$ decreases in $D$. Summing over
$t'$ bounds the removed mass by $\eps'_Q\mu(\F)$.

For the numerical bound, $L+1\le2^\ell$ gives $c_Y\le2^{2\ell+10}/Y$. So $2E_0/(Q+1)\le2^{1-123\ell}$,
$E_0c_Y\le2^{7\ell+10-\ell^3}$, $32(L+1)\sqrt{E_0/Q}\le2^{5+\ell-60\ell}$ and
$(L+1)^3(1+\log E_0)/(2Q)\le2^{8\ell-128\ell}$, using $1+\log E_0\le E_0$. Each of these four is at most
$2^{-(\ell+11)}$ for $\ell\ge24$, and $4\cdot2^{-(\ell+11)}\le2^{-(\ell+8)}$.
\end{proof}

\subsection{The events of the conditional family}
The sieve of \S\S\ref{sec:sieve}--\ref{sec:KP} and the cluster estimates of \S\ref{sec:free} use a family of events
only through the following properties, for constants $c_0,c_{\rm r},c_1,c_2,\delta\ge0$:
\begin{enumerate}[label=\textup{(E\arabic*)}]
\item every member $\langle m,c,a\rangle$ is an event: $m\ge1$, $\gcd(m,c)=1$;
\item every modulus is squarefree, $\ge2$, with all prime factors $>Q$ and $\omega(m)\le K_0$;
\item there are no nested pairs: if $m_{\eps'}\mid m_\eps$ and $m_{\eps'}\ne m_\eps$, the events $\eps,\eps'$ hold at no
common $j$;
\item for every $D\ge2$, every $n_0$ and every $j\ge1$,
$\sum D/m_\eps\le c_0+2c_{\rm r}(Q+1)^{-j}+c_1D^{-1/2}+c_2D^{-1}$, the sum running over the events with
$D\mid m_\eps$, $D\mid a_\eps+n_0c_\eps$ and $m_\eps\ge(Q+1)^jD$;
\item two events with equal moduli that hold at a common $j$ are equal;
\item $\sum_{m_\eps=q}1/m_\eps\le\delta$ for every $q$.
\end{enumerate}
The condition $D\mid a_\eps+n_0c_\eps$ says that $\eps$ holds at $n_0$ modulo $D$; it is used with $n_0$ a common point
of a cluster. By (E2), $m_\eps\ne D$ already gives $m_\eps\ge(Q+1)D$, so the case $j=1$ of (E4) covers all events with
$D\mid m_\eps$, $m_\eps\ne D$.

\begin{proposition}[Conditional events]\label{prop:condev}
Let $r\ge0$ be coprime to $P_0$, and let $\F\subseteq\F_r$ be finite, with $v_t\ge2$ good for $t\in\F$ and no
nested pairs in $\F$. Then $t\mapsto\hat e_t$ is injective on $\F$, and $\{\hat e_t:t\in\F\}$ satisfies (E1)--(E6)
with
\begin{gather*}
c_0=E_0c_Y,\qquad c_{\rm r}=E_0,\qquad c_1=32(L+1)\sqrt{E_0},\\
c_2=\frac{(L+1)^3(1+\log E_0)}2,\qquad\delta=\frac{E_0(E_0+1)}{4Y},
\end{gather*}
so that the bound in (E4) is $\eps'_D(j)$.
\end{proposition}

\begin{proof}
(E1) is Lemma~\ref{lem:condev}, (E2) holds as $v_t$ is good, (E3) by assumption and Lemma~\ref{lem:condev}, and (E5) and injectivity are
Lemma~\ref{lem:nodup}. For (E4),
\[
D\mid(rc+a)+n_0\,(P_0c)\iff D\mid a+\rho c,\qquad \rho=r+P_0n_0,
\]
so the class $n_0$ of $\hat e_t$ in the variable $j$ is the class $\rho$ of the triple, and Lemma~\ref{lem:condkey}
applies at the class $\rho$ (if $\gcd(D,P_0)>1$ the sum is empty, since $\gcd(v_t,P_0)=1$). For (E6), let $v_t=q$. Then
$u_t=sq$ with $s\le E_0$, and $a^2Y\le acd=(sq+1)/K\le(sq+1)/4$, likewise for $c$, so for each $s$ there are at most
$(sq+1)/(4Y)\le sq/(2Y)$ such triples, each of weight $1/q$; sum over $s\le E_0$.
\end{proof}

The floor $c_0\le2^{7\ell+10-\ell^3}$ is tiny. This matters in \S\ref{sec:free}, where a vertex with $s$ earlier
neighbours is charged the floor up to $2^{K_0s}$ times.

\section{The mass of the conditional family}\label{sec:mass}

We now show that the conditional family has mass of the full order $L^3/K$. The argument has three steps. Along one
progression in $d$, Brun's pure sieve shows that a proportion $\gg\lambda_K/\ell$ of the cofactors $(Kacd-1)/e$ are good
(\S\ref{sec:brun}). The coprime pairs $(a,c)$ in a class modulo $e$ carry weight $\varphi(e)/e^2$
(Lemma~\ref{lem:wpairs}). Summed over $e$, these weights give a factor $\asymp\ell/\lambda_K$, which cancels the
roughness loss exactly (Lemma~\ref{lem:gain}). The elementary prime estimates with explicit constants used below
(Lemmas~\ref{lem:mertens} and~\ref{lem:mertlow}), in the form of Mertens \cite{Mertens}, are collected in
Appendix~\ref{app:primes}.

\subsection{Brun's pure sieve along affine forms}\label{sec:brun}
We use the pure sieve of Brun \cite{Brun1915,Brun} along affine forms. Put
\[
r_B=2\lfloor\ell/2\rfloor+43,\qquad J_0^\flat=30+\ell+(128\ell+1)\,r_B .
\]
Throughout this subsection $a_1,b_1\ge1$ are coprime integers with $K\mid a_1$, and $I$ is a set of
$X\ge2^{J_0^\flat}$ consecutive integers $\nu\ge0$. Put
\[
\lambda_K=\prod_{p\mid K,\ p\le y_0}\Bigl(1-\frac1p\Bigr)^{-1}\ \ge1 .
\]

\begin{lemma}[Rough values: Brun's pure sieve]\label{lem:brun}
At least $0.0034\,\lambda_KX/\ell$ of the numbers $a_1\nu+b_1$, $\nu\in I$, have no prime factor $\le Q$.
\end{lemma}

\begin{proof}
\emph{Setup.} Let $P$ be the set of primes $\le Q$. Put $g(p)=1/p$ for $p\nmid a_1$, and $g(p)=0$ for
$p\mid a_1$; in the latter case $p\nmid b_1$, so $p\nmid a_1\nu+b_1$. For $S\subseteq P$ let $A_S$ be the number of
$\nu\in I$ with $\prod S\mid a_1\nu+b_1$. This condition is one residue class modulo $\prod S$, or no class if
$\gcd(\prod S,a_1)>1$. Hence $|A_S-X\prod_Sg|\le1$.

\emph{Bonferroni.} For odd $r$ and any finite set $T$, the identity
$\sum_{m\le r}(-1)^m\binom{n+1}m=(-1)^r\binom nr$ gives
$\sum_{S\subseteq T,\,|S|\le r}(-1)^{|S|}\le\mathbf 1[T=\emptyset]$. Apply this with $T$ the set of primes
of $P$ dividing $a_1\nu+b_1$, and sum over $\nu$:
\[
\#\{\text{rough values}\}\ge\sum_{S\subseteq P,\,|S|\le r}(-1)^{|S|}A_S
\ge X\sum_{|S|\le r}(-1)^{|S|}\prod_Sg-\#\{S\subseteq P:|S|\le r\}.
\]

\emph{Main term.} We have $\sum_{S\subseteq P}(-1)^{|S|}\prod_Sg=V:=\prod_P(1-g)$. By Rankin's trick
\cite{Rankin},
\[
\sum_{|S|>r}\prod_Sg\le2^{-r}\sum_S\prod_S(2g)=2^{-r}\prod_P(1+2g).
\]

\emph{Estimates.} The primes $p\mid K$, $p\le y_0\le Q$, divide $a_1$, so $g(p)=0$ there, and $g(p)\le1/p$
elsewhere. Since $Q\ge101$ and $\log Q=128\ell\log2\ge139$, Lemma~\ref{lem:mertens}(c) gives
\[
V\ge\lambda_K\prod_{p\le Q}\Bigl(1-\frac1p\Bigr)\ge\frac{\lambda_K}{3.3\log Q}\ge\frac{\lambda_K}{293\,\ell},\qquad\prod_P(1+2g)\le412\,(128\ell)^2.
\]
The primes of $K$ never divide a value: Brun's density is $0$ there, and $\lambda_K$ records the gain.

\emph{Choice of $r$.} The odd number $r=r_B\ge\ell+42$ gives $2^{r}\ge2^{41}\ell^3$, since
$\ell^3\le2^\ell$ for $\ell\ge10$. As $2^{10}\cdot293\cdot412\cdot128^2\le2^{41}$, the tail is at most
$2^{-10}/(293\ell)$, which is at most $2^{-10}$ times the lower bound for $V$, as $\lambda_K\ge1$.

\emph{Remainder.} It is at most $(r+1)Q^r\le(2Q)^r=2^{(128\ell+1)r}$, which is at most
$X/(2^{30}\ell)$ because $X\ge2^{J_0^\flat}$ and $\ell\le2^\ell$.

\emph{Conclusion.} $X\bigl(\lambda_K(1-2^{-10})/(293\ell)-2^{-30}/\ell\bigr)\ge0.0034\,\lambda_KX/\ell$.
\end{proof}

\begin{lemma}[Non-squarefree values]\label{lem:nsq}
If $a_1\nu+b_1\le U$ for $\nu\in I$, then the number of $\nu\in I$ for which $a_1\nu+b_1$ has no prime factor
$\le Q$ and is not squarefree is at most $X/Q+\sqrt U$.
\end{lemma}

\begin{proof}
If $q^2\mid a_1\nu+b_1$ with $q\ge2$, then $q>Q$, because a prime factor of $q$ divides $a_1\nu+b_1$. Also
$q\le\sqrt U$, and $\gcd(q,a_1)=1$. For each $q$ the admissible $\nu$ form at most one class modulo $q^2$, so there
are at most $X/q^2+1$ of them. Finally $\sum_{q>Q}X/q^2\le X/Q$.
\end{proof}

\begin{lemma}[Many prime factors]\label{lem:omega}
If $a_1\nu+b_1\le X^2$ for $\nu\in I$ and $X\le2^t$, then
\[
2^{K_0+1}\,\#\{\nu\in I:\ \omega(a_1\nu+b_1)>K_0\}\le2X(t+2).
\]
\end{lemma}

\begin{proof}
Since $2^{\omega(n)}\le\tau(n)$, the left side is at most $\sum_{\nu\in I}\tau(a_1\nu+b_1)$. Pairing $e$ with $n/e$ gives
$\tau(n)\le2\#\{e\mid n:e^2\le n\}$, and $n\le X^2$. Hence
\[
\sum_\nu\tau(a_1\nu+b_1)\le2\sum_{e\le X}\#\{\nu:e\mid a_1\nu+b_1\}\le2\sum_{e\le X}\Bigl(\frac Xe+1\Bigr)\le2X(t+2),
\]
since $e\mid a_1\nu+b_1$ is one class of $\nu$ modulo $e$, or none, and $\sum_{e\le X}1/e\le1+t$.
\end{proof}

The average of $\tau$ over the progression replaces Shiu's Brun--Titchmarsh theorem for multiplicative functions \cite{Shiu} here;
this suffices because $2^{K_0}=2^{3\ell}$ beats the logarithmic loss.

\begin{corollary}[Good values]\label{cor:good}
Suppose $a_1\nu+b_1\le U$ for $\nu\in I$, $U(2^{25}\ell^2)^2\le X^2$, and $X\le2^t$ with $t+2\le2^\ell$. Then at
least $\lambda_KX/(300\ell)$ of the numbers $a_1\nu+b_1$, $\nu\in I$, are good. The same holds for the numbers $(Kacd-1)/e$
with $d$ running over $X$ consecutive terms of a progression $d\equiv d_0\pmod e$ with $e\mid Kacd_0-1$, where $U$
bounds these numbers.
\end{corollary}

\begin{proof}
We have $2^{25}\ell^2\le2^{25+2\ell}\le Q$, so $X/Q\le X/(2^{25}\ell^2)$, and $\sqrt U\le X/(2^{25}\ell^2)$. By
Lemma~\ref{lem:nsq} at most $2X/(2^{25}\ell^2)$ rough values are not squarefree. Since $U\le X^2$,
Lemma~\ref{lem:omega} and $2^{34}\ell^2(t+2)\le2^{2\ell}2^\ell$ (as $\ell\ge24$) show that at most $X/(2^{34}\ell^2)$
values have $\omega>K_0=3\ell$. With Lemma~\ref{lem:brun}, $0.0034\lambda_K-2\cdot2^{-25}-2^{-34}\ge\lambda_K/300$ for $\lambda_K\ge1$.

For the progression write $d=d_0+e\nu$. Then $(Kacd-1)/e=Kac\,\nu+b_1$ with $b_1=(Kacd_0-1)/e$, and
$\gcd(Kac,b_1)=1$ because $\gcd(Kac,Kacd_0-1)=1$; and $K\mid Kac$.
\end{proof}

\begin{lemma}[Mass along one progression]\label{lem:massprog}
Let $a,c,e\ge1$ with $\gcd(Kac,e)=1$, let $X\ge4e$, $X'=\lfloor X/e\rfloor-1\ge2^{J_0^\flat}$,
$\lfloor 2KacX/e\rfloor(2^{25}\ell^2)^2\le X'^2$, and $X\le2^t$ with $t+2\le2^\ell$. Then
\[
\sum_{d}\frac e{Kacd-1}\ \ge\ \frac{\lambda_K}{1200\,Kac\,\ell},
\]
the sum running over $d\in[X,2X)$ with $e\mid Kacd-1$ and $(Kacd-1)/e$ good.
\end{lemma}

\begin{proof}
Since $Kac$ is invertible modulo $e$, the $d$ with $e\mid Kacd-1$ form one class modulo $e$, and $[X,2X)$ contains
$X'$ consecutive terms of it. The values $(Kacd-1)/e$ there are at most $2KacX/e$. By
Corollary~\ref{cor:good} at least $\lambda_KX'/(300\ell)$ of them are good, and each term of the sum is at least
$e/(2KacX)$. Finally $eX'\ge X-2e\ge X/2$.
\end{proof}

\subsection{Weighted coprime pairs in a class}

\begin{lemma}\label{lem:wpairs}
Let $A,B,e\ge1$ and $r\in\mathbb Z$. Then
\[
\sum_{\substack{(a,c)\in[A,2A)\times[B,2B)\\ \gcd(a,c)=1,\ e\mid a+rc}}\frac1{ac}
\ \ge\ \frac{\varphi(e)}{e^2}\Bigl(\frac18-\frac{e(A+B)}{4AB}-\frac{e(3A+B)(1+\log2A)}{4AB}-\frac{3e^2}{2B}\Bigr).
\]
In particular the sum is at least $\varphi(e)/(16e^2)$ if $64e(2+\log2A)\le\min(A,B)$ and $96e^2\le B$.
\end{lemma}

\begin{proof}
\emph{Pairs with $\gcd(e,c)=1$.} There are at least $\varphi(e)\lfloor B/e\rfloor$ integers $c\in[B,2B)$ prime to
$e$, and for each of them $e\mid a+rc$ is one class of $a$ modulo $e$, with at least $\lfloor A/e\rfloor$ members
in $[A,2A)$. So these pairs number at least $\varphi(e)\lfloor A/e\rfloor\lfloor B/e\rfloor$.

\emph{Non-coprime pairs.} Such a pair has a common prime $p\le2A$, and $p\nmid e$ since $\gcd(e,c)=1$. Writing
$a=pa'$, $c=pc'$, the $c'$ lie in an interval of length $B/p+2$ and are prime to $e$, and for each of them $a'$ lies
in one class modulo $e$ (as $\gcd(p,e)=1$). So there are at most
$\varphi(e)(\lfloor(\lfloor B/p\rfloor+2)/e\rfloor+1)(\lfloor A/(pe)\rfloor+1)$ of them. Summing over $p\le2A$ with
$\sum_pp^{-2}\le\frac12$ (Lemma~\ref{lem:mertens}(d)), $\sum_{p\le2A}1/p\le1+\log2A$ and at most $2A$ primes,
the coprime pairs of the class number at least $\varphi(e)\frac{AB}{e^2}(\frac12-4\delta)$, where $\delta$ is the
sum of the three error terms in the statement.

\emph{Weights.} Each pair has $1/(ac)>1/(4AB)$. Under the conditions of the simple form, $\delta\le\frac1{64}+
\frac1{64}+\frac1{256}<\frac1{16}$.
\end{proof}

Every pair counted here has $\gcd(ac,e)=1$ when $\gcd(r,e)=1$, since a prime of $e$ dividing $c$ divides $a$. The factor
$\varphi(e)/e$ is intrinsic, not a loss of the method: the coprime pairs of the class have density
$e^{-1}\prod_{p\mid e}(1-\frac1p)\prod_{p\nmid e}(1-\frac1{p^2})$.

\subsection{The mass before pruning}

\begin{proposition}\label{prop:massF0}
Let $r$ be an integer coprime to $P_0$. There is a finite set $\F_r^\flat\subseteq\F_r$, all of whose members have
$v_t\ge2$ good, with
\[
\mu(\F_r^\flat)\ \ge\ \frac{\lambda_KL^3}{1382400\,K\ell}\sum_{e\mid P_0,\ e\le E_0}\frac{\varphi(e)}{e^2}.
\]
\end{proposition}

\begin{proof}
\emph{Scales.} Put $G_0=16\ell$, $H_0=3\ell+4$, $J_1=J_0^\flat+5\ell+1$, $\Delta=L-H_0-J_1-G_0-4\ell^3$ and
$w=\lfloor\Delta/4\rfloor$. Take
\[
i,i'\in[\ell^3,\ell^3+w),\qquad j\in[J_0^\flat+G_0+5\ell+1+s,\ L-H_0-s),\qquad s=i+i'.
\]
Put $A=2^i$, $B=2^{i'}$, $X=2^j$. For $e\mid P_0$ with $e\le E_0$ let $T(i,i',j,e)$ be the set of triples $(a,c,d)$
with $(a,c)\in[A,2A)\times[B,2B)$, $\gcd(a,c)=1$, $e\mid a+rc$, $d\in[X,2X)$, $e\mid Kacd-1$ and $(Kacd-1)/e$ good.
Let $\F_r^\flat$ be the union of all these sets.

\emph{Membership.} A triple of $T(i,i',j,e)$ lies in $\F_r$ with $s_t=e$ and $v_t=u_t/e$: $a,c\ge2^{\ell^3}=Y$;
$a,c<2^{j}\le d$ as $i,i'<j$; $Kacd<2^{3\ell+s+j+3}\le2^L$ as $K\le2^{3\ell}$ and $j+s<L-H_0$; $u_t/e$ is $Q$-rough,
hence prime to $P_0$, so $\gcd(u_t,P_0)=e$ (as $e\mid P_0$ is squarefree); and $e\mid a+rc$, $e\le E_0$. Also
$v_t\ge2$, since $u_t\ge KY^3-1>2E_0$.

\emph{Size conditions.} For $1\le e\le E_0=2^{5\ell}$ the hypotheses of Lemmas~\ref{lem:massprog}
and~\ref{lem:wpairs} hold on these scales, with $t=j$: $X/e\ge2^{j-5\ell}\ge2^{J_0^\flat+1}$; the condition on
$2KacX/e$ follows from $j\ge s+G_0+5\ell+1+J_0^\flat$, since $2KacX/e\cdot(2^{25}\ell^2)^2\cdot(2e)^2/X^2\le
2^{8\ell+s+55+4\log_2\ell}/X\le1$; $j+2\le L<2^\ell$; and $A,B\ge Y=2^{\ell^3}$ is much larger than $96E_0^2$ and
$64E_0(2+\log2A)$.

\emph{One scale and one $e$.} Since $\gcd(r,P_0)=1$, every pair of the class has $\gcd(ac,e)=1$, and
$\gcd(K,e)=1$ as $e\mid P_0$. By Lemma~\ref{lem:massprog} each pair contributes $\sum_d1/v_t=\sum_de/u_t\ge
\lambda_K/(1200Kac\ell)$, and by Lemma~\ref{lem:wpairs}
\[
\sum_{t\in T(i,i',j,e)}\frac1{v_t}\ \ge\ \frac{\varphi(e)}{16e^2}\cdot\frac{\lambda_K}{1200K\ell}=\frac{\lambda_K\varphi(e)}{19200\,K\ell e^2}.
\]

\emph{Counting scales.} The sets $T(i,i',j,e)$ are pairwise disjoint: $e=s_t$, and the dyadic blocks of $a,c,d$
determine $i,i',j$. The $j$-range above is the interval $[J_1+G_0+s,\,L-H_0-s)$. With $i=\ell^3+x$, $i'=\ell^3+x'$,
this interval has $\Delta-2x-2x'$ elements, and
\[
\sum_{x,x'<w}(\Delta-2x-2x')=w^2(\Delta-2w+2)\ge\frac{(\Delta-3)^2\Delta}{32}.
\]
For $\ell\ge24$ one has $100(H_0+J_1+G_0+4\ell^3)\le23\cdot2^{\ell-1}\le23L$, so $\Delta\ge0.77L$ and, as
$\Delta\ge6\cdot10^6$, this is at least $L^3/71\ge L^3/72$. Multiply, with $72\cdot19200=1382400$.
\end{proof}

\subsection{The divisor sum}

\begin{lemma}[The gain]\label{lem:gain}
Let $y\ge3$ and $\vartheta>0$ be real with $\log \vartheta\ge5\log y$, let $\mathcal K$ be a finite set of primes, and let $P$ be
the product of the primes $p\le y$ not in $\mathcal K$. Then
\[
\sum_{e\mid P,\ e\le \vartheta}\frac{\varphi(e)}{e^2}\ \ge\ \frac{e^{-7}\log y}{6\prod_{p\in\mathcal K}(1+1/p)}.
\]
\end{lemma}

\begin{proof}
Let $S$ be the set of primes of $P$. By multiplicativity the sum is
$\sum_{R\subseteq S,\ \prod R\le \vartheta}\prod_{p\in R}\frac{p-1}{p^2}$.

\emph{Full sum.} $\prod_S(1+\frac{p-1}{p^2})=\prod_S(1+\frac1p)(1-\frac1{p(p+1)})\ge\frac{37}{72}\prod_S(1+\frac1p)$,
by $\prod(1-t_i)\ge1-\sum t_i$ and $\sum_pp^{-2}\le\frac{35}{72}$.

\emph{Tail (Rankin).} Let $\sigma=1/\log y$. Since $\frac{p-1}{p^2}\le\frac1p$, the terms with $\prod R>\vartheta$ weigh at
most $\vartheta^{-\sigma}\prod_S(1+p^{\sigma-1})$. Here $1+p^{\sigma-1}\le(1+\frac1p)\exp(\frac{p^\sigma-1}p)$, and
$p^\sigma-1\le(e-1)\sigma\log p$ for $p\le y$ by convexity of $e^x$ on $[0,1]$. With
Lemma~\ref{lem:mertens}(a), the tail is at most
\[
\prod_S\Bigl(1+\frac1p\Bigr)\exp\Bigl((e-1)\Bigl(1+\frac{\log4}{\log y}\Bigr)-\frac{\log \vartheta}{\log y}\Bigr)
\le e^{3.9-5}\prod_S\Bigl(1+\frac1p\Bigr)\le0.335\prod_S\Bigl(1+\frac1p\Bigr),
\]
since $(e-1)(1+\log4/\log3)\le3.9$.

\emph{Conclusion.} The truncated sum is at least $(\frac{37}{72}-0.335)\prod_S(1+\frac1p)\ge\frac16\prod_S(1+\frac1p)$,
and $\prod_S(1+\frac1p)\ge e^{-7}\log y/\prod_{\mathcal K}(1+\frac1p)$ by Lemma~\ref{lem:mertlow}(c).
\end{proof}

The sum is $\asymp\log y$: the upper bound $\prod_{p\le y}(1+1/p)\ll\log y$ follows from Lemma~\ref{lem:mertens}(b).

\subsection{The mass}

\begin{proposition}[Mass]\label{prop:massF}
Let $r$ be coprime to $P_0$. Then the pruned family $(\F_r^\flat)^\ast$ satisfies
\[
\mu\bigl((\F_r^\flat)^\ast\bigr)\ge X_m:=\frac{L^3}{2^{37}K} .
\]
\end{proposition}

\begin{proof}
By Lemma~\ref{lem:condprune}, $\mu((\F_r^\flat)^\ast)\ge(1-2^{-(\ell+8)})\mu(\F_r^\flat)\ge\frac12\mu(\F_r^\flat)$. By
Lemma~\ref{lem:gain} with $y=y_0=2^\ell$, $\vartheta=E_0$ and $\mathcal K=\{p\le y_0:\ p\mid K\}$, for which
$\log E_0=5\log y_0$, the divisor sum in Proposition~\ref{prop:massF0} is at least
$e^{-7}\ell\log2/(6\prod_{p\mid K,\,p\le y_0}(1+\frac1p))$. Since $(1-\frac1p)(1+\frac1p)\le1$,
$\lambda_K\ge\prod_{p\mid K,\,p\le y_0}(1+\frac1p)$, and the factor $\ell$ cancels:
\[
\mu\bigl((\F_r^\flat)^\ast\bigr)\ge\frac12\cdot\frac{\lambda_KL^3}{1382400\,K\ell}\cdot\frac{e^{-7}\ell\log2}{6\lambda_K}
\ge\frac{L^3}{2^{37}K},
\]
using $e^7\le4096$ and $\log2\ge\frac12$.
\end{proof}

The roughness at $Q=2^{128\ell}$ keeps the fraction $\asymp1/\ell$ of each scale (Lemma~\ref{lem:brun}), and the
averaging over the primes up to $y_0=2^\ell$ recovers the factor $\ell$ (Lemma~\ref{lem:gain}). The primes of $K$ are
not in $P_0$, which costs the divisor sum $\prod_{p\mid K}(1+\frac1p)$; but they never divide a modulus, which Brun's
sieve pays back as $\lambda_K$. The mass is therefore of the order of the number $\asymp L^3$ of scales, divided by $K$,
with no other loss.

\section{Cluster estimates}\label{sec:free}

In this section $\Ev$ is a finite family of events satisfying (E1)--(E6) of \S\ref{sec:avg} with nonnegative constants,
$k\ge1$, and $\mu\ge\sum_{\eps\in\Ev}1/m_\eps$. The gas is that of \S\ref{sec:sieve}, on $\U=\{1,2\}\times\Ev$, and
$\Mlcm=\lcm_\eps m_\eps$. For $v\in\U$ we write $m_v$, $a_v$, $c_v$ for the data of its event. Everything reduces to a
multi-site bound (defined below): the sum of the weights of the connected sets through a vertex, with at least two
elements, is a small multiple of the weight of the vertex.

\subsection{Polymers}
For $C\subseteq\U$ let $\mathrm{cnt}(C)$ be the number of $0\le j<\Mlcm$ at which every event of $C$ holds; $C$ is
\emph{solvable} if $\mathrm{cnt}(C)>0$. A polymer $C$ has $W_{y,z}(C)\ne0$ only if it is independent in each copy and
has at most $2k$ vertices. The only connected such set containing both copies of an event $\eps$ is the
\emph{double site} $\{(1,\eps),(2,\eps)\}$. Indeed, a vertex adjacent to either copy shares a prime with $m_\eps$,
and hence with the copy of $\eps$ in its own half, which contradicts independence. The double site has
$W=yz/m_\eps$, and singletons have $W=y/m_\eps$ or $z/m_\eps$, since $\mathrm{cnt}=\Mlcm/m_\eps$. Call $C$ \emph{good} if it
is independent in each copy, contains no double site, and has $\le2k$ vertices. In a good set every prime divides
at most two vertices, one in each copy.

\emph{Link factor.} Let $C$ be solvable, let $x\in\U$, and put $D=\gcd(m_x,\lc C)$.
\begin{enumerate}[label=\textup{(\roman*)}]
\item By the Chinese remainder theorem, $m_x\,\mathrm{cnt}(C\cup x)=D\,\mathrm{cnt}(C)$ if $C\cup x$ is
solvable, and $\mathrm{cnt}(C\cup x)=0$ otherwise.
\item If $n_0$ is any common point of $C$ and $C\cup x$ is solvable, then $D\mid a_x+n_0c_x$.
\end{enumerate}

\emph{Weights.} For $\eta>0$ and $\theta\ge1$ put
\[
g(C)=g_{\eta,\theta}(C)=\mathbf 1[C\text{ good}]\ \eta^{|C|}\,\theta^{\omega(\lc C)}\,\frac{\mathrm{cnt}(C)}{\Mlcm}.
\]
If $|y|,|z|\le\eta$ and $C$ is good, then $|W_{y,z}(C)|\theta^{\omega(\lc C)}\le g(C)$, and
$g(\{v\})=\eta\theta^{\omega(m_v)}/m_v\le\eta\theta^{K_0}/m_v$. We say that \emph{$(\eta,\theta,\alpha)$ is a multi-site
bound} if for every $v\in\U$
\[
\sum_{\substack{C\ni v\ \text{connected}\\|C|\ge2}}g_{\eta,\theta}(C)\le\frac\alpha{1-\alpha}\,g_{\eta,\theta}(\{v\}).
\]

\begin{corollary}[Per-site polymer sums]\label{cor:site}
Let $(\eta,\theta,\alpha)$ be a multi-site bound with $0\le\alpha<1$, and $|y|,|z|\le\eta$. Then for every $v\in\U$,
\[
\sum_{C\ni v}|W_{y,z}(C)|\,\theta^{\omega(\lc C)}\le\theta^{K_0}\Bigl(\eta+\eta^2+\eta\frac{\alpha}{1-\alpha}\Bigr)\frac1{m_v},
\]
the sum running over polymers $C$.
\end{corollary}

\begin{proof}
The polymers through $v$ with $W\ne0$ are $\{v\}$, the double site of $v$, and good connected sets with at least two
vertices. The first two contribute at most $\theta^{K_0}(\eta+\eta^2)/m_v$. For the last kind $|W|\theta^{\omega}\le g$,
and the multi-site bound together with $g(\{v\})\le\eta\theta^{K_0}/m_v$ gives the third term.
\end{proof}

For an integer $q\ge2$, let $\mu_q=\sum_{\eps\in\Ev,\ q\mid m_\eps}1/m_\eps$ be the mass through $q$.

\begin{lemma}\label{lem:Mq}
For every $q\ge2$, $\mu_q\le c_0+2c_{\rm r}(Q+1)^{-1}+c_1q^{-1/2}+c_2q^{-1}+\delta$.
\end{lemma}

\begin{proof}
\emph{Events with $m_\eps\ne q$.} If $q\mid m_\eps$, then $c_\eps$ is invertible modulo $q$ and $q\mid a_\eps+\rho c_\eps$
with $\rho\equiv-a_\eps c_\eps^{-1}$. For each of the $q$ classes $\rho$, (E4) with $D=q$ and $j=1$ bounds $\sum q/m_\eps$ by
$c_0+2c_{\rm r}(Q+1)^{-1}+c_1q^{-1/2}+c_2q^{-1}$ (by (E2), $m_\eps\ne q$ gives $m_\eps\ge(Q+1)q$). Summing over the $q$
classes gives the first four terms. \emph{Events with $m_\eps=q$.} By (E6) they weigh at most $\delta$.
\end{proof}

\begin{proposition}\label{prop:polymer}
Let $\tau\ge0$, $\theta=e^\tau$ and $a(C)=\tau\,\omega(\lc C)$. Suppose that $(1,\theta,\alpha)$ and $(2,1,\alpha_2)$ are
multi-site bounds with $0\le\alpha,\alpha_2<1$.
\begin{enumerate}[label=\textup{(\alph*)}]
\item \textup{(KP condition.)} If
$2\theta^{K_0}\bigl(2+\tfrac{\alpha}{1-\alpha}\bigr)m_P\le\tau$ with $m_P\ge\sup_{q>Q}\mu_q$, then the
activities $W_{-1,-1}$ satisfy the hypothesis of Corollary~\ref{cor:KP} with weight $a$.
\item \textup{(First order.)} With $\mu_\Ev=\sum_{\eps\in\Ev}1/m_\eps$,
\[
\sum_CW_{-1,-1}(C)+\sum_C|W_{-1,-1}(C)|(e^{a(C)}-1)\le
-\mu_\Ev+3(\theta^{K_0}-1)\mu_\Ev+4\theta^{K_0}\tfrac{\alpha}{1-\alpha}\mu_\Ev.
\]
\item \textup{(Positive activities.)} We have
$\log Z(2,2)\le\sum_CW_{2,2}(C)\le2\mu_\Ev\bigl(4+2\tfrac{\alpha_2}{1-\alpha_2}\bigr)$.
\end{enumerate}
\end{proposition}

\begin{proof}
(a) If $D\not\sim C$, some vertex of $D$ has modulus divisible by a prime $q\mid\lc C$, and $q>Q$ by (E2).
Summing Corollary~\ref{cor:site} (with $\eta=1$) over the $\le\omega(\lc C)$ such primes, over the events
$\eps$ with $q\mid m_\eps$, and over both copies, gives
\[
\sum_{D\not\sim C}|W(D)|e^{a(D)}\le\omega(\lc C)\cdot2\theta^{K_0}\bigl(2+\tfrac{\alpha}{1-\alpha}\bigr)m_P\le a(C).
\]

(b) The singletons contribute $-2\mu_\Ev$ and the double sites contribute $+\mu_\Ev$ to $\sum W$. Their terms
$|W|(e^{a}-1)$ add at most $3(\theta^{K_0}-1)\mu_\Ev$, since $\omega(m_\eps)\le K_0$. Larger polymers contribute at most
$2\sum|W|e^{a}$, which by Corollary~\ref{cor:site} is at most
$2\sum_v\theta^{K_0}\frac{\alpha}{1-\alpha}\frac1{m_v}=4\theta^{K_0}\frac{\alpha}{1-\alpha}\mu_\Ev$.

(c) For nonnegative activities,
$Z\le\prod_C(1+W(C))\le\exp\sum_CW(C)$. Singletons and double sites give $4\mu_\Ev+4\mu_\Ev$. Larger polymers
give at most $2\cdot2\mu_\Ev\frac{\alpha_2}{1-\alpha_2}$, by distributing each polymer over its sites and
using the multi-site bound $(2,1,\alpha_2)$.
\end{proof}

It remains to produce multi-site bounds. A step-by-step bound, which attaches one vertex at a time to a connected
cluster, has to handle \emph{forced} vertices: a vertex $x$ attached through $D=\gcd(m_x,\lc C)$ contributes $D/m_x$,
and when $m_x=D$ this is $1$, with no decay. Prepaying such attachments by a potential in the open primes of the
cluster multiplies all other terms by a factor $2^{\Theta(\ell^2)}$, which the ray term $1/Q$ of the key lemma can
absorb only if $\log_2Q\gg\ell^2$; with $\log Q\asymp\ell^2$ the roughness would keep only the fraction $\asymp\ell^{-2}$
of the mass, and the averaging of \S\ref{sec:avg} would recover just one factor $\ell$. Below there is no potential and
no prepayment. Each vertex is charged only its own factor $D/m_x$, the vertices of a cluster are summed in an order
chosen for its shape, and a vertex that is forced in one order has later neighbours in others, which make its
cofactor large; (E4) turns a large cofactor into decay.

\subsection{Shapes and the free-order Fubini}\label{sec:fubini}

\emph{The prime graph.} For a good set $C$ join $x,z\in C$ when $D_{xz}=\gcd(m_x,m_z)>1$. Every prime divides at most two
vertices of a good set, so the $D_{xz}$ at a fixed $x$ are pairwise coprime. As the moduli are
squarefree by (E2), $\lc C=\prod_xm_x/\prod_{xz}D_{xz}$. If $C$ is solvable its common points form one class modulo
$\lc C$, hence
\begin{equation}\label{eq:crt}
\frac{\mathrm{cnt}(C)}{\Mlcm}=\prod_{x\in C}\frac1{m_x}\prod_{xz}D_{xz}.
\end{equation}
Fix an order of $C$ with the root $v$ first, and attribute each $D_{xz}$ to the later endpoint. For $x\in C$ let
$C_{<x}$ be the set of vertices before $x$ and $D_<(x)=\prod_{z\in C_{<x},\,z\sim x}D_{xz}=\gcd(m_x,\lc C_{<x})$. Then
\eqref{eq:crt} becomes $\mathrm{cnt}(C)/\Mlcm=\prod_xD_<(x)/m_x$, and, since $\omega(\lc C)\le\omega(m_v)+K_0(|C|-1)$,
\begin{equation}\label{eq:factor}
g(C)\le g(\{v\})\prod_{x\in C\setminus v}\eta\theta^{K_0}\,\frac{D_<(x)}{m_x}.
\end{equation}
A vertex costs $D_<(x)/m_x$, the reciprocal of its \emph{cofactor} $m_x/D_<(x)$. The prefixes $C_{<x}$ need not be
connected. Formally \eqref{eq:factor} is obtained by applying the link factor (i) along the order; every prefix of a good
solvable set is good and solvable.

\emph{Shapes.} A \emph{shape} $H$ on $m$ non-root vertices is a graph on $\{0,1,\dots,m\}$, where $0$ is the root,
together with a bit $p(i)\in\{0,1\}$ for each $i\ge1$ (and $p(0)=0$). A tuple $x=(x_1,\dots,x_m)$ of elements of $\U$
\emph{realizes} $H$ with root $v=x_0$ if $x_0,\dots,x_m$ are distinct, $C=\{x_0,\dots,x_m\}$ is good and solvable,
$i\sim_Hi'$ exactly when $\gcd(m_{x_i},m_{x_{i'}})>1$, and $p(i)=1$ exactly when $m_{x_i}$ has a prime factor dividing
no other modulus of $C$ (a \emph{private} prime). For $i\ge1$ let $t_e(i)$ and $t_l(i)$ be the numbers of neighbours
of $i$ in $H$ with smaller and with larger index, so the root counts as an earlier neighbour, and put
$j(i)=t_l(i)+p(i)$, the \emph{cofactor exponent}.

\begin{lemma}[Cofactors]\label{lem:cof}
Let $x$ realize $H$, $1\le i\le m$, $P=\{x_0,\dots,x_{i-1}\}$ and $D=\gcd(m_{x_i},\lc P)$. Then
$D\ge(Q+1)^{t_e(i)}$ and $(Q+1)^{j(i)}D\le m_{x_i}$, and $m_{x_i}=D$ if $j(i)=0$.
\end{lemma}

\begin{proof}
Each neighbour of $x_i$ shares a prime with it, and distinct neighbours share distinct primes, since no prime divides
three vertices. The primes shared with earlier neighbours divide $D$. A prime shared with a later neighbour divides no
third vertex, so it does not divide $\lc P$; nor does a private prime. So $m_{x_i}/D$ is divisible by at least
$j(i)$ distinct primes, and $D$ by at least $t_e(i)$, all greater than $Q$. If $j(i)=0$, every prime of $m_{x_i}$
divides another vertex, necessarily an earlier one, so $m_{x_i}\mid\lc P$ and $m_{x_i}=D$.
\end{proof}

So a vertex is cheap in two ways: through a large $D_<$, when it has many earlier neighbours (the $c_1D^{-1/2}$ and
$c_2D^{-1}$ terms of (E4)), or through a large cofactor, when it has many later ones (the ray term $(Q+1)^{-j}$). It is
expensive only when it is \emph{free} ($t_e=0$; it costs the whole mass) or \emph{exact} ($j=0$; it is determined by its
predecessors and costs $1$).

\begin{lemma}[Slot bound]\label{lem:slot}
Let $P\subseteq\U$, $E\subseteq P$, $s=|E|$ and $j\ge0$. Call $y\in\U$ a \emph{candidate} if $y\notin P$, $P\cup\{y\}$ is
good and solvable, $\gcd(m_y,m_z)>1$ exactly for the $z\in P$ that lie in $E$, and $D_y=\gcd(m_y,\lc P)$ satisfies
$(Q+1)^jD_y\le m_y$, with $m_y=D_y$ if $j=0$. Then $\sum_yD_y/m_y\le\Sigma_{s,j}$, the sum running over candidates,
where
\[
\Sigma_{0,j}=2\mu,\qquad\Sigma_{s,0}=2\cdot2^{K_0s},\qquad
\Sigma_{s,j}=2\cdot2^{K_0s}\Bigl(c_0+\frac{2c_{\rm r}}{(Q+1)^j}+\frac{c_1}{(Q+1)^{s/2}}+\frac{c_2}{(Q+1)^s}\Bigr)
\]
for $s,j\ge1$. No assumption on $P$ is made.
\end{lemma}

\begin{proof}
If there is no candidate there is nothing to prove. Otherwise $P$ is solvable; fix a common point $n_0$ of its events.
For every candidate $y$ the link factor (ii) gives $D_y\mid a_y+n_0c_y$.

\emph{Free vertex, $s=0$.} No prime of $m_y$ divides a vertex of $P$, so $D_y=1$, and the sum is at most
$\sum_y1/m_y\le2\mu$ over both copies.

\emph{$s\ge1$.} A prime of $D_y$ divides some $w\in P$ with $\gcd(m_y,m_w)>1$, so $w\in E$. Hence $D_y$ is the product of
a subset of the at most $K_0s$ primes of the moduli of $E$, and takes at most $2^{K_0s}$ values. Fix one, $D$. As in
Lemma~\ref{lem:cof}, $D\ge(Q+1)^s\ge2$. If $j=0$, the candidates with $D_y=D$ have modulus $D$ and class $n_0$ modulo
$D$; by (E5) they are the two copies of one event, and each contributes $1$. If $j\ge1$, they have $m_y\ne D$,
$(Q+1)^jD\le m_y$ and $D\mid a_y+n_0c_y$, so (E4), summed over the two copies, bounds their contribution by
$2(c_0+2c_{\rm r}(Q+1)^{-j}+c_1D^{-1/2}+c_2D^{-1})$, and $D\ge(Q+1)^s$.
\end{proof}

The bound is uniform in the earlier vertices, whether or not they are connected to $y$. So the cluster sum can be done as
an iterated sum.

\begin{proposition}[Free-order Fubini]\label{prop:fubini}
For every shape $H$ on $m$ non-root vertices and every $v\in\U$,
\[
\sum_{x\text{ realizes }H}g(\{x_0,\dots,x_m\})\ \le\ g(\{v\})\prod_{i=1}^m\eta\theta^{K_0}\,\Sigma_{t_e(i),\,j(i)} .
\]
\end{proposition}

\begin{proof}
For $1\le i\le m$ put $P_i=\{x_0,\dots,x_{i-1}\}$ and let $E_i$ be the set of $x_{i'}$, $i'<i$, with $i'\sim_Hi$. Let
$f_i(x_1,\dots,x_i)=\eta\theta^{K_0}\gcd(m_{x_i},\lc P_i)/m_{x_i}$ if $|E_i|=t_e(i)$ and $x_i$ is a candidate for
$(P_i,E_i,j(i))$, and $f_i=0$ otherwise. If $x$ realizes $H$, then $|E_i|=t_e(i)$, every $x_i$ is a candidate by
Lemma~\ref{lem:cof}, and \eqref{eq:factor}, for the order of the indices, bounds the summand by
$g(\{v\})\prod_if_i$. The $f_i$ are nonnegative and $f_i$ depends only on $x_1,\dots,x_i$. Drop the realization
condition, sum over all of $\U^m$, and sum over $x_m$ first, then over $x_{m-1}$, and so on. By
Lemma~\ref{lem:slot} each inner sum is at most $\eta\theta^{K_0}\Sigma_{t_e(i),j(i)}$, uniformly in the outer variables.
\end{proof}

Nothing here uses an attachment sequence: a vertex is bounded given all earlier vertices, whether or not they are
connected to it, and a vertex without earlier neighbours costs the whole mass $2\mu$. The price is that the order is
fixed by the indexing of $H$; the next subsection chooses it.

\subsection{Averaging over orders}\label{sec:avgorders}

\begin{lemma}[Admissible shapes]\label{lem:admshape}
Let $v\in\U$ and let $C\ni v$ be connected, good and solvable, with $m+1\ge2$ elements. Then every enumeration $x$ of
$C\setminus\{v\}$ realizes its shape, the shape $H$ defined by the moduli of $x$. Moreover $H$ is \emph{admissible}:
every non-root vertex has degree between $1$ and $K_0$, and a non-root vertex of degree $1$ has $p=1$. Hence the sum
of $g(C)$ over such $C$ is at most the sum over admissible $H$ of the left side of Proposition~\ref{prop:fubini}.
This passage to tuples counts each $C$ once for each of its $m!$ enumerations; we keep the upper bound and discard
the factor $1/m!$, which is slack.
\end{lemma}

\begin{proof}
Degrees are positive by connectivity. Distinct neighbours of $y$ share distinct primes of $m_y$, so the degree is at
most $\omega(m_y)\le K_0$. Suppose $y$ has a single neighbour $z$ and no private prime. Then every prime of $m_y$
divides $m_z$, so $m_y\mid m_z$, and the events of $y$ and $z$ meet. If $m_y=m_z$, then by (E5) $y$ and $z$ are the
two copies of one event, a double site, which a good set does not contain. Otherwise $y$ and $z$ form a nested pair,
which (E3) excludes.
\end{proof}

An \emph{order} is a permutation $o$ of $\{1,\dots,m\}$; it puts the non-root vertex $i$ of $H$ at position $o(i)$.
Let $H_o$ be the shape in which the vertex $i$ of $H$ gets the index $o(i)$, with its bit $p(i)$, the root staying at
$0$. Then $x\mapsto x\circ o^{-1}$ is a bijection from the realizations of $H$ onto those of $H_o$ that preserves
$\{x_0,\dots,x_m\}$. So Proposition~\ref{prop:fubini}, applied to $H_o$, bounds the
realization sum of $H$, for every $o$: the order of summation may be chosen separately for every shape. For a
non-root vertex $i$ of $H$ let $d_i$ be its degree, $\varrho_i=\mathbf 1[i\sim_H0]$ its root adjacency and $N_i$ its set
of non-root neighbours, and let $\mathrm{rk}_o(N,i)=\#\{b\in N:\ o(b)<o(i)\}$. In $H_o$ the vertex $i$ has
$t_e=\varrho_i+\mathrm{rk}_o(N_i,i)$ and $t_l=d_i-t_e$.

\begin{lemma}[Uniform rank]\label{lem:rank}
Let $N\subseteq\{1,\dots,m\}$ and $i\notin N$. For every $s\in\{0,\dots,|N|\}$ exactly $m!/(|N|+1)$ permutations $o$
have $\mathrm{rk}_o(N,i)=s$. Consequently, in a uniformly random order, $t_e$ of a non-root vertex of degree $d$ and
root adjacency $\varrho$ is uniform on $\{\varrho,\dots,d\}$.
\end{lemma}

\begin{proof}
Let $S=N\cup\{i\}$. For $a\in S$, composing with the transposition of $i$ and $a$ gives
$\mathrm{rk}_{o\circ(i\,a)}(N,i)=\mathrm{rk}_o(S,a)$, and $o\mapsto o\circ(i\,a)$ is a bijection. So for any
function $G$,
\[
(|N|+1)\sum_oG(\mathrm{rk}_o(N,i))=\sum_o\sum_{a\in S}G(\mathrm{rk}_o(S,a))=m!\sum_{s=0}^{|N|}G(s),
\]
because the ranks of the elements of $S$ in a fixed order are $0,\dots,|S|-1$, each once. Take $G$ an indicator.
\end{proof}

\begin{lemma}\label{lem:avgorder}
Let $N_1,\dots,N_m\subseteq\{1,\dots,m\}$ with $i\notin N_i$, and $G_1,\dots,G_m:\mathbb N\to\mathbb R$. Some $o$ satisfies
\[
\sum_{i=1}^mG_i(\mathrm{rk}_o(N_i,i))\le\sum_{i=1}^m\frac1{|N_i|+1}\sum_{s=0}^{|N_i|}G_i(s).
\]
\end{lemma}

\begin{proof}
By Lemma~\ref{lem:rank} the right side is the average of the left side over all $o$.
\end{proof}

This is where the step-by-step method fails and the random order succeeds: a vertex that is exact in one order has
later neighbours, and hence a large cofactor, in most others, so no vertex has to be prepaid.

\begin{lemma}[Counting shapes]\label{lem:shapecount}
If $0\le\lambda$ and $\lambda(m+1)\le1$, then $\sum_H\lambda^{\sum_{i\ge1}\deg(i)}\le4^m$, the sum running over all
shapes on $m$ non-root vertices.
\end{lemma}

\begin{proof}
Every edge has a non-root endpoint, so $\sum_{i\ge1}\deg(i)$ is at least the number of edges, and $\lambda\le1$. A
shape is determined by its edge set and its $m$ bits, so the sum is at most $2^m(1+\lambda)^{m(m+1)/2}$, and
$(1+\lambda)^{m(m+1)/2}\le e^{\lambda(m+1)m/2}\le e^{m/2}\le2^m$.
\end{proof}

\begin{proposition}[Multi-site bounds from shapes]\label{prop:sitefree}
Let $A_*>0$ and $0<\lambda\le1/(2k)$, and suppose that for all integers $1\le d\le K_0$ and $p,\varrho\in\{0,1\}$ with
$\varrho\le d$, and $p=1$ if $d=1$,
\begin{equation}\label{eq:W}
\frac1{d+1-\varrho}\sum_{s=\varrho}^d\log\bigl(\eta\theta^{K_0}\,\Sigma_{s,\,d-s+p}\bigr)\le\log A_*+d\log\lambda .
\end{equation}
Then for every $v\in\U$ and $m\ge1$ the connected sets $C\ni v$ with $m+1$ elements satisfy
$\sum_Cg(C)\le(4A_*)^mg(\{v\})$. If $\alpha=4A_*<1$, then $(\eta,\theta,\alpha)$ is a multi-site bound.
\end{proposition}

\begin{proof}
Only good solvable $C$ contribute, and a good set has at most $2k$ vertices, so $\lambda(m+1)\le1$. By
Lemma~\ref{lem:admshape} it suffices to bound, for each admissible $H$, the realization sum of $H$ by
$A_*^m\lambda^{\sum_{i\ge1}\deg(i)}g(\{v\})$, and then to apply Lemma~\ref{lem:shapecount}. For a non-root vertex $i$ of
$H$ let $d_i,\varrho_i,p_i,N_i$ be as above, so $|N_i|=d_i-\varrho_i$, and
$G_i(s)=\log(\eta\theta^{K_0}\Sigma_{\varrho_i+s,\,d_i-\varrho_i-s+p_i})$. Admissibility says that $(d_i,p_i,\varrho_i)$ is one of the
types in \eqref{eq:W}. Take $o$ from Lemma~\ref{lem:avgorder}, and bound the realization sum of $H$ by
Proposition~\ref{prop:fubini} for $H_o$:
\[
g(\{v\})\prod_ie^{G_i(\mathrm{rk}_o(N_i,i))}\le g(\{v\})\exp\sum_i\frac1{d_i+1-\varrho_i}\sum_{s=\varrho_i}^{d_i}
\log\bigl(\eta\theta^{K_0}\Sigma_{s,d_i-s+p_i}\bigr)\le g(\{v\})\prod_iA_*\lambda^{d_i}.
\]
(All costs $\eta\theta^{K_0}\Sigma_{s,j}$ are positive when $\mu>0$, since $c_{\rm r}>0$; if $\mu=0$ the family is empty
and there is nothing to prove.) The order $o$ depends only on the shape $H$, through the degrees, bits and
neighbourhoods of its vertices, and not on the realization $x$; this is what allows Proposition~\ref{prop:fubini},
whose order is fixed in advance, to be applied with it. Summing $(4A_*)^m$ over $m\ge1$ gives
$\frac\alpha{1-\alpha}g(\{v\})$.
\end{proof}

\subsection{The numerics}\label{sec:num}
From here on the constants are those of Proposition~\ref{prop:condev}, and $\ell\ge24$. By (E4) and
Lemma~\ref{lem:slot}, the cost of a vertex with $s$ earlier neighbours and cofactor exponent $j$ decays like
$2^{-64\ell\min(2j,s)}$: through the ray term when $j$ is small and through $D\ge(Q+1)^s$ otherwise. The decay is
capped at about $2^{-\ell^3}$ by the floor $c_0$, and it is absent when $s=0$ (free vertex) or $j=0$ (exact vertex).
Summed over the $d+1-\varrho$ positions of a vertex of degree $d$, the decay $\gtrsim\frac{32}3\ell d^2$ beats the cost
$\frac32\ell d^2$ of the factors $2^{K_0t_e}$, with room to pay $(2k)^d$ for the edges and the factor $4$ of
Lemma~\ref{lem:shapecount}. Appendix~\ref{app:num} turns this into the following statement.

\begin{lemma}[The averaged step cost]\label{lem:avgcostmain}
Let $0<\eta\le2$, $\theta\ge1$, $\theta^{K_0}\le\frac98$, $0\le\mu\le2^{3\ell}$ and $2k\le2^{3\ell+6}$. Then \eqref{eq:W}
holds with $A_*=2^{-7}$ and $\lambda=1/(2k)$.
\end{lemma}

\begin{proof}
This is Lemma~\ref{lem:avgcost} in Appendix~\ref{app:num}.
\end{proof}

\subsection{The gas}

\begin{proposition}[The gas]\label{prop:gasfree}
Let $r\ge0$ be coprime to $P_0$, let $\F\subseteq\F_r$ be as in Proposition~\ref{prop:condev}, $\mu=\mu(\F)\le2^{3\ell}$, and
let $k\ge1$ with $2k\le2^{3\ell+6}$. For the two-copy gas of \S\ref{sec:sieve} built on the events $\hat e_t$,
$t\in\F$,
\[
\log Z(-1,-1)\le-0.47\,\mu,\qquad\log Z(2,2)\le\sum_CW_{2,2}(C)\le9\mu .
\]
\end{proposition}

\begin{proof}
The events satisfy (E1)--(E6) by Proposition~\ref{prop:condev}, and $\mu=\mu_\Ev$ by injectivity. Take
$\tau=1/(16K_0)=1/(48\ell)$ and $\theta=e^\tau$, so $\theta^{K_0}=e^{1/16}\le\frac98$. By Lemma~\ref{lem:avgcostmain} and
Proposition~\ref{prop:sitefree}, $(1,e^\tau,\frac1{32})$ and $(2,1,\frac1{32})$ are multi-site bounds, as
$4\cdot2^{-7}=\frac1{32}$.

\emph{Prime masses.} For primes $q>Q$, Lemma~\ref{lem:Mq} and the monotonicity in $q$ give
$\mu_q\le\eps'_{Q+1}+E_0(E_0+1)/(4Y)$. Each of the five pieces, $2E_0/(Q+1)$, $E_0c_Y$, $c_1(Q+1)^{-1/2}$,
$c_2(Q+1)^{-1}$ and $E_0(E_0+1)/(4Y)$, is at most $2^{-(\ell+11)}$ for $\ell\ge24$, so $m_P=2^{-(\ell+8)}$ is admissible
in Proposition~\ref{prop:polymer}(a).

\emph{The three bounds.} $2\theta^{K_0}(2+\frac1{31})m_P\le5\cdot2^{-(\ell+8)}\le\frac1{48\ell}=\tau$, so
Proposition~\ref{prop:polymer}(a) and Corollary~\ref{cor:KP} apply, and by Proposition~\ref{prop:polymer}(b)
$\log Z(-1,-1)\le(-1+3(\theta^{K_0}-1)+4\theta^{K_0}/31)\mu\le-0.47\mu$. By Proposition~\ref{prop:polymer}(c),
$\log Z(2,2)\le\sum_CW_{2,2}(C)\le2\mu(4+\frac2{31})\le9\mu$.
\end{proof}

\section{One level and the choice of the level}\label{sec:level}

The parameters are those of \S\ref{sec:avg}: $L\ge2^{23}$, $\ell\ge24$, $K\le2^{3\ell}$. In this section and the next
we prove Theorem~\ref{thm:A}(b), with $M=2^{37}K$; \S\ref{sec:full} explains the changes for (a).

\subsection{One class and one level}

\begin{proposition}[One class]\label{prop:oneclass}
Let $r$ be coprime to $P_0$, and put
\[
X_m=\frac{L^3}{2^{37}K},\qquad T=\lfloor X_m\rfloor,\qquad A=\Bigl\lfloor\frac{2\max(T-1,0)}3\Bigr\rfloor,\qquad B=14T,\qquad k=15T.
\]
Then for every $N'\ge0$,
\[
2^A\,\#\{0\le j\le N':\ r+P_0j\ \text{not representable}\}\le2N'+2^A\bigl(2^L(L+1)^2\bigr)^{2k}.
\]
\end{proposition}

\begin{proof}
If $T=0$ then $A=k=0$ and the bound is trivial. Otherwise choose $\F\subseteq(\F_r^\flat)^\ast$ with
$T-1\le\mu(\F)\le T$, adding members one at a time, which is possible by Proposition~\ref{prop:massF} as $T\le X_m$
and each $1/v_t\le\frac12$. $\F$ satisfies the hypotheses of Proposition~\ref{prop:condev}:
its members have $v_t\ge2$ good, and it has no nested pairs, being a subset of a pruned family. Since
$X_m\le L^3<2^{3\ell}$, we have $\mu\le T\le2^{3\ell}$ and $2k=30T\le2^{3\ell+6}$, so Proposition~\ref{prop:gasfree}
applies: $\log Z(-1,-1)\le-0.47(T-1)\le-A\log2$, as $3A\le2(T-1)$ and $\log2\le0.7$, and
$\log Z(2,2)\le9T\le14T\log2=B\log2$, as $\log2\ge0.69$. Also $A+B\le k+1$.

Let $\Ev=\{\hat e_t:t\in\F\}$, so $\#\Ev=\#\F$. Proposition~\ref{prop:cond} bounds $2^A$ times the number of
$0\le j\le N'$ at which no event of $\Ev$ holds by $2N'+2^A(\#\F+1)^{2k}$. If $r+P_0j$ is not representable, then it
lies in no $E_t$ by Lemma~\ref{lem:typeII}, so no $\hat e_t$ holds at $j$ by Lemma~\ref{lem:condev}. Finally
$\#\F+1\le2^L(L+1)^2$. Indeed, a triple with $u_t\le2^L$ has $a,c\le2^L$, and for fixed $(a,c)$ there are at most
$(2^L+1)/(Kac)$ values of $d$; summing over $a,c\le2^L$ with $\sum_{a\le2^L}1/a\le1+L\log2$ and $K\ge4$ gives
$\#\F\le\frac14(2^L+1)(1+L\log2)^2\le2^L(L+1)^2-1$.
\end{proof}

\begin{proposition}[One level]\label{prop:onelevel}
Let $K\le L^3$, form $A$, $k$ as in Proposition~\ref{prop:oneclass}, put $E_{\rm th}=A+6L+3+2k(L+2\ell)$, and let $R$
be as in Lemma~\ref{lem:reduce}. For every integer $X\ge0$:
\begin{enumerate}[label=\textup{(\alph*)}]
\item $2^AR(X)\le2X\,\varphi(P_0)/P_0+2^{E_{\rm th}}$;
\item $R(X)\le X+P_0$.
\end{enumerate}
\end{proposition}

\begin{proof}
(a) $K\le L^3<2^{3\ell}$. Proposition~\ref{prop:oneclass} with $N'=\lfloor X/P_0\rfloor$, summed over the $\varphi(P_0)$
unit classes (in the formal version each class carries the extra term $2^A2^{2L+2}$, from the size condition of
Lemma~\ref{lem:typeII}; see Appendix~\ref{app:lean}), gives
\[
2^AR(X)\le2\varphi(P_0)\lfloor X/P_0\rfloor+2^A\varphi(P_0)\,\Xi,\qquad\Xi=\bigl(2^L(L+1)^2\bigr)^{2k}+2^{2L+2}.
\]
Here $\varphi(P_0)\lfloor X/P_0\rfloor\le X\varphi(P_0)/P_0$. Also $P_0\le\prod_{p\le y_0}p\le4^{y_0}\le2^{4L}$, as
$y_0=2^\ell\le2L$, and $(L+1)^2\le2^{2\ell}$, so $\Xi\le2^{2k(L+2\ell)+2L+3}$ and $2^AP_0\Xi\le2^{E_{\rm th}}$.
(b) Each of the $\varphi(P_0)\le P_0$ classes contributes at most $\lfloor X/P_0\rfloor+1$.
\end{proof}

The proposition is used only at levels $L\ge L_0$ with $2^9M\le L_0^3$ (Lemma~\ref{lem:chooseL}).

\subsection{The choice of the level}
The transfer forces $L^4\asymp M\log N$, and the sieve gains $A\asymp L^3/M$; so the constant in the exponent depends on
$M$ as $M^{-1/4}$. For an integer $M\ge1$ put $c(M)=1/(22M^{1/4})$.

\begin{lemma}\label{lem:chooseL}
Let $M=2^{37}K$, and let $L_0\ge2^{23}$ be an integer with $2^9M\le L_0^3$. Let
$h$ be an integer with $121L_0^4\le4Mh$, and put
\[
L=\bigl\lfloor\lfloor4Mh/121\rfloor^{1/4}\bigr\rfloor .
\]
Then $L\ge L_0$, $K\le L^3$, $121L^4\le4Mh<121(L+1)^4$, $2L\le h$ and $h\ge1024$. If $T,A,k$ are formed at the level $L$ as in
Proposition~\ref{prop:oneclass}, then $2L^3\le M(3A+6)$, $3MA\le2L^3$, $Mk\le15L^3$ and $E_{\rm th}\le h$.
\end{lemma}

\begin{proof}
From $L_0^4\le\lfloor4Mh/121\rfloor$ we get $L\ge L_0$, and the two bounds on $121L^4$ are the definition of $L$. Also
$K\le M\le2^9M\le L_0^3\le L^3$. As $121\cdot2^9L\cdot L^3\le121L^4\cdot2^9\le4h\cdot2^9M\le4hL^3$, we get
$121\cdot2^9L\le4h$; this gives $2L\le h$ and $h\ge1024$.

With $T=\lfloor L^3/M\rfloor$ we have $MT\le L^3<M(T+1)$. Since $A=\lfloor2(T-1)/3\rfloor$, $2T\le3A+4$ and $3A\le2T$,
which give the first two inequalities, and $k=15T$ gives the third. For $E_{\rm th}$ we use these, $2^9M\le L^3$,
$L\ge2^{23}$ and $1024\ell\le L$:
\begin{align*}
4M\,E_{\rm th}&=4MA+24ML+12M+8(L+2\ell)Mk\\
&\le3L^3+\tfrac{24}{512}L^4+\tfrac{12}{512}L^4+120L^4+\tfrac{240}{1024}L^4<121L^4\le4Mh.\qedhere
\end{align*}
\end{proof}

\begin{lemma}[The exponent]\label{lem:exponent}
In Lemma~\ref{lem:chooseL}, $A\log2\ge c(M)(h\log2)^{3/4}+\log3$.
\end{lemma}

\begin{proof}
Put $\lambda=h\log2$. As $h\ge1024$, $\lambda<(h+1)\log2\le\frac{1025}{1024}h\log2$, and $(L+1)^4\le L^4+15L^3\le\frac{122}{121}L^4$ as
$L\ge2^{23}$, so $4Mh<122L^4$ and
\[
M\lambda\le\tfrac{1025}{1024}\cdot\tfrac{0.6931471808}4\cdot122\,L^4\le21.17\,L^4 .
\]
Put $\nu=\lambda^{3/4}$ and $m=M^{1/4}$. Then $(\nu m^3)^4=(M\lambda)^3\le21.17^3L^{12}\le10^4L^{12}$, so $\nu m^3\le10L^3$
and $c(M)\nu=\nu m^3/(22M)\le\frac{10}{22}\cdot\frac{L^3}M$. On the other hand $2L^3\le M(3A+6)$ gives
$A\log2\ge(\frac23\frac{L^3}M-2)\log2$. Since $\frac23\log2-\frac{10}{22}\ge0.0075$ and $L^3/M\ge2^9$,
\[
A\log2-c(M)\nu\ge0.0075\cdot2^9-2\log2\ge2.45\ge\log3 .\qedhere
\]
\end{proof}

\section{\texorpdfstring{Divisors of exceptions: proof of Theorem~\ref{thm:A}(b)}{Divisors of exceptions: proof of Theorem 1.1(b)}}\label{sec:boot}

Lemma~\ref{lem:reduce} bounds $E_K(N)$ by a sum over the $P_0$-part $s$ of $n$. For $s\le N^{1/16}$ we use
Proposition~\ref{prop:onelevel}(a) at the level chosen for $\frac{15}{16}\log_2N$; the main terms add up to exactly
$2N$. The larger $s$ are few.

\begin{lemma}[Euler products]\label{lem:euler}
Let $P$ be a product of distinct primes and $F$ a finite set of integers $s\ge1$ with $s\mid P^\infty$. Then
\[
\sum_{s\in F}\frac1s\le\prod_{p\mid P}\frac p{p-1}=\frac P{\varphi(P)},\qquad
\sum_{s\in F}\frac1{\sqrt s}\le4^{\omega(P)}.
\]
\end{lemma}

\begin{proof}
For a completely multiplicative $f\ge0$ with $\sum_kf(p^k)<\infty$, $\sum_{s\mid P^\infty}f(s)=\prod_{p\mid P}\sum_{k\ge0}f(p^k)$
(Euler's product over a finite set of primes). With $f(s)=1/s$ the factors are $p/(p-1)$, and $\varphi(P)=\prod(p-1)$.
With $f(s)=s^{-1/2}$ they are $(1-p^{-1/2})^{-1}\le(1-2^{-1/2})^{-1}<4$.
\end{proof}

\begin{proposition}[Sum over $s$]\label{prop:sums}
In Proposition~\ref{prop:onelevel}, let $N$, $h$, $w$ be integers with $2^h\le N$,
\[
w+4L+E_{\rm th}+1\le h\qquad\text{and}\qquad A+8L+2\le w .
\]
Then $2^AE_K(N)\le3N$.
\end{proposition}

\begin{proof}
Let $Z=4^w$ and split the sum of Lemma~\ref{lem:reduce} at $s\le Z$. Note $\omega(P_0)\le y_0\le2L$, so
$4^{\omega(P_0)}\le2^{4L}$, and $P_0\le2^{4L}$.

\emph{$s\le Z$.} By Proposition~\ref{prop:onelevel}(a) with $X=\lfloor N/s\rfloor\le N/s$ and Lemma~\ref{lem:euler},
\[
2^A\sum_{s\le Z}R(\lfloor N/s\rfloor)\le2N\,\frac{\varphi(P_0)}{P_0}\sum_{s\le Z}\frac1s+\#\{s\le Z\}\,2^{E_{\rm th}}
\le2N+2^{w+4L+E_{\rm th}},
\]
since $\#\{s\le Z\}\le\sum_{s\le Z}\sqrt{Z/s}\le2^w4^{\omega(P_0)}$. The last term is at most $2^{h-1}\le N/2$.

\emph{$Z<s\le N$.} By Proposition~\ref{prop:onelevel}(b) and $\lfloor N/s\rfloor\ge1$,
$R(\lfloor N/s\rfloor)\le(1+P_0)N/s$, and $1/s\le2^{-w}s^{-1/2}$. So
\[
2^A\sum_{Z<s\le N}R(\lfloor N/s\rfloor)\le2^A(1+P_0)N\,2^{-w}4^{\omega(P_0)}\le2^{A+8L+1-w}N\le N/2 .\qedhere
\]
\end{proof}

\begin{lemma}[The base level]\label{lem:L0K}
Let $L_*$ be the least integer with $2^{46}K\le L_*^3$, and $L_0=\max(2^{23},L_*)$. If $h'\ge1$ satisfies
$121\cdot2^{53}\le Kh'$ and $2^{89}K\le h'^3$, then $L_0$ satisfies the hypotheses of Lemma~\ref{lem:chooseL} for
$h'$ and $M=2^{37}K$: $L_0\ge2^{23}$, $2^9M\le L_0^3$ and $121L_0^4\le4Mh'$.
\end{lemma}

\begin{proof}
$2^9M=2^{46}K\le L_*^3\le L_0^3$. If $L_*\le2^{23}$, then $121\cdot2^{92}=4\cdot2^{37}\cdot121\cdot2^{53}\le4Mh'$.
Otherwise $L_0=L_*=m+1$ with $m\ge2^{23}$ and $m^3<2^{46}K$. From $24(m+1)\le25m$ and $25^{12}\le2\cdot24^{12}$,
$(m+1)^{12}\le2m^{12}$, so
\[
(121L_0^4)^3\le121^3\cdot2\cdot(2^{46}K)^4=2^{117}K^3\cdot121^3\cdot2^{68}K\le2^{117}K^3\cdot2^{89}K\le(4Mh')^3,
\]
as $121^3\le2^{21}$.
\end{proof}

\begin{proof}[Proof of Theorem~\ref{thm:A}(b)]
Let $M=2^{37}K$, $h=\lfloor\log_2N\rfloor$, $h'=\lfloor15h/16\rfloor$ and $w=\lfloor h/32\rfloor$. As $h\ge4096$,
$100h'\ge93(h+1)$, so $Kh\ge2^{61}$ gives $Kh'\ge121\cdot2^{53}$ and $2^{90}K\le h^3$ gives $2^{89}K\le h'^3$.
Take $L_0$ from Lemma~\ref{lem:L0K} and $L$ from Lemma~\ref{lem:chooseL} for $h'$. Then $E_{\rm th}\le h'$.

\emph{The level is small.} From $121L^4\le4Mh$ and $2^{28}M\le2^{65}K\le h^3$ we get
$121\cdot2^{28}(265L)^4\le265^4\cdot4\,h^4\le121\cdot2^{28}h^4$, so $265L\le h$. From $3MA\le2L^3$ and $121L^4\le4Mh'$,
$363AL\le8h'$, and $L\ge2^{23}$ gives $A\le h/2^{28}$. With $h\ge4096$ these give $w+4L+E_{\rm th}+1\le h$ and
$A+8L+2\le w$, and Proposition~\ref{prop:sums} gives $2^AE_K(N)\le3N$.

\emph{The exponent.} By Lemma~\ref{lem:exponent} at $h'$, $A\log2\ge c(M)(h'\log2)^{3/4}+\log3$. Since
$\log N<(h+1)\log2$ and $100h'\ge93(h+1)$, $h'\log2\ge0.93\log N$; and $0.93^{3/4}\ge0.94$ (as $0.94^4\le0.93^3$), with
$0.94/22\ge1/24$. Hence
\[
E_K(N)\le3N\,2^{-A}\le N\exp\Bigl(-\frac{(\log N)^{3/4}}{24M^{1/4}}\Bigr).\qedhere
\]
\end{proof}

The splitting at $Z=N^{1/16}$ costs the factor $(15/16)^{3/4}$ in the exponent, which $22\to24$ absorbs. Nothing in the
proof depends on which small primes are representable, and $K$ enters only through the mass $L^3/(2^{37}K)$ and the
condition $K\le L^3$; this is why the constant is absolute in $((\log N)^3/K)^{1/4}$ and the range reaches
$K\asymp(\log N)^3$.

\section{\texorpdfstring{The whole Type~II family: proof of Theorem~\ref{thm:A}(a)}{The whole Type II family: proof of Theorem 1.1(a)}}\label{sec:full}

The family of \S\S\ref{sec:arith}--\ref{sec:mass} keeps only the triples with $a,c\le d$. This is used in the uniqueness
of a triple given its modulus and class (Lemmas~\ref{lem:oneray} and~\ref{lem:nodup}), in the diagonal bound (E6) of
Proposition~\ref{prop:condev}, and in the off-ray step of Lemma~\ref{lem:key}, where $u_0\ge2ac\max(a,c)$ needs
$d\ge\max(a,c)$. This section replaces all three. Lemma~\ref{lem:keyfull} bounds the cells of the whole family
$Y\le a,c,d$ by a slope separation inside dyadic strips. Duplicates do occur without $a,c\le d$: for $K=4$, the triples
$(65,528,2)$ and $(528,65,2)$ have $c^2-a^2=u$, so they give the same class. So the sieve keeps one triple per event,
with priority to a core $\F_0\cup W_a\cup W_c$ that has no duplicates (Lemmas~\ref{lem:dedup} and~\ref{lem:coreinj}). The
wide regions $W_a$, $W_c$ each carry at least $1.144$ times the lower bound for the mass of $\F_0$
(Proposition~\ref{prop:masscore}), and the diagonal is Lemma~\ref{lem:diagfull}. The parameters are those of
\S\ref{sec:avg}: $L\ge2^{23}$, $\ell\ge24$, $K\le2^{3\ell}$. The triples outside the core may enter the sieve as events,
but no mass is counted for them.

\subsection{The key lemma for the full family}

A triple $t=(a,c,d)$ is \emph{fully admissible} if $Y\le a$, $Y\le c$, $Y\le d$, $\gcd(a,c)=1$ and $u_t\le2^L$; no
order among $a,c,d$ is assumed. For $D\ge1$ and $B\ge0$ put
\[
\eps^{\rm full}_D(B)=(L+1)\Bigl(B+\frac{2(2L+3)}Y+\frac{L+2}2\Bigl(\frac1{Y^2}+\frac{4(2L+3)}D\Bigr)\Bigr).
\]

\begin{lemma}[Key lemma, full family]\label{lem:keyfull}
Let $D\ge1$, $\rho\in\mathbb Z$ and $B\ge0$. Let $\mathcal C$ be a set of fully admissible triples with $D\mid u_t$,
$u_t\ne D$, $D\mid a+\rho c$ and $D/u_t\le B$ for every $t\in\mathcal C$. Then
$\sum_{t\in\mathcal C}D/u_t\le\eps^{\rm full}_D(B)$.
\end{lemma}

For fixed $(a,c)$ only the least $d$ matters, as before; it costs $\min(B,D/(2acY))$ since $d\ge Y$. Two points
satisfy $D\mid ac'-a'c$, so inside a dyadic strip $c\in[C,2C)$ the slopes $a/c$ are $D/(4C^2)$ apart. Summing over a
strip costs the logarithm of the range, and there are $L+1$ strips.

\begin{proof}
\emph{(0)} Since $K\ge4$, $u_t\ge2acd$, so $a,c,d\le2^L$. Two members have $D\mid ac'-a'c$, because
$ac'-a'c=c'(a+\rho c)-c(a'+\rho c')$.

\emph{(1) Fibres.} Fix a point $p=(a,c)$ and let $G_p$ be the members with these $a,c$. Since $D\mid Kacd-1$, $Kac$ is a
unit modulo $D$, so the $d$ of $G_p$ lie in one class modulo $D$ and are distinct. Let $d^*\ge Y$ be the least. Its
term is at most $B$ and at most $D/(2acd^*)\le D/(2acY)$. The other members have $d>D$ in one class, $d\le2^L$, so
$\sum1/d\le(L+2)/D$ and their terms add up to at most $\frac D{2ac}\sum1/d\le(L+2)/(2ac)$. Hence
\[
\sum_{G_p}\frac D{u_t}\le w(p)+\frac{L+2}{2ac},\qquad w(p)=\min\Bigl(B,\frac D{2acY}\Bigr).
\]

\emph{(2) Strips.} Group the points with $c\in[C,2C)$, $C=2^i$, $0\le i\le L$, and put $\beta(p)=\lfloor4C^2a/(cD)\rfloor$.
If $\beta(p)=\beta(p')=b$, write $4C^2a=cDb+r$ and $4C^2a'=c'Db+r'$ with $0\le r<cD$, $0\le r'<c'D$. Then
$4C^2(ac'-a'c)=rc'-r'c$ and $|rc'-r'c|<cc'D<4C^2D$, so $|ac'-a'c|<D$. By (0) $ac'=a'c$, and coprimality gives $p=p'$.
So $\beta$ is injective on a strip. If $\beta(p)\ge1$, then $cD\beta(p)\le4C^2a\le4c^2a$, so $ac\ge D\beta(p)/4$. Hence a
weight with $w\le B'$ and $w\le\kappa/(ac)$ sums over a strip to at most
$B'+\sum_{1\le b\le2^{2L+2}}4\kappa/(Db)\le B'+4\kappa(2L+3)/D$. For $w(p)$ take $B'=B$, $\kappa=D/(2Y)$. For $1/(ac)$
take $B'=1/Y^2$, $\kappa=1$.

\emph{(3)} There are $L+1$ strips. Summing (1) over the points and (2) over the strips gives
$(L+1)\bigl(B+2(2L+3)/Y\bigr)+\frac{L+2}2(L+1)\bigl(1/Y^2+4(2L+3)/D\bigr)$.
\end{proof}

Compared with Lemma~\ref{lem:key}, the term $16(L+1)/\sqrt D$ disappears and the others are larger by powers of $L$:
$(L+1)B$ against one ray term, and about $4L^3/D$ against $(L+1)^3/(2D)$. No condition on $D$ is needed.

\begin{definition}\label{def:fullfam}
For $r\in\mathbb Z$, the \emph{full conditional family} $\F^{\rm full}_r$ is the set of fully admissible triples with
$s_t\le E_0$, $s_t\mid a+rc$ and $v_t\ge2$ good. Let $\F^{\rm all}_r$ be the finite set of its members with
$a,c,d\le2^L$; by step (0) above this is all of $\F^{\rm full}_r$.
\end{definition}

\begin{lemma}[Conditional key lemma, full family]\label{lem:condkeyfull}
Let $r,\rho\in\mathbb Z$, $j\ge1$ and $D\ge1$ with $\gcd(D,P_0)=1$, and let $\F\subseteq\F^{\rm full}_r$. Summing over
$t\in\F$ with $D\mid v_t$, $v_t\ne D$, $D\mid a+\rho c$ and $v_t\ge(Q+1)^jD$,
\[
\sum_t\frac D{v_t}\le E_0\,\eps^{\rm full}_D\bigl((Q+1)^{-j}\bigr).
\]
\end{lemma}

\begin{proof}
Group by $e=s_t\le E_0$. As in Lemma~\ref{lem:condkey}, the group lies in a cell at modulus $De$ with a class given by
the Chinese remainder theorem, and $De/u_t=D/v_t\le(Q+1)^{-j}$. Lemma~\ref{lem:keyfull} bounds the group by
$\eps^{\rm full}_{De}((Q+1)^{-j})\le\eps^{\rm full}_D((Q+1)^{-j})$, and there are at most $E_0$ values of $e$.
\end{proof}

\subsection{The events of a duplicate-free full family}

\begin{lemma}[Diagonal]\label{lem:diagfull}
Let $q\ge1$ and $\F\subseteq\F^{\rm full}_r$. Then $\sum_{t\in\F,\ v_t=q}1/v_t\le E_0(E_0+1)(L+1)/(KY)$.
\end{lemma}

\begin{proof}
For these $t$, $Kacd=qs_t+1$, so $t$ is determined by $(s_t,c,d)$, and $KYcd\le Kacd\le qE_0+1$. The pairs
$(c,d)\in[1,2^L]^2$ with $cd\le X=(qE_0+1)/(KY)$ number at most $\sum_{c\le2^L}X/c\le X(1+\log2^L)\le X(L+1)$. So
there are at most $E_0X(L+1)$ such triples, each of weight $1/q$, and $(qE_0+1)/q\le E_0+1$.
\end{proof}

For $\F\subseteq\F^{\rm full}_r$ let $\F^\ast$ be the members without a nested partner in $\F$, as in \S\ref{sec:avg}.

\begin{lemma}[Pruning, full family]\label{lem:prunefull}
$\mu(\F^\ast)\ge(1-\eps_F)\mu(\F)$ with $\eps_F=E_0\,\eps^{\rm full}_{Q+1}\bigl((Q+1)^{-1}\bigr)\le2^{-(\ell+6)}$.
\end{lemma}

\begin{proof}
As in Lemma~\ref{lem:condprune}, with Lemma~\ref{lem:condkeyfull} at $D=v_{t'}\ge Q+1$, $j=1$, and the monotonicity of
$\eps^{\rm full}_D$ in $D$. The numerical bound is part of Lemma~\ref{lem:numfull}.
\end{proof}

\begin{proposition}[Events of the full family]\label{prop:evfull}
Let $r\ge0$ be coprime to $P_0$ and let $\F\subseteq\F^{\rm full}_r$ be finite, without nested pairs, and with no two
members $t\ne t'$ such that $v_t=v_{t'}=v$ and $v\mid ac'-a'c$. Then $t\mapsto\hat e_t$ is injective on $\F$, and
$\{\hat e_t:t\in\F\}$ satisfies (E1)--(E6) of \S\ref{sec:avg} with
\begin{gather*}
c_0=E_0(L+1)\Bigl(\frac{2(2L+3)}Y+\frac{L+2}{2Y^2}\Bigr),\qquad c_{\rm r}=E_0(L+1),\qquad c_1=0,\\
c_2=2E_0(L+1)(L+2)(2L+3),\qquad\delta=\frac{E_0(E_0+1)(L+1)}{KY}.
\end{gather*}
\end{proposition}

\begin{proof}
(E1) and (E2) as in Proposition~\ref{prop:condev}; (E3) by assumption. (E5) and injectivity: two events with equal
moduli and a common point come from triples with $v_t=v_{t'}$ and $v\mid ac'-a'c$, which the hypothesis excludes. (E4):
by Lemma~\ref{lem:condkeyfull}, with $E_0(L+1)(Q+1)^{-j}\le2c_{\rm r}(Q+1)^{-j}$ and
$E_0(L+1)\frac{L+2}2\cdot\frac{4(2L+3)}D=c_2/D$. (E6) is Lemma~\ref{lem:diagfull}.
\end{proof}

\subsection{One triple per event}

\begin{lemma}[One triple per event]\label{lem:dedup}
Let $S$ be a finite set of triples, $P\subseteq S$, and $\kappa$ an injective map from triples to $\mathbb N$. Write
$t\approx t'$ if $v_t=v_{t'}=v$ and $v\mid ac'-a'c$. Let $\mathrm{Ded}(S,P)$ consist of the $t\in S$ with $t\in P$, or with
$t'\notin P$ and $\kappa(t)<\kappa(t')$ for every $t'\in S$, $t'\ne t$, $t'\approx t$. Then $P\subseteq\mathrm{Ded}(S,P)$,
and if no two distinct members of $P$ satisfy $\approx$, then no two distinct members of $\mathrm{Ded}(S,P)$ do.
\end{lemma}

\begin{proof}
Let $t\ne t'$ be kept with $t\approx t'$; the relation is symmetric. If both are in $P$, this contradicts the
hypothesis. If $t\in P$ and $t'\notin P$, then $t'$ was kept although $t\approx t'$ lies in $P$. If neither is in $P$,
then $\kappa(t)<\kappa(t')<\kappa(t)$.
\end{proof}

For members of $\F^{\rm full}_r$, $t\approx t'$ holds exactly when $\hat e_t=\hat e_{t'}$ as sets of $j$: by
Lemma~\ref{lem:condev} each is one class modulo $v_t$, the class $j\equiv-(rc+a)(P_0c)^{-1}$, and two such classes agree
if and only if $v\mid ac'-a'c$. So $\approx$ is an equivalence relation and $\mathrm{Ded}(S,P)$ keeps one member of each
class. Only Lemma~\ref{lem:dedup} itself is used below.

\subsection{The core}

\begin{definition}\label{def:core}
Let $\F_0=\F^\flat_r$ be the family of Proposition~\ref{prop:massF0}. Let $W_a$ be the set of $t\in\F^{\rm all}_r$ with
$E_0^3c^2\le a$ and $d<a$, and $W_c$ the set of $t\in\F^{\rm all}_r$ with $E_0^3a^2\le c$ and $d<c$. The \emph{core} is
$\F_0\cup W_a\cup W_c$.
\end{definition}

The three parts are disjoint: $\F_0$ has $a,c\le d$, and $E_0^3c^2\le a$, $E_0^3a^2\le c$ with $a,c\ge2$ would give
$a\ge c^2>c\ge a^2>a$. And $\F_0\subseteq\F^{\rm all}_r$, since its members are admissible, hence fully admissible, with
$v_t\ge2$ good.

\begin{lemma}[The core has no duplicates]\label{lem:coreinj}
Let $t,t'$ be in the core with $v_t=v_{t'}=v$ and $v\mid ac'-a'c$. Then $t=t'$. More precisely, let $t,t'$ be fully
admissible with $u_t=s_tv$, $u_{t'}=s_{t'}v$, $v\ge1$, $v\mid ac'-a'c$ and $1\le s_t,s_{t'}\le E_0$, and assume only
$K\ge4$, $Y\ge2$, $E_0\ge2$, $E_0<KY$ and $E_0<KY^2$. Then:
\begin{enumerate}[label=\textup{(\alph*)}]
\item if $E_0c\le a$ and $E_0c'\le a'$, then $t=t'$;
\item the case $E_0^3c^2\le a$, $d<a$, $a'\le d'$, $c'\le d'$ cannot occur;
\item the case $E_0^3c^2\le a$, $E_0^3a'^2\le c'$ cannot occur.
\end{enumerate}
\end{lemma}

\begin{proof}
Write $ac'-a'c=qv$, $u=sv=Kacd-1$, $u'=s'v=Ka'c'd'-1$, $1\le s,s'\le E_0$.

(a) If $v\le c^2d$, then $Kacd-1=sv\le E_0c^2d$, while $Kacd\ge4E_0c^2d$; so $c^2d<v$, and likewise $c'^2d'<v$. Now
\[
K(c'^2d'-c^2d)=Kcd\,Kc'd'(ac'-a'c)-Kv\,(sc'^2d'-s'c^2d),
\]
and $\gcd(v,K)=1$, so $v\mid c'^2d'-c^2d$, and with both terms in $[0,v)$ we get $c^2d=c'^2d'=:n_1$. Then
$Kan_1=c(sv+1)$ and $Ka'n_1=c'(s'v+1)$ give $Kn_1q=cc'(s-s')$. If $q\ne0$, then $Kn_1\le cc'(E_0-1)$, while
$n_1^2=c^2d\,c'^2d'\ge(cc'Y)^2$ gives $Kn_1\ge Kcc'Y>cc'E_0$. So $q=0$, $ac'=a'c$, coprimality gives $(a,c)=(a',c')$, and
$n_1$ gives $d=d'$.

(b) Note $a\ge E_0^3c^2\ge c$. \emph{If $q=0$}, then $(a,c)=(a',c')$ and $v\mid Kac(d-d')$; as $\gcd(v,Kac)=1$, $d\ne d'$
would give $|s-s'|v=Kac|d-d'|\ge Kacv$ and $E_0>|s-s'|\ge KY^2$; and $d=d'$ contradicts $d<a=a'\le d'$. \emph{If $q\ge1$},
then $v<ac'$. From $s'v\ge KY^2c'-1>E_0c'$ we get $c'\le v$, so $q<a$. From $Ka'^2c'-1\le Ka'c'd'-1=s'v\le E_0v<E_0ac'$ we get
$Ka'^2\le E_0a$, hence $K(sa'c)^2\le E_0^2c^2\cdot E_0a\le a^2<Ka^2$ and $0<sa'c<a$. The identity
\[
sa'c-q=a\,(sc'-qKcd)
\]
(from $s\cdot qv=sac'-sa'c$ and $q\cdot sv=qKacd-q$) and $|sa'c-q|<a$ give $sa'c=q$ and $sc'=qKcd$, so $c'=Ka'c^2d$ and
$a'\mid\gcd(a',c')=1$, against $a'\ge Y$. \emph{If $q\le-1$}, then $v<a'c$, so $Kacd\le E_0a'c$ and $Ka'c'd'\le E_0a'c$,
that is $Kad\le E_0a'$ and $KYa'\le Kc'd'\le E_0c$. Hence $K^2Y^2a\le KY\,Kad\le E_0KYa'\le E_0^2c\le E_0^2a$, against
$E_0<KY$.

(c) Here $a\ge2c$ and $c'\ge2a'$, so $ac'>cc'>ca'$ and $q\ge1$. As in (b), $c'\le v$ and $q<a$. Symmetrically,
$sv\ge KY^2a-1\ge(E_0+1)a-1\ge E_0a$ gives $a\le v$, and $qv<ac'\le vc'$ gives $q<c'$. If $sa'c<a$, the identity of (b)
gives $c'=Ka'c^2d$, a contradiction. If $s'a'c<c'$, the identity $s'a'c-q=c'(s'a-qKa'd')$ gives $a=Ka'^2cd'$ and
$c\mid\gcd(a,c)$, a contradiction. Otherwise $E_0a'c\ge sa'c\ge a\ge E_0^3c^2$ and $E_0a'c\ge c'\ge E_0^3a'^2$, so
$a'\ge E_0^2c$ and $c\ge E_0^2a'$, which is impossible for $E_0\ge2$.

\emph{The core.} Two members of $W_a$ satisfy (a), as $E_0c\le E_0^3c^2$. A member of $W_a$ and one of $\F_0$ satisfy (b).
A member of $W_a$ and one of $W_c$ satisfy (c). The cases with $W_c$ follow by the swap $(a,c,d)\mapsto(c,a,d)$, which
keeps $u$, $s$, $v$ and changes only the sign of $ac'-a'c$. Two members of $\F_0$ are equal by Lemma~\ref{lem:nodup}.
At our parameters $E_0=2^{5\ell}\ge2$, $2E_0<KY$ and $E_0<KY^2$.
\end{proof}

The quadratic condition is needed for the arithmetic of (c). For $K=4$, $Y=2$ and $E_0=3$, the triples $(11,2,4)$ and
$(2,11,4)$ have $u=351=3\cdot117$ and $117\mid11^2-2^2$. With $s=s'=3$ and $v=117$ they satisfy the hypotheses of the
arithmetic in (c) with $E_0c\le a$, $E_0a'\le c'$ in place of the quadratic conditions. ($v=117$ is not good, so they
are not members of the family.)

\subsection{The mass of the wide regions}

Let $T_0=J_0^\flat+10\ell+64$, with $J_0^\flat$ as in \S\ref{sec:brun}, and $H=3\ell+4$. A \emph{wide scale} is a
triple $(i',j,i)$ with
\[
\ell^3\le i'<L,\qquad\ell^3\le j<L,\qquad 3i'+j+T_0\le i<L-H-i'-j .
\]
Its box consists of the triples with $a\in[2^i,2^{i+1})$, $c\in[2^{i'},2^{i'+1})$, $d\in[2^j,2^{j+1})$. The lower edge
of $i$ is what a sieve in $a$ needs: along a progression $a=a_0+ecm$ the cofactor $(Kacd-1)/e$ is an affine form in
$m$ with leading coefficient $Kc^2d$, and the range of $m$ must be long compared with it.

\begin{lemma}[Count of wide scales]\label{lem:countW}
There are at least $L^3/60$ wide scales.
\end{lemma}

\begin{proof}
Put $V_1=L-H-T_0$ and $U_1=V_1-6\ell^3$. Since $2^{\ell-1}\le L$ and $J_0^\flat\le128\ell^2+5506\ell+73$, we have
$40(H+T_0)\le40(128\ell^2+5519\ell+141)\le2^{\ell-1}\le L$ and $512\ell^3\le2^{\ell-1}\le L$ for $\ell\ge24$ (at
$\ell=24$: $40(H+T_0)=8{,}253{,}000\le2^{23}$); these give $U_1\ge0.9632L$. For fixed $i'$,
$j$ the number of $i$ is at least $V_1-4i'-2j$. Restrict to $i'=\ell^3+k$, $0\le k<n_1=\lfloor(U_1-1)/4\rfloor$, and
$j=\ell^3+m$, $0\le m<\lfloor x/2\rfloor$, where $x=U_1-4k\ge1$. The inner sum over $j$ is at least
$\lfloor x/2\rfloor\bigl(x-\lfloor x/2\rfloor\bigr)\ge(x^2-x)/4$. With $f(y)=y^3/12-y^2/8$ we have
$f(y)-f(y-4)=y^2-5y+\frac{22}3\le y^2-y$ for $y\ge\frac{11}6$, and $x\ge U_1-4(n_1-1)\ge5$ on the range, so the outer
sum telescopes:
\[
\sum_{k<n_1}\frac{(U_1-4k)^2-(U_1-4k)}4\ge\frac{f(U_1)-f(U_1-4n_1)}4\ge\frac14\Bigl(\frac{U_1^3}{12}-\frac{U_1^2}8-\frac{16}3\Bigr),
\]
as $1\le U_1-4n_1\le4$. With $U_1\ge0.9632L$ and $L\ge2^{23}$ this is at least $L^3/60$.
\end{proof}

The exact number of wide scales is $1.125\cdot L^3/60$ at $L=2^{23}$ and tends to $L^3/48$.

\begin{lemma}[Coprime pairs in a box]\label{lem:pairsbox}
Let $B,e\ge1$ be integers with $64e(5+3\log2B)\le B$. Then
\[
\sum_{\substack{1\le x<B,\ B\le c<2B\\ \gcd(x,c)=\gcd(c,e)=1}}\frac1{c^2}\ \ge\ \frac{31}{256}\cdot\frac{\varphi(e)}e .
\]
\end{lemma}

\begin{proof}
\emph{Pairs with $\gcd(c,e)=1$.} Any $B$ consecutive integers contain at least $\varphi(e)\lfloor B/e\rfloor$ integers
prime to $e$. So there are at least $(B-1)\varphi(e)\lfloor B/e\rfloor\ge\varphi(e)(B-1)(B/e-1)$ pairs $(x,c)$ with
$1\le x<B$, $B\le c<2B$ and $\gcd(c,e)=1$.

\emph{Non-coprime pairs.} If $\gcd(x,c)>1$, some prime $g\le2B$ divides $x$ and $c$. For fixed $g$ there are at most
$B/g$ such $x$. The $c=gy$ with $B\le c<2B$ and $\gcd(c,e)=1$ have $y$ prime to $e$ in an interval of length
at most $B/g+2$, so there are at most $\varphi(e)(\lfloor(\lfloor B/g\rfloor+2)/e\rfloor+1)\le\varphi(e)(B/(ge)+3)$ of
them. Hence at most $\varphi(e)(B^2/(g^2e)+3B/g)$ pairs have $g\mid\gcd(x,c)$. Summing over the primes $g\le2B$,
with $\sum_gg^{-2}\le\frac12$ (Lemma~\ref{lem:mertens}(d)) and $\sum_{g\le2B}1/g\le1+\log2B$, at most
$\varphi(e)(B^2/(2e)+3B(1+\log2B))$ pairs are not coprime.

\emph{Conclusion.} The coprime pairs number at least $\varphi(e)$ times
\[
(B-1)\Bigl(\frac Be-1\Bigr)-\frac{B^2}{2e}-3B(1+\log2B)\ \ge\ \frac{B^2}{2e}-B(5+3\log2B)\ \ge\ \frac{31}{64}\cdot\frac{B^2}e,
\]
using $(B-1)(B/e-1)\ge B^2/e-2B$ and the hypothesis. Each term of the sum is at least $1/(4B^2)$.
\end{proof}

\begin{proposition}[Mass of the core]\label{prop:masscore}
Let $r$ be coprime to $P_0$, and let
\[
\Phi=\frac{\lambda_KL^3}{1382400K\ell}\sum_{e\mid P_0,\,e\le E_0}\frac{\varphi(e)}{e^2}
\]
be the lower bound of Proposition~\ref{prop:massF0}. Then
\[
\mu(W_a)\ge1.144\,\Phi,\qquad\mu(W_c)\ge1.144\,\Phi,\qquad\mu(\F_0\cup W_a\cup W_c)\ge\tfrac{13}4\Phi .
\]
\end{proposition}

\begin{proof}
\emph{One scale, one $e$.} Fix a wide scale, write $B=2^{i'}$, $D=2^j$, $X=2^i$, and fix $e\mid P_0$, $e\le E_0$. Take the
triples $(x,c,d)$ with $1\le x<B$, $c\in[B,2B)$, $\gcd(x,c)=\gcd(c,e)=1$, $d\in[D,2D)$ and $e\mid Krc^2d+1$; and for each,
the $a\in[X,2X)$ with $a\equiv x\pmod c$, $e\mid Kacd-1$ and $(Kacd-1)/e$ good. These $a$ form one progression
$a=a_0+ecm$, whose cofactor is $Kc^2d\,m+(Kcda_0-1)/e$, affine with coprime coefficients. Brun's sieve
(Corollary~\ref{cor:good}), with the conditions on $X$ of a wide scale, gives $\sum e/(Kacd-1)\ge\lambda_K/(1200Kc^2d\ell)$
over the good $a$. Then $e\mid Kacd-1$ and $e\mid Krc^2d+1$ give $e\mid a+rc$ (as $\gcd(K,e)=\gcd(c,e)=1$, and $e\mid Krc^2d+1$ forces
$\gcd(d,e)=1$), $\gcd(a,c)=1$ by the class of $x$, $s_t=e$ by
goodness, and $t\in W_a$, since $a\ge X\ge E_0^3(2B)^2$ and $a\ge X>2D$ follow from $i\ge3i'+j+T_0$. For fixed $c$ prime to
$e$, the $d\in[D,2D)$ with $e\mid Krc^2d+1$ form one class modulo $e$, so there are at least $D/e-1$ of them, each with
$1/d>1/(2D)$; as $128e\le D$, $\sum_d1/d\ge\frac1{2e}-\frac1{2D}\ge\frac{63}{128}\cdot\frac1e$. The pairs $(x,c)$ have
$\sum1/c^2\ge\frac{31}{256}\cdot\frac{\varphi(e)}e$ by Lemma~\ref{lem:pairsbox}, as $64e(5+3\log2B)\le B$ (since $i'\ge\ell^3$ and $e\le2^{5\ell}$). Together, one scale and one $e$ carry at least
$\frac{\lambda_K}{1200K\ell}\cdot\frac{1953}{32768}\cdot\frac{\varphi(e)}{e^2}$.

\emph{All scales.} Different scales and different $e$ give disjoint sets, so by Lemma~\ref{lem:countW}
\[
\mu(W_a)\ge\frac{L^3}{60}\cdot\frac{\lambda_K}{1200K\ell}\cdot\frac{1953}{32768}\sum_e\frac{\varphi(e)}{e^2}
=\frac{1382400\cdot1953}{60\cdot1200\cdot32768}\,\Phi\ge1.144\,\Phi .
\]
\emph{The mirror.} The swap $(a,c,d)\mapsto(c,a,d)$ of the same construction, run with the class $r_e^\ast$ inverse to
$r$ modulo $e$ for each $e$, lands in $W_c$: $e\mid a+r_e^\ast c$ becomes $e\mid c+ra$ after the swap, and $u$, $s$, $v$
are unchanged. The three parts are disjoint, and $1+2\cdot1.144=3.288\ge\frac{13}4$.
\end{proof}

\subsection{\texorpdfstring{The numerics with per-step exponent $11\ell$}{The numerics with per-step exponent 11l}}
The constants of Proposition~\ref{prop:evfull} are larger than those of Proposition~\ref{prop:condev}: $c_2$ is about
$2^{8\ell}$. Lemmas~\ref{lem:stepcost} and~\ref{lem:avgcost} hold with $5\ell$ replaced by $11\ell$. (About $8\ell$ would
suffice; the formal proof uses $11\ell$.)

\begin{lemma}[Numerics, full family]\label{lem:numfull}
Let $\ell\ge24$ and take the constants of Proposition~\ref{prop:evfull}.
\begin{enumerate}[label=\textup{(\alph*)}]
\item $c_02^\Gamma\le2^{5\ell}\le2^{11\ell+3}$ and $2c_{\rm r}+c_1+c_2\le2^{11\ell+3}$.
\item If $0<\eta\theta^{K_0}\le4$ and $0\le\mu\le2^{3\ell}$, then for $0\le s\le d\le K_0$, $p\in\{0,1\}$, $j=d-s+p$,
\[
\log_2\bigl(\eta\theta^{K_0}\Sigma_{s,j}\bigr)\le3\ell s+11\ell+7-\min\bigl(64\ell\min(2j,s),\Gamma\bigr).
\]
\item The average of the right side of (b) over $\varrho\le s\le d$ is at most $-7-d(3\ell+6)$, for $1\le d\le3\ell$,
$\varrho\le\min(1,d)$ and $p=1$ if $d=1$. Hence \eqref{eq:W} holds with $A_*=2^{-7}$, $\lambda=1/(2k)$ under the
hypotheses of Lemma~\ref{lem:avgcostmain}.
\item The prime masses satisfy $\mu_q\le2^{-(\ell+6)}$ for $q>Q$, and $\eps_F\le2^{-(\ell+6)}$.
\end{enumerate}
Here $\Gamma=\ell^3-2\ell-7$, as in Appendix~\ref{app:num}.
\end{lemma}

\begin{proof}
Let $x=2^\ell\ge L+1$. (a) $c_0\le E_0x(6x/Y+x/Y)=7x^7/Y$ and $Y=2^\Gamma x^2\,2^7$, so $c_02^\Gamma\le7x^5/128$. Also
$c_{\rm r}\le x^6$ and $c_2\le x^5\cdot x\cdot12x^2=12x^8$, and $2x^6+12x^8\le8x^{11}$.

(b) As Lemma~\ref{lem:stepcost}, with $2^{5\ell+3}$ replaced by $2^{11\ell+3}$ in the two floors.

(c) As in Lemma~\ref{lem:avgcost} it suffices that
$3\ell\frac{d(d+1)}2+(d+1)\bigl(11\ell+14+d(3\ell+6)\bigr)\le S$. For $d=1,2$ the term $s=1$ gives $S\ge64\ell$, against
$31\ell+40$ and $60\ell+78$. For $d=3$ the terms $s=1,2$ give $S\ge192\ell$ (as $128\ell\le\Gamma$), against $98\ell+128$. For $d\ge4$, the bound
$S\ge\frac{64\ell(d^2-1)}6-\frac{35d^3}\ell$ of Lemma~\ref{lem:avgcost} reduces the claim to
\[
\tfrac{37}6\ell^2d^2-\tfrac{31}2\ell^2d-\tfrac{65}3\ell^2-6\ell d^2-20\ell d-14\ell-35d^3\ge0 .
\]
For $4\le d\le11$ this is checked for each $d$ (it is a quadratic in $\ell$, positive at $\ell\ge24$). For $d\ge12$, use
$d^3\le3\ell d^2$: the left side is at least $\ell\bigl(24(\frac{37}6d^2-\frac{31}2d-\frac{65}3)-111d^2-20d-14\bigr)=
\ell(37d^2-392d-534)>0$.

(d) By Lemma~\ref{lem:Mq}, $\mu_q\le c_0+2c_{\rm r}(Q+1)^{-1}+c_2(Q+1)^{-1}+\delta$ for $q\ge Q+1$. Each of
$2c_{\rm r}/(Q+1)\le2x^6/2^{128\ell}$, $c_0\le7x^7/Y$, $c_2/(Q+1)\le12x^8/2^{128\ell}$ and $\delta\le x^{11}/Y$ is at most
$12x^{11}/x^{14}\le1/(256x)=2^{-(\ell+8)}$, as $Y,Q+1\ge x^{14}$, so $\mu_q\le2^{-(\ell+6)}$. Likewise
$\eps_F=c_0+E_0(L+1)/(Q+1)+c_2/(Q+1)\le2^{-(\ell+6)}$.
\end{proof}

With these, Proposition~\ref{prop:gasfree} holds for every finite $\F$ as in Proposition~\ref{prop:evfull} with
$\mu(\F)\le2^{3\ell}$, with the same proof: $m_P=2^{-(\ell+6)}$ satisfies $2\theta^{K_0}(2+\frac1{31})m_P\le5\cdot2^{-(\ell+6)}\le
\frac1{48\ell}$.

\subsection{\texorpdfstring{Proof of Theorem~\ref{thm:A}(a)}{Proof of Theorem 1.1(a)}}\label{sec:proofA}

\begin{proposition}[Mass, full family]\label{prop:massfull}
Let $r$ be coprime to $P_0$, let $\mathcal D=\mathrm{Ded}(\F^{\rm all}_r,\F_0\cup W_a\cup W_c)$ with an injective key
$\kappa$, and let $\mathcal D^\ast$ be its pruning. Then $\mathcal D^\ast$ satisfies the hypotheses of
Proposition~\ref{prop:evfull}, and
\[
\mu(\mathcal D^\ast)\ge\frac{L^3}{2^{32}K}.
\]
\end{proposition}

\begin{proof}
By Lemmas~\ref{lem:dedup} and~\ref{lem:coreinj}, $\mathcal D$ has no two members with $t\approx t'$, and neither has
$\mathcal D^\ast\subseteq\mathcal D$; it has no nested pairs by construction. The core lies in $\mathcal D$, so by
Lemma~\ref{lem:prunefull}, Proposition~\ref{prop:masscore} and Lemma~\ref{lem:gain} (as in Proposition~\ref{prop:massF}),
\[
\mu(\mathcal D^\ast)\ge(1-2^{-30})\cdot\frac{13}4\cdot\frac{L^3}{1382400\,K}\cdot\frac{e^{-7}\log2}6\ge\frac{L^3}{2^{32}K},
\]
using $\eps_F\le2^{-(\ell+6)}\le2^{-30}$, $e^7\le1097$ (from $e<2.7182818286$) and $\log2\ge0.6931471803$:
$(1-2^{-30})\frac{13}4\cdot\frac{0.6931471803}{1097\cdot8294400}\ge2.475\cdot10^{-10}\ge2^{-32}$.
\end{proof}

\emph{One level.} Propositions~\ref{prop:oneclass} and~\ref{prop:onelevel} hold with $X_m=L^3/(2^{32}K)$, with the same
proofs: choose $\F\subseteq\mathcal D^\ast$ with $T-1\le\mu(\F)\le T$, which satisfies the hypotheses of
Proposition~\ref{prop:evfull} as a subset, and apply the gas bounds of Lemma~\ref{lem:numfull}. The count
$\#\F+1\le2^L(L+1)^2$ holds as before, since $a,c,d\le2^L$ and $u_t\le2^L$. Lemmas~\ref{lem:chooseL}
and~\ref{lem:exponent} and Proposition~\ref{prop:sums} hold for $M=2^{32}K$ with the same proofs, which use only $K\le M$,
$2^9M\le L_0^3$ and $L\ge2^{23}$.

\begin{lemma}[The base level, $M=2^{32}K$]\label{lem:L0Z}
Let $L_*$ be the least integer with $2^{41}K\le L_*^3$, and $L_0=\max(2^{23},L_*)$. If $h'\ge1$ satisfies
$121\cdot2^{58}\le Kh'$ and $2^{84}K\le h'^3$, then $L_0\ge2^{23}$, $2^9M\le L_0^3$ and $121L_0^4\le4Mh'$ for $M=2^{32}K$.
\end{lemma}

\begin{proof}
As Lemma~\ref{lem:L0K}: if $L_*\le2^{23}$, then $121\cdot2^{92}=4\cdot2^{32}\cdot121\cdot2^{58}\le4Mh'$. Otherwise
$L_0=m+1$ with $m\ge2^{23}$ and $m^3<2^{41}K$, so $(m+1)^{12}\le2m^{12}$, and, as $121^3\le2^{21}$, $(121L_0^4)^3\le121^3\cdot2\cdot(2^{41}K)^4=2^{102}K^3\cdot121^3\cdot2^{63}K\le
2^{102}K^3\cdot2^{84}K\le(4Mh')^3$.
\end{proof}

\begin{proof}[Proof of Theorem~\ref{thm:A}(a)]
As the proof of Theorem~\ref{thm:A}(b) in \S\ref{sec:boot}, with $M=2^{32}K$ and Lemma~\ref{lem:L0Z} in place of
Lemma~\ref{lem:L0K}. As $100h'\ge93(h+1)$, $Kh\ge2^{66}$ gives $Kh'\ge121\cdot2^{58}$ and $2^{85}K\le h^3$ gives
$2^{84}K\le h'^3$; and $2^{28}M=2^{60}K\le h^3$.
\end{proof}

Dividing $M$ by $32$ multiplies the linear threshold by $32$ ($2^{61}\to2^{66}$), through $121L_0^4\le4Mh'$ at
$L_0=2^{23}$; this is why Theorem~\ref{thm:A}(a) does not contain (b).

\section{Discussion}\label{sec:disc}

\emph{Where $\frac34$ comes from.} The following heuristic explains the exponent. Two rates fix it. The \emph{mass}
of the Type~II classes with moduli up to $2^L$ grows like $L^3/K$, one logarithm for each free scale $a,c,d$, and
Proposition~\ref{prop:massF} attains this order. Elsholtz and Tao \cite[Theorem~1.1]{ElsholtzTao} show that the average
number of solutions per prime lies between $(\log p)^3$ and $(\log p)^3\log\log p$ (for $K=4$), so no system of Type~I
or~II classes grows much faster. (In the terminology of \cite{ElsholtzTao}, these are the solutions in which $n$
divides exactly one, respectively exactly two, of $x,y,z$.) The \emph{transfer} (Lemma~\ref{lem:transfer}) evaluates intersections of
up to $2k\asymp\mu$ events, and an event of modulus near $2^j$ uses about $j$ of the $\log_2N$ available bits. Scale
$j$ therefore brings mass $\asymp j^2/K$ per unit of $j$ at a cost of $j$ bits per event. The bits are best spent on all
scales $j\le L$ with $\int_0^Lj^3\,dj\asymp K\log N$, giving $A\asymp\mu\asymp L^3/K\asymp((\log N)^3/K)^{1/4}$. In
general, mass $L^\beta$ at a cost of $L^\gamma$ bits per event gives the exponent $\beta/(\beta+\gamma)$, with $\beta=3$
and $\gamma=1$ here. Heuristically, Selberg's sieve and the large sieve follow the same law with $\gamma=1$, since
$(\log N)^{\beta/(\beta+1)}$ is the Tauberian exponent of $\prod_t(1+u_t^{-s})\approx\exp(c(s-1)^{-\beta})$. Vaughan's
system has $\beta=2$, whence the large-sieve ceiling $\frac23$ of \cite{Dahan}. The large sieve sifts by classes modulo
primes, and the Type~II classes with prime modulus carry mass only $\asymp L^2/K$ (for $K=4$, see
\cite[Proposition~4.32]{Dahan}). Those with squarefree rough composite moduli carry, after averaging over the small
primes, $\asymp L^3/K$, but using them requires control of the correlations between
classes that share prime factors; this is the role of the polymer gas (\S\S\ref{sec:KP} and~\ref{sec:free}). The
exponent $\frac57$ announced in \cite{Dahan57} also comes from composite witnesses.

\emph{A barrier.} In a model with independent events, square-sieve weights of the kind used here cannot improve the
exponent.

\begin{proposition}\label{prop:barrier}
Let $(E_t)$ be independent events, one for each Type~II triple $t=(a,c,d)$, with $\mathbf P(E_t)=1/u_t$,
$u_t=Kacd-1$, where $K\ge4$, and let $K\log\xi\ge48$. If $\Psi$ is a linear combination of the products $\prod_{t\in T}\mathbf 1_{E_t}$
over finite sets $T$ with $\prod_{t\in T}u_t\le\xi$, and $\Psi=1$ where no event occurs, then
\[
\mathbf E[\Psi^2]\ \ge\ \exp\Bigl(-\tfrac43\bigl(\tfrac{48}K\bigr)^{1/4}(\log\xi)^{3/4}\Bigr).
\]
\end{proposition}

\begin{proof}
Only the finitely many $t$ with $u_t\le\xi$ enter $\Psi$. With $p_t=1/u_t$, the functions
$\phi_T=\prod_{t\in T}(\mathbf 1_{E_t}-p_t)/\sqrt{p_t(1-p_t)}$ are orthonormal and span the same space as the products,
since the family of sets $T$ is closed under subsets. Where no event occurs, $\phi_T^2=\prod_{t\in T}g_t$ with
$g_t=1/(u_t-1)$, so Cauchy--Schwarz gives $1\le\mathbf E[\Psi^2]\,G$ with $G=\sum_T\prod_{t\in T}g_t$. For $\sigma\ge0$,
$G\le\xi^\sigma\prod_t(1+g_tu_t^{-\sigma})\le\exp(\sigma\log\xi+\sum_tg_tu_t^{-\sigma})$. Since $Kacd\ge4$,
$g_tu_t^{-\sigma}\le\frac2K(acd)^{-1-\sigma}$, and $\sum_{n\ge1}n^{-1-\sigma}\le1+\frac1\sigma\le\frac2\sigma$ for
$0<\sigma\le1$; hence $\sum_tg_tu_t^{-\sigma}\le16/(K\sigma^3)$. Take $\sigma=(48/(K\log\xi))^{1/4}\le1$; then
$16/(K\sigma^3)=\frac13\sigma\log\xi$.
\end{proof}

\begin{remark}
As $\xi\to\infty$ the constant $\frac43(48/K)^{1/4}$ improves to $(\frac43+o(1))(3/K)^{1/4}$: in the proof,
$\sum_tg_tu_t^{-\sigma}=(1+o(1))K^{-1}\zeta(1+\sigma)^3=(1+o(1))K^{-1}\sigma^{-3}$ as $\sigma\to0$, and we take
$\sigma=(3/(K\log\xi))^{1/4}$. (This asymptotic statement is not part of the formal proof.)
\end{remark}

The weight $F_k^2$ of the proof has this form with $\xi=2^{kL}\le\sqrt N$ (the choice of the level
gives $2kL\le E_{\rm th}\le\log_2N$), and so do Selberg's $\Lambda^2$ weights; heuristically, by duality, so does the
large sieve. For such weights the main term is $\mathbf E[\Psi^2]$ even when the events are exactly independent. So in
the independent model of the Type~II classes, no such weight proves Theorem~\ref{thm:A} with a constant larger than $\frac43(\frac{48}K)^{1/4}2^{-3/4}\approx2.09\,K^{-1/4}$
(asymptotically $(\frac43+o(1))(\frac3K)^{1/4}2^{-3/4}\approx1.04\,K^{-1/4}$), or with an exponent above $\frac34$. With
all Type~I and~II classes one expects, by the count of \cite{ElsholtzTao}, the same up to $(\log N)^{o(1)}$. This
suggests that a larger exponent needs information about the exceptional set beyond the distribution of $n$ in residue
classes of modulus at most $N$. In the heuristic of \cite[\S1]{PW} a prime $p$ is exceptional with probability
$\exp(-(\log p)^3/K)$. Our exponent is a constant times the fourth root of the exponent $(\log N)^3/K$ of this heuristic,
because the sieve works at the level $L$ rather than at $\log_2N$.

\emph{What is avoided.} The sieve is combinatorial, and we never show that the expansion converges: Dobrushin's
induction gives the real upper bound directly. The mass needs only Brun's pure sieve and Mertens-type bounds
\cite{Mertens}, because $\log Q\asymp\log L$ against scales $\asymp L$. The proof uses neither the fundamental lemma of
sieve theory \cite{HalberstamRichert,FriedlanderIwaniec} nor Shiu's theorem \cite{Shiu}, and everything is effective.
The random order keeps $\log Q\asymp\ell$, so that roughness and averaging cancel up to a constant factor; a
step-by-step cluster bound would lose a power of $\log\log N$.

\emph{Constants.} The $22$ in $c(M)$ comes from the gas bounds $0.47\mu$ and $9\mu$ through $A\approx\frac23T$, $k=15T$,
and the divisor reduction turns it into $24$. The $2^{32}$ in $M$ collects pruning, the count of scales, coprime pairs,
the sieve fraction and the divisor sum. Against the $2^{37}$ of (b), the whole family gains a factor $\frac{13}4$, and
sharper values of $e^7$, $\log2$ and the pruning loss a factor $10.35$. None is optimized. Theorem~\ref{thm:A}(a) has
$1/c_K=24\cdot2^8K^{1/4}=6144K^{1/4}$, against the limit $1/c\approx K^{1/4}/2.09$ of Proposition~\ref{prop:barrier} at
$\xi=\sqrt N$. The power of $K$ is the right one, but the absolute constant is not.

\emph{Questions.} Which arithmetic information beyond residue classes of modulus at most $N$ could pass the exponent
$\frac34$? Theorem~\ref{thm:A} is non-trivial for $K$ up to a constant times $(\log N)^3$, while
$E_K(N)\gg N/\log N$ for $K\ge(\log N)^3(\log\log N)^{2+o(1)}$ \cite{PW}. What happens in between? Can the thresholds be
brought to a computationally relevant size?

\appendix

\section{Elementary prime estimates}\label{app:primes}
\begin{lemma}\label{lem:mertens}
For integers $n,z\ge2$:
\begin{enumerate}[label=\textup{(\alph*)}]
\item $\sum_{p\le n}\log p\le n\log4$ and $\sum_{p\le n}\frac{\log p}p\le\log n+\log4$;
\item $\sum_{p\le z}\frac1p\le\log\log z-\log\log2+3+\frac{\log4}{\log z}$;
\item if $z\ge101$, then $\prod_{p\le z}(1-1/p)\ge\frac1{3.3\log z}$; if $\log z\ge139$, then
$\prod_{p\le z}(1+2/p)\le412\bigl(\frac{\log z}{\log2}\bigr)^2$;
\item $\sum_{p\ne2}p^{-2}\le\frac19+\frac18$, so $\sum_pp^{-2}\le\frac14+\frac19+\frac18$.
\end{enumerate}
\end{lemma}

\begin{proof}
(a) The first bound is $\prod_{p\le n}p\le4^n$. For the second,
\[
\log n!\ge\sum_{p\le n}\Bigl\lfloor\frac np\Bigr\rfloor\log p\ge n\sum_{p\le n}\frac{\log p}{p}-\sum_{p\le n}\log p,
\]
and $\log n!\le n\log n$.

\emph{Abel summation.} For integers $2\le m\le z$, reals $a_n$, $A(n)=\sum_{m\le k\le n}a_k$ and $f$ positive and
non-increasing,
\begin{equation}\label{eq:abel}
\sum_{m\le n\le z}a_nf(n)=A(z)f(z)+\sum_{n=m}^{z-1}A(n)\bigl(f(n)-f(n+1)\bigr),
\end{equation}
and the coefficients $f(n)-f(n+1)$ are nonnegative, so an upper (lower) bound for every $A(n)$ gives an upper (lower)
bound for the sum. We use it with $f(n)=1/\log n$ and $a_n=\log n/n$ at primes $n$, $a_n=0$ otherwise. Put
$\delta_n=\log\log(n+1)-\log\log n$; with $t=\log(n+1)/\log n$, the inequalities $\log t-\frac12(\log t)^2\le1-1/t\le\log t$
give
\begin{equation}\label{eq:delta}
\delta_n-\tfrac12\delta_n^2\ \le\ 1-\frac{\log n}{\log(n+1)}\ \le\ \delta_n,\qquad\sum_{n=m}^{z-1}\delta_n=\log\frac{\log z}{\log m}.
\end{equation}

(b) Take $m=2$ in \eqref{eq:abel}, and $A(n)\le\log n+\log4$ from (a). Then
\[
\sum_{p\le z}\frac1p\le1+\frac{\log4}{\log z}+\sum_{n=2}^{z-1}\Bigl(1-\frac{\log n}{\log(n+1)}\Bigr)
+\log4\Bigl(\frac1{\log2}-\frac1{\log z}\Bigr)\le1+\log\frac{\log z}{\log2}+2,
\]
by \eqref{eq:delta} and $\log4=2\log2$; this implies (b).

(c) We treat the primes up to $100$ exactly. By direct computation over the $25$ primes,
$P_0'=\prod_{p\le100}(1-\frac1p)\ge0.1203$ and $S_0=\sum_{p\le100}\frac{\log p}p\ge3.36$; for the latter,
each $\log p\ge\frac sr\log2$ follows from an integer inequality $2^s\le p^r$ with $r\le12$. For the primes
in $(100,z]$ take $m=101$ in \eqref{eq:abel}, with $A(n)=\sum_{100<p\le n}\frac{\log p}p\le\log n+c_*$ by (a), where
$c_*=\log4-S_0\le-1.9737$. As in (b),
\[
\sum_{100<p\le z}\frac1p\le1+\frac{c_*}{\log z}+\log\frac{\log z}{\log101}+c_*\Bigl(\frac1{\log101}-\frac1{\log z}\Bigr)
=\log\frac{\log z}{\log101}+1+\frac{c_*}{\log101}\le\log\frac{\log z}{\log101}+0.573,
\]
since $c_*<0$ and $\log101\le\frac{20}3\log2$. The Taylor remainder
$|\log(1-x)+x+\frac{x^2}2|\le\frac{x^3}{1-x}$ gives $\log(1-x)\ge-x-x^2$ for $0\le x\le\frac13$, and
$\sum_{p>100}p^{-2}\le\frac1{100}$, so
\[
\prod_{p\le z}\Bigl(1-\frac1p\Bigr)\ge P_0'\,\frac{\log101}{\log z}\,e^{-0.583}\ge\frac{0.1203\cdot4.5999}{1.826\,\log z}\ge\frac1{3.3\log z},
\]
using $\log101\ge\frac{73}{11}\log2\ge4.5999$ and $e^{0.583}\le1.826$. (The true constant is
$e^{\gamma_{\rm E}}\approx1.78$.) For the second claim, (b) with $\log4/\log z\le0.01$ gives
$\sum_{p\le z}\frac1p\le\log\frac{\log z}{\log2}+3.01$; use $1+2x\le e^{2x}$ and $e^{6.02}\le412$.

(d) $\sum_{k\ge1}(2k+1)^{-2}\le\frac19+\sum_{k\ge2}\frac1{4k(k+1)}=\frac19+\frac18$.
\end{proof}

\begin{lemma}[Mertens from below]\label{lem:mertlow}
\leavevmode
\begin{enumerate}[label=\textup{(\alph*)}]
\item $\sum_{p\le n}\frac{\log p}p\ge\log n-3$ for integers $n\ge1$;
\item $\sum_{p\le x}\frac1p\ge\log\log x-6.5$ for real $x\ge2$;
\item $\prod_{p\le x}(1+\frac1p)\ge e^{-7}\log x$ for real $x\ge2$.
\end{enumerate}
\end{lemma}

\begin{proof}
(a) By $n^n/n!\le e^n$ and Legendre's formula with $v_p(n!)\le n/(p-1)$,
\[
n\log n-n\le\log n!\le\sum_{p\le n}\frac{n\log p}{p-1}=n\sum_{p\le n}\frac{\log p}p+n\sum_{p\le n}\frac{\log p}{p(p-1)},
\]
and the last sum is at most $\sum_{m\ge2}\log m/(m(m-1))\le2$. (b) Let $z=\lfloor x\rfloor$. In \eqref{eq:abel} with $m=2$
use $A(n)\ge\log n-3$ from (a), and then \eqref{eq:delta}:
\[
\sum_{p\le z}\frac1p\ge1-\frac3{\log z}+\sum_{n=2}^{z-1}\Bigl(\delta_n-\frac{\delta_n^2}2\Bigr)-3\Bigl(\frac1{\log2}-\frac1{\log z}\Bigr)
\ge\log\log z-\log\log2+1-\frac3{\log2}-\frac12\sum_{n\ge2}\delta_n^2 .
\]
Here $\delta_n\le\log(n+1)/\log n-1\le1/(n\log n)$, and $\sum_{n\ge2}(n\log n)^{-2}\le0.75$, so the right side is at least
$\log\log z-3.35$. Finally $\log\log x-\log\log z\le\log(\log3/\log2)\le0.47$, since $\log x/\log z\le\log(z+1)/\log z$, so
$\sum_{p\le x}1/p\ge\log\log x-3.82$. (c) $\log(1+t)\ge t-t^2$ for $t\ge0$, (b), and
$\sum_pp^{-2}\le\frac{35}{72}$ (Lemma~\ref{lem:mertens}(d)).
\end{proof}

\section{The numerics of the cluster bound}\label{app:num}
The constants are those of Proposition~\ref{prop:condev}, and $\ell\ge24$. Put
\[
\Gamma=\ell^3-2\ell-7.
\]

\begin{lemma}[Step costs]\label{lem:stepcost}
Suppose $0<\eta\theta^{K_0}\le4$ and $0\le\mu\le2^{3\ell}$. Then for $0\le s\le d\le K_0$ and $p\in\{0,1\}$, with
$j=d-s+p$,
\[
\log_2\bigl(\eta\theta^{K_0}\,\Sigma_{s,\,j}\bigr)\le3\ell s+5\ell+7-\min\bigl(64\ell\min(2j,s),\ \Gamma\bigr).
\]
\end{lemma}

\begin{proof}
\emph{Constants.} From $L+1\le2^\ell$ we get $c_Y\le2^{2\ell+10}/Y$, so
$c_02^\Gamma\le2^{5\ell+2\ell+10+\Gamma-\ell^3}=2^{5\ell+3}$. Also $c_{\rm r}=2^{5\ell}$,
$c_1=32(L+1)\sqrt{E_0}\le2^{4\ell+5}\le2^{5\ell}$ and $c_2=(L+1)^3(1+\log E_0)/2\le2^{3\ell}\cdot3\cdot2^\ell\le2^{5\ell}$,
so $2c_{\rm r}+c_1+c_2\le2^{5\ell+3}$.

\emph{Cases.} Recall $K_0=3\ell$. If $s=0$, then the minimum is $0$ and
$\eta\theta^{K_0}\Sigma_{0,j}\le8\mu\le2^{3\ell+3}\le2^{5\ell+7}$. If $s\ge1$ and $j=0$, the minimum is again $0$ and
$\eta\theta^{K_0}\Sigma_{s,0}\le2^{3\ell s+3}$. If $s,j\ge1$, put $N_*=\min(64\ell\min(2j,s),\Gamma)$. From $N_*\le128\ell j$,
$N_*\le64\ell s$ and $Q+1\ge2^{128\ell}$ we get $2^{N_*}\le(Q+1)^j$ and $2^{N_*}\le(Q+1)^{s/2}\le(Q+1)^s$, and $N_*\le\Gamma$
gives $c_02^{N_*}\le2^{5\ell+3}$. So $2^{N_*}$ times the bracket in $\Sigma_{s,j}$ is at most
$c_02^{N_*}+2c_{\rm r}+c_1+c_2\le2^{5\ell+4}$, and $\eta\theta^{K_0}\Sigma_{s,j}\le4\cdot2\cdot2^{3\ell s}\cdot2^{5\ell+4}\cdot2^{-N_*}$.
\end{proof}

The step cost decays like $2^{-64\ell\min(2j,s)}$: through the ray term when the cofactor exponent $j$ is small,
through $D\ge(Q+1)^s$ otherwise. The decay is capped at $2^{-\Gamma}$ by the floor $c_0\approx2^{-\ell^3}$, and it
is absent when $s=0$ (free) or $j=0$ (exact); the averaging has to pay for these.

\begin{lemma}[Average over $t_e$]\label{lem:avgcost}
Let $\ell\ge24$, $1\le d\le3\ell$, $p,\varrho\in\{0,1\}$, $\varrho\le d$, and $p=1$ if $d=1$. Then
\[
\sum_{s=\varrho}^d\bigl(3\ell s+5\ell+7-\min(64\ell\min(2(d-s+p),s),\ \Gamma)\bigr)\le(d+1-\varrho)\bigl(-7-d(3\ell+6)\bigr).
\]
Consequently \eqref{eq:W} holds with $A_*=2^{-7}$ and $\lambda=1/(2k)$ whenever $0<\eta\le2$, $\theta\ge1$,
$\theta^{K_0}\le\frac98$, $0\le\mu\le2^{3\ell}$ and $2k\le2^{3\ell+6}$.
\end{lemma}

\begin{proof}
Let $S=\sum_{s=\varrho}^d\min(64\ell\min(2(d-s+p),s),\Gamma)$. Since $\sum_{s=\varrho}^ds\le d(d+1)/2$,
$d+1-\varrho\le d+1$ and $-7-d(3\ell+6)<0$, it suffices to show
\begin{equation}\label{eq:Sneed}
3\ell\,\frac{d(d+1)}2+(d+1)\bigl(5\ell+14+d(3\ell+6)\bigr)\le S.
\end{equation}
Note $64\ell\le\Gamma$.

\emph{$d\le2$.} The term $s=1$ has $\min(2(d-1+p),1)=1$ (for $d=1$ as $p=1$), so $S\ge64\ell$. The left side of
\eqref{eq:Sneed} is $19\ell+40$ for $d=1$ and $42\ell+78$ for $d=2$, both at most $64\ell$.

\emph{$d\ge3$.} For $1\le s\le d$ put $x_s=\frac{64\ell}d\,s(d-s)$. Since $s(d-s)\le d\min(2(d-s+p),s)$, we have
$x_s\le64\ell\min(2(d-s+p),s)$, and $\min(x,\Gamma)\ge x-x^2/(4\Gamma)$ for $x\ge0$ (if $x>\Gamma$ the difference is
$-(x-2\Gamma)^2/(4\Gamma)$). With $\sum_{s=0}^ds(d-s)=(d^3-d)/6$ and $\sum_{s=0}^d(s(d-s))^2=(d^5-d)/30$,
\[
S\ge\frac{64\ell}d\cdot\frac{d^3-d}6-\Bigl(\frac{64\ell}d\Bigr)^2\frac{d^5-d}{120\,\Gamma}
\ge\frac{64\ell(d^2-1)}6-\frac{(64\ell)^2d^3}{120\,\Gamma}.
\]
Since $104\ell^3\ge8400\ell+29400$, we have $4096\ell^3\le4200\Gamma$, so the last term is at most $35d^3/\ell$. Multiplying
by $\ell$, \eqref{eq:Sneed} follows from
\[
\ell\Bigl(\tfrac32\ell d(d+1)+(d+1)\bigl(5\ell+14+d(3\ell+6)\bigr)\Bigr)+35d^3\le\tfrac{32}3\ell^2(d^2-1).
\]
The left side is $\frac92\ell^2d^2+6\ell d^2+\frac{19}2\ell^2d+20\ell d+5\ell^2+14\ell+35d^3$. From
$d(d-3)(3\ell-d)\ge0$ we get $d^3\le(3+3\ell)d^2-9\ell d$, so the difference of the two sides is at least
\[
d^2\bigl(\tfrac{37}6\ell^2-111\ell-105\bigr)-d\bigl(\tfrac{19}2\ell^2-295\ell\bigr)-\tfrac{47}3\ell^2-14\ell .
\]
For $\ell\ge24$ the first bracket is positive, so with $d^2\ge3d$ this is at least
$d(9\ell^2-38\ell-315)-\frac{47}3\ell^2-14\ell$, and with $d\ge3$ at least $\frac{34}3\ell^2-128\ell-945>0$.

For the consequence: $\eta\theta^{K_0}\le\frac94\le4$, so Lemma~\ref{lem:stepcost} applies to every term, and dividing
by $d+1-\varrho$ gives the average $\le-7-d(3\ell+6)$ in $\log_2$. This is at most
$\log_22^{-7}+d\log_2(1/(2k))$, as $2k\le2^{3\ell+6}$.
\end{proof}

Summed over the positions of a vertex of each type, the decay $\gtrsim\frac{32}3\ell d^2$ beats the cost $\frac32\ell d^2$ of the factors $2^{K_0t_e}$, with enough
margin to pay $\lambda^{-d}=(2k)^d$ for the edges and the factor $4$ per vertex of Lemma~\ref{lem:shapecount}. A vertex
without a private prime is exact only in the orders that put all its neighbours before it, and only its average cost
over the $d+1-\varrho$ positions enters. The cap at $\Gamma\approx\ell^3$ costs a term of order $\ell d^3/\ell^3$, which
the margin absorbs for $d\le3\ell$.

\section{The formal proof}\label{app:lean}
Every lemma, proposition, corollary and theorem of this paper has been formalized in Lean~4 \cite{Lean4} (version
4.22.0) with Mathlib \cite{mathlib} (tag v4.22.0, commit \texttt{79e94a093aff4a60fb1b1f92d9681e407124c2ca}) and checked by Lean's kernel. The formalization is the package
\texttt{erdos-straus-34-supplement-lean}, version~2: $114$ modules with $25{,}504$ lines. Of these, $68$ build the sieve,
the polymer gas and the mass of \S\S\ref{sec:sieve}--\ref{sec:mass}; $9$ (\texttt{Uniform/}) the divisor reduction, the
factor $\lambda_K$ and the proof of Theorem~\ref{thm:A}(b); $36$ (\texttt{BoxA/}) \S\ref{sec:full} and the proof of
Theorem~\ref{thm:A}(a); and one (\texttt{Barrier}) Proposition~\ref{prop:barrier}. A root audit file
\texttt{ImprovedStatementAudit.lean} of $122$ lines restates the headline statements. The package is archived on
Zenodo, doi:\,\leanDOItwo\ (all versions: doi:\,\leanDOIall). The script \texttt{compile.sh}
installs the toolchain, fetches Mathlib, builds the package and runs every check below (about $30$--$40$ minutes on two
cores; \texttt{--lean4checker} adds the replay of every module). Lean names below omit the namespace prefixes
\lean{ErdosStraus}, \lean{ErdosStraus.Annealed}, \lean{ErdosStraus.Free}, \lean{ErdosStraus.FreeBoot},
\lean{ErdosStraus.Uniform} and \lean{ErdosStraus.BoxA}.

\emph{Headline and audit.} Theorem~\ref{thm:A}(a) is \lean{erdos_straus_EK_uniformZ} (\texttt{BoxA/StageZ.lean}) and
Theorem~\ref{thm:A}(b) is \lean{erdos_straus_EK_uniform} (\texttt{Uniform/Stage2.lean}), for the count \lean{EKu} of the
$2\le n\le N$ for which no positive integers $x,y,z$ satisfy $K/n=1/x+1/y+1/z$ in $\mathbb Q$. The audit file restates
the objects of the statements from Mathlib definitions only (\lean{MyRep}, \lean{MyEK}, \lean{MycU}, each equal to the
development's object by \texttt{rfl}), and derives \lean{main_bound_full} (Theorem~\ref{thm:A}(a)), \lean{main_bound}
(Theorem~\ref{thm:A}(b)), \lean{main_bound_full_ES} and \lean{main_bound_ES} (Corollary~\ref{cor:ES}), and
\lean{barrier_bound} (Proposition~\ref{prop:barrier}, with the events, the family of sets $T$ and $\Psi$ restated as
\lean{MyU}, \lean{MyFam}, \lean{MyQ}, and $\mathbf E[\Psi^2]$ written out as a sum over configurations). The asymptotic
improvement of the constant in Proposition~\ref{prop:barrier} is not formalized.

\emph{Checks.} All declarations ($3171$) depend only on \texttt{propext}, \texttt{Classical.choice} and
\texttt{Quot.sound}, and the source contains no \texttt{sorry}, \texttt{admit}, \texttt{native\_decide}, \texttt{axiom},
\texttt{implemented\_by}, \texttt{unsafe} or \texttt{extern}. Every module is replayed by \texttt{lean4checker}; the whole
environment of version~1, with its import closure in Lean's core library and Mathlib, was replayed by
\texttt{lean4checker --fresh}. The numerical inequalities are proved in the kernel from rational bounds on $\log2$ and
$e$, without floating-point evaluation.

\emph{Differences from the text.} In a few places the formal statements differ from the text; none changes a
conclusion.
\begin{itemize}[leftmargin=*]
\item Lemma~\ref{lem:typeII} carries the redundant hypothesis $n\ge(a+1)u_t$, so the formal sieve excludes from $E_t$
only the $n\ge(a+1)u_t$, and Proposition~\ref{prop:oneclass} has an extra term $2^A(2^L+1)2^L$.
\item For $j=1$ the formal key lemmas use the ray term $1/Q$ instead of $2/(Q+1)$ (the remark after
Lemma~\ref{lem:key}), and state (E4) twice: without the cofactor condition, with the ray term $E_0/Q$ absorbed in $c_0$
(used for $\mu_q$ in Lemma~\ref{lem:Mq}, which is formally slightly sharper), and for $j\ge1$ with ray term
$2E_0/(Q+1)^j$ (used in Lemma~\ref{lem:slot}).
\item Propositions~\ref{prop:onelevel} and~\ref{prop:sums} are formally stated for $L\ge L_0$ with $2^9M\le L_0^3$,
the only case used, rather than for $K\le L^3$.
\item The proposition of \S\ref{sec:disc} is formalized with $K\ge4$ and with an arbitrary finite set of positive triples
in place of the Type~II triples.
\item The formal proof of Proposition~\ref{prop:massF0} counts a smaller interval of scales; events and triples range over
$\mathbb N$ rather than $\mathbb Z$.
\item In \S\ref{sec:full} the formal proof uses an injective key $\kappa$ built with \lean{Nat.pair}. The package also
contains intermediate versions of the argument, with two boxes $B_a$, $B_c$ in place of $W_a$, $W_c$; the mass of one wide
scale (Proposition~\ref{prop:masscore}) is the box computation \lean{scale_massBa}, which does not use the box.
\end{itemize}

\emph{Statements and declarations.}
{\small
\begin{longtable}{@{}l>{\raggedright\arraybackslash}p{0.74\textwidth}@{}}
Theorem~\ref{thm:A}(a) & \lean{erdos_straus_EK_uniformZ}, \lean{erdos_straus_uniformZ}, \lean{erdos_straus_uniform_L0Z}, \lean{main_bound_full}\\
Theorem~\ref{thm:A}(b) & \lean{erdos_straus_EK_uniform}, \lean{erdos_straus_uniform}, \lean{final_atU}, \lean{L_small}, \lean{c093}, \lean{cU_le_cF}, \lean{cU}, \lean{main_bound}\\
Corollary~\ref{cor:ES} & \lean{main_bound_full_ES}, \lean{main_bound_ES}\\
Lemma~\ref{lem:typeII} & \lean{representable_of_event}\\
Lemma~\ref{lem:transfer} & \lean{transferE}, \lean{sumF2E_coeff}, \lean{count_class}, \lean{rsum_mul_rsum}\\
Lemma~\ref{lem:gas} & \lean{Polymer.decomposition}, \lean{EvGasPolymer.WtEv_mul}, \lean{EvGasPolymer.polymer_rep}, \lean{EvGasPolymer.genfun}, \lean{cTrE_nonneg}, \lean{cTrE_agree}\\
Lemma~\ref{lem:tail} & \lean{tail_bound}, \lean{truncated_bound}\\
Proposition~\ref{prop:cond} & \lean{sumF2E_bound}, \lean{FkE_of_noEvent}, \lean{EvReduction.gasBoundCE_of_KP}, \lean{EvCond.cond_prop_E_all}\\
Proposition~\ref{prop:dob} & \lean{KP.dobrushin}\\
Corollary~\ref{cor:KP} & \lean{KP.kp_real}\\
Lemma~\ref{lem:lattice} & \lean{Lattice.thue}, \lean{Lattice.collinear}, \lean{Lattice.offRay}, \lean{Lattice.box_count}, \lean{Lattice.lemma_73_1}\\
Lemma~\ref{lem:key} & \lean{key_lemma_cof}, \lean{key_lemma_cofJ}, \lean{KeyLemmaB.ray_unique}\\
Definition~\ref{def:condfam} & \lean{CondFam}, \lean{coprime_vPart}\\
Lemma~\ref{lem:reduce} & \lean{decompS}, \lean{E_le_reduce}, \lean{card_coprime_le_Rcls}, \lean{Rcls}, \lean{L8_rep_mul}, \lean{P0u}\\
Lemma~\ref{lem:condev} & \lean{condEv_holds_iff}, \lean{condEv_valid}\\
Lemma~\ref{lem:condkey} & \lean{A2_conditional}, \lean{A2_explicit}, \lean{A2_cofJ}\\
Lemma~\ref{lem:oneray} & \lean{L14_one_ray_triples}, \lean{E0F_small}\\
Lemma~\ref{lem:nodup} & \lean{A1_no_dup}\\
Lemma~\ref{lem:condprune} & \lean{A1_prune}, \lean{mass_prunedC}, \lean{epsPrime_F}\\
Proposition~\ref{prop:condev} & \lean{EvSite.FamOKE}, \lean{KeyJE}, \lean{FamOKJD}, \lean{famOKE_cond}, \lean{keyJE_famE}, \lean{diag_cond}\\
Lemma~\ref{lem:brun} & \lean{rough_affineFU}, \lean{V_geK}, \lean{lamK}, \lean{lamK_ge_one}, \lean{Sieve5.rB}, \lean{J0F}\\
Lemma~\ref{lem:nsq} & \lean{AffineBrun.nsq_affine}\\
Lemma~\ref{lem:omega} & \lean{AffineBrun.omega_affine}, \lean{AffineBrun.tau_sum_affine}\\
Corollary~\ref{cor:good} & \lean{good_affineFU}, \lean{good_progressionFU}, \lean{AffineBrun.progression_to_affine}\\
Lemma~\ref{lem:massprog} & \lean{mass_progressionFU}\\
Lemma~\ref{lem:wpairs} & \lean{W_weighted}, \lean{W_simple}, \lean{coprime_ac_of_class}\\
Proposition~\ref{prop:massF0} & \lean{mass_fullFU}, \lean{scale_massFU}, \lean{scale_mass_eFU}, \lean{FamFullF}, \lean{condFam_of_goodF}, \lean{vPart_ge_twoF}, \lean{scaleOK_blocksF}, \lean{condOK_fullF}, \lean{count4F}, \lean{J0F_le}\\
Lemma~\ref{lem:gain} & \lean{DivSumPhi.divsum_phi_gain_divisors}, \lean{DivSumGain.trunc_divsum}\\
Proposition~\ref{prop:massF} & \lean{mass_full_prunedFU}, \lean{mass_ge_XFK}, \lean{eset_gainFU}, \lean{prodKpU_le_lamK}, \lean{gain_aux}\\
Corollary~\ref{cor:site} & \lean{MultiB}, \lean{site_absF}\\
Lemma~\ref{lem:Mq} & \lean{EvSums.Mq_le}, \lean{EvSums.Mq_uniform}\\
Proposition~\ref{prop:polymer} & \lean{kp_condF}, \lean{firstOrder_boundF}, \lean{r2_boundF}, \lean{EvReduction.famZ_le_exp}\\
Lemma~\ref{lem:cof} & \lean{cof_of_real}, \lean{pow_le_DD}\\
Lemma~\ref{lem:slot} & \lean{slot_boundD}\\
Proposition~\ref{prop:fubini} & \lean{fubini_shapeD}, \lean{fubini_abs}, \lean{gC_le_prod}, \lean{Shape}, \lean{Realizes}, \lean{te}, \lean{tl}\\
Lemma~\ref{lem:admshape} & \lean{conn_sum_le_shapes}, \lean{sum_reals_relabel}\\
Lemma~\ref{lem:rank} & \lean{sum_rk}, \lean{rk_swap}\\
Lemma~\ref{lem:avgorder} & \lean{exists_perm_le}\\
Lemma~\ref{lem:shapecount} & \lean{sum_shape_le}, \lean{card_edges_le}, \lean{one_add_pow_pairs_le}\\
Proposition~\ref{prop:sitefree} & \lean{site_sum_pow_freeD}, \lean{site_sum_multi_freeD}, \lean{shape_boundD}\\
Lemma~\ref{lem:avgcostmain} & \lean{hW_free24}, \lean{GF}\\
Proposition~\ref{prop:gasfree} & \lean{multiB_cond24}, \lean{epsF_le24}, \lean{kp_free_cond24}, \lean{r2_free_cond24}, \lean{KP.kp_real}\\
Proposition~\ref{prop:oneclass} & \lean{class_boundF24K}, \lean{window_FK}, \lean{condOK_fullF_P0u}, \lean{ClassBound.greedy}, \lean{ClassBound.card_le_base}\\
Proposition~\ref{prop:onelevel} & \lean{rcls_bound}, \lean{rcls_le}, \lean{card_units_eq}, \lean{P0u_le_pow}, \lean{K_le_of}\\
Lemma~\ref{lem:chooseL} & \lean{final_atU} (the choice of $L$), \lean{EthU_le}, \lean{TcF_bounds}, \lean{hypUK}, \lean{clsK}, \lean{ell1024}\\
Lemma~\ref{lem:exponent} & \lean{core_realM}, \lean{lN_le_real}, \lean{cF}\\
Lemma~\ref{lem:euler} & \lean{sum_le_euler}, \lean{sum_inv_le_euler}, \lean{sum_isqrt_le}, \lean{totient_prod_primes}, \lean{phi_mul_euler}\\
Proposition~\ref{prop:sums} & \lean{levelU}, \lean{num_T3}, \lean{four_pow_SPU_le}, \lean{card_SPU_le}\\
Lemma~\ref{lem:L0K} & \lean{cubeK}, \lean{cubeK_spec}, \lean{cubeK_min}, \lean{L0K}, \lean{L0K_props}\\
Lemma~\ref{lem:keyfull} & \lean{key_lemma_fullB}, \lean{strip_sum}, \lean{fibre_full}, \lean{AdmFull}\\
Definition~\ref{def:fullfam} & \lean{CondFamFull}, \lean{FullAll}\\
Lemma~\ref{lem:condkeyfull} & \lean{A2Full_condB}, \lean{A2Full_explicit}, \lean{A2Full_cofJ}\\
Lemma~\ref{lem:diagfull} & \lean{diag_Full}, \lean{hyperbola_card}, \lean{diag_condF}\\
Lemma~\ref{lem:prunefull} & \lean{prunedC_byF}, \lean{mass_prunedCF}, \lean{epsFP}, \lean{epsFP_F}\\
Proposition~\ref{prop:evfull} & \lean{CondOKF}, \lean{famOKE_condF}, \lean{keyJE_famEF}, \lean{condEv_injF}, \lean{c0FJ}, \lean{c2FE}, \lean{BoxA.crF}, \lean{deltaF}\\
Lemma~\ref{lem:dedup} & \lean{Dedup}, \lean{dedup_inj}, \lean{prio_sub_dedup}, \lean{keyT}, \lean{EqvT}\\
Definition~\ref{def:core} & \lean{WaF}, \lean{WcF}, \lean{CoreZ}, \lean{FamFullF}, \lean{disj_F0_Wa}, \lean{disj_F0Wa_Wc}\\
Lemma~\ref{lem:coreinj} & \lean{inj_WW}, \lean{nW_lt_v}, \lean{noDup_WF0}, \lean{noDup_WF0_arith}, \lean{noDup_WWc}, \lean{noDup_WWc_arith}, \lean{inj_BaBa_arith}, \lean{coreZ_inj}, \lean{A1_no_dupX}\\
Lemma~\ref{lem:countW} & \lean{count_W}, \lean{blockW_card_ge}, \lean{sum_sq_teleW}, \lean{blocksCW}, \lean{blocksAW}\\
Lemma~\ref{lem:pairsbox} & \lean{pairsXC_inv_sq}, \lean{pairsXC_card_ge}, \lean{pairsC_card_ge}, \lean{pairsC_split}, \lean{badXC_card_le}\\
Proposition~\ref{prop:masscore} & \lean{mass_FamWaS}, \lean{scaleOK_W}, \lean{scale_massBa}, \lean{mass_progBa}, \lean{dclass_sum_ge}, \lean{FamWaS_subW}, \lean{FamWcS_subW}, \lean{class_swap}, \lean{mass_coreZ_raw}, \lean{union_algZ}\\
Lemma~\ref{lem:numfull} & \lean{c0FJ_small24}, \lean{CstarF_le24}, \lean{costF_log_le2X}, \lean{arith_main24X}, \lean{poly24X}, \lean{hW_gen}, \lean{hW_free24F}, \lean{epsF_le24F}, \lean{kp_free_cond24F}, \lean{r2_free_cond24F}\\
Proposition~\ref{prop:massfull} & \lean{DedZ}, \lean{condOK_dedZ}, \lean{mass_dedZ}, \lean{mass_algZ}, \lean{mass_ge_XFKZ}, \lean{exp_neg7_ge}\\
\S\ref{sec:proofA}, one level & \lean{class_bound_genF}, \lean{window_genF}, \lean{class_boundF24KZ}, \lean{TcF_le32}, \lean{TcF_boundsZ}, \lean{hypUKZ}\\
Lemma~\ref{lem:L0Z} & \lean{cubeZ}, \lean{L0Z}, \lean{L0Z_props}\\
Proposition~\ref{prop:barrier} & \lean{Barrier.barrier}, \lean{Barrier.christoffel}, \lean{Barrier.E_Q_sq}, \lean{Barrier.rankin}, \lean{Barrier.typeII_sum_le}, \lean{Barrier.harmonic_sigma_le}, \lean{barrier_bound}\\
Lemma~\ref{lem:mertens} & \lean{Mertens.theta_le}, \lean{Mertens.mertens1_upper}, \lean{Mertens2.recip_sum_le}, \lean{Mertens2.recip_sum_le'}, \lean{Mertens5.V_lower5}, \lean{Sieve5.V_ge5}, \lean{MassB.prod_two_le2}\\
Lemma~\ref{lem:mertlow} & \lean{MertensUpper.mertens1_lower}, \lean{DivSumGain.recip_sum_ge}, \lean{DivSumGain.prod_one_add_all}\\
Lemma~\ref{lem:stepcost} & \lean{costF_log_le2}, \lean{c0F_small24}, \lean{c12_le24}, \lean{Cstar2_le24}\\
Lemma~\ref{lem:avgcost} & \lean{arith_main24}, \lean{poly24}, \lean{hW_free24}\\
\end{longtable}}

\subsection*{Data availability}
{\sloppy The Lean~4 formalization of every statement of this paper (version~2 of the package) is the file
\texttt{erdos-straus-34-supplement-lean-v2.zip}, SHA-256 checksum \texttt{\leanSHA}, archived on Zenodo,
doi:\,\leanDOItwo. Version~1, for Theorem~\ref{thm:A}(b), is doi:\,\leanDOI; all versions are cited as
doi:\,\leanDOIall. The script \texttt{compile.sh} installs the toolchain,
builds the package and runs the axiom audits and the kernel replays.\par}

\subsection*{Funding and competing interests}
This research received no specific funding. The author declares no competing interests.

\end{document}